\RequirePackage{fix-cm}
\documentclass[smallextended]{svjour3}

\usepackage[final]{graphicx}
\usepackage{booktabs}
\usepackage{multirow}
\usepackage{amsmath,amssymb,amsfonts}
\usepackage{mathrsfs}
\usepackage[title]{appendix}
\usepackage{xcolor}
\usepackage{manyfoot}
\usepackage{algorithm}
\usepackage{algorithmicx}
\usepackage{algpseudocode}
\usepackage[justification=centering]{caption}
\usepackage{hyperref}
\hypersetup{hidelinks}
\usepackage{geometry}
\usepackage{nicefrac}
\usepackage[final]{showlabels}
\usepackage{float}
\usepackage{enumitem}
\usepackage{soul}
\usepackage{tikz}
\usepackage{pgfplots}
\pgfplotsset{compat=1.18}
\usepgfplotslibrary{groupplots}
\usepackage{placeins}
\usetikzlibrary{positioning,arrows.meta,fit,calc,backgrounds}

\smartqed

\let\oldproof\proof
\let\oldendproof\endproof

\renewenvironment{proof}
  {\oldproof}
  {\qed\oldendproof}
\newcommand{\Sn}{\mathbb{S}^{n}}
\newcommand{\Sp}{\mathbb{S}^{n}_{+}}
\newcommand{\Sm}{\mathbb{S}^{n}_{-}}

\newcommand{\On}{\mathbb{O}^{n}}

\newcommand{\tr}{\operatorname{tr}}
\newcommand{\Diag}{\operatorname{Diag}}
\newcommand{\diag}{\operatorname{diag}}
\newcommand{\dist}{\operatorname{dist}}
\newcommand{\sym}{\operatorname{sym}}
\newcommand{\grad}{\operatorname{grad}}

\newcommand{\argmin}{\operatorname*{arg\,min}}

\newcommand{\ip}[2]{\left\langle #1,#2\right\rangle}
\newcommand{\norm}[1]{\left\lVert #1\right\rVert}
\newcommand{\Nlim}{N^{\operatorname{lim}}}
\newcommand{\Nhat}{\widehat N}
\newcommand{\cmin}{\underline c}
\newcommand{\cmax}{\overline c}

\begin{document}

\title{An Inexact Augmented Lagrangian Method for Semidefinite Complementarity Optimization with Applications}

\author{
    Wanqi Hou \and
    Yanzhen Peng \and
    Chao Ding
}

\institute{
    Wanqi Hou \at
    Academy of Mathematics and Systems Science,
    Chinese Academy of Sciences,
    Beijing, China \\
    School of Advanced Interdisciplinary Sciences,
    University of Chinese Academy of Sciences,
    Beijing, China \\
    \email{houwanqi19@mails.ucas.ac.cn}
    \and
    Yanzhen Peng \at
    Academy of Mathematics and Systems Science,
    Chinese Academy of Sciences,
    Beijing, China \\
    School of Advanced Interdisciplinary Sciences,
    University of Chinese Academy of Sciences,
    Beijing, China \\
    \email{pengyanzhen23@mails.ucas.ac.cn}
    \and
    Chao Ding \at
    Academy of Mathematics and Systems Science,
    Chinese Academy of Sciences,
    Beijing, China \\
    \email{dingchao@amss.ac.cn}
}
\date{This version: November 20, 2026}

\maketitle

\begin{abstract}

Semidefinite complementarity constraints arise in applications such as inverse semidefinite programming and rank-constrained matrix optimization. We develop a structure-preserving augmented Lagrangian method (SPALM) for mathematical programs with these constraints. The method augments the smooth constraints and retains matrix complementarity in its subproblems. We handle these subproblems by an inexact proximal alternating framework with verifiable descent and stationarity conditions. With vanishing absolute subproblem errors, feasible accumulation points of infinite outer sequences are asymptotically M-stationary. Under AM-regularity, these points are M-stationary; we prove that active-set SDCMPCC-LICQ is sufficient for this regularity. For each AL subproblem, the full stationarity residual converges to zero along infinite accepted inner sequences if the initial sublevel set is bounded, the proximal parameters are bounded above and away from zero, and the block errors vanish. If the block errors are also square-summable, an $O(T^{-1/2})$ best-residual bound holds over $T$ accepted iterations. Huber--Stiefel and fixed-basis procedures generate feasible matrix candidates, whose acceptance is verified by descent tests and computable limiting-normal certificates. For rank-constrained correlation matrix approximation, we construct a lifting that preserves the projected feasible set and global optimal value and satisfies active-set SDCMPCC-LICQ at every feasible point. The stationarity guarantees apply to this lifted model. Finite-budget experiments on this lifting and a single-equality observation-fitting inverse SDP meet the prescribed numerical feasibility and stationarity tests in all $120$ runs.

    \vspace{10pt}
    \noindent \textbf{Key words:} semidefinite cone complementarity constraints; augmented Lagrangian method; asymptotic stationarity; proximal alternating method

\end{abstract}

\section{Introduction}

Mathematical programs with equilibrium constraints (MPECs) arise in equilibrium models, variational inequalities, and bilevel programming, with applications in engineering and economics. Mathematical programs with complementarity constraints (MPCCs) form a basic subclass \cite{LuoPangRalph1996,ScheelScholtes2000}. Complementarity constraints also enter optimality-based reformulations of constrained bilevel problems. Bilevel models appear, for example, in hyperparameter optimization and meta-learning \cite{FranceschiEtAl2018}, and in domain reweighting for large language model pretraining \cite{FanPagliardiniJaggi2024}. These connections motivate the study of complementarity-constrained optimization. Although optimality theory and numerical methods for vector MPCCs are well developed, optimization with matrix-valued complementarity constraints has received less attention.

Matrix complementarity arises naturally in inverse optimization, where observed decisions are used to infer parameters of a forward optimization model \cite{ChanMahmoodZhu2025}. For a linear semidefinite forward problem, the Karush--Kuhn--Tucker (KKT) conditions consist of primal--dual feasibility and complementarity between a positive semidefinite primal matrix $X$ and a positive semidefinite dual slack $S$:
\begin{equation*}
X\succeq0,\qquad S\succeq0,\qquad XS=0.
\end{equation*}
Imposing these KKT conditions in an inverse model introduces a semidefinite complementarity constraint; see, for example, \cite{WuZhangZhangLu2016}. 

Ghatrani and Ghate~\cite{GhatraniGhate2024} study the recovery of constraint matrices in inverse SDPs and include truss topology design as an application. In robotics, semidefinite relaxations have been developed for trajectory planning \cite{kang2026fast}, and inverse optimality conditions have been used to learn costs from manipulation demonstrations \cite{englert2017inverse}. Family~$\mathrm{II}$ in Section~\ref{subsec:family-2-prime} considers cost calibration with fixed constraint data and jointly fits a cost matrix and a primal optimizer to observations.

Low-rank correlation matrix approximation provides another application. In financial model calibration, one seeks a correlation matrix that approximates a given target under a prescribed rank bound \cite{Wu2003,BrigoMercurio2002,LiQi2011}. The approximating matrix must be positive semidefinite and have unit diagonal. An auxiliary matrix allows the rank restriction to be expressed through semidefinite complementarity and additional constraints \cite{DingSunYe2014}; Section~\ref{subsec:family-1-prime} develops the lifting used here. Further applications, including certain robust optimization models, can be found in \cite{DingSunYe2014}.

Motivated by these applications, we develop a numerical method for mathematical programs with semidefinite cone complementarity constraints (SDCMPCCs) of the form
\begin{equation}\label{P}
     \min_{w} f(w) \quad \text{s.t.} \quad g(w) \leq 0,\ h(w) = 0,\ \Sp \ni G(w) \perp H(w) \in \Sm,
\end{equation}
Here $w \in \mathcal{Z}$, where $\mathcal{Z}$ is a finite-dimensional real Euclidean space. The mappings $f: \mathcal{Z} \to \mathbb{R}$, $g: \mathcal{Z} \to \mathbb{R}^p$, $h: \mathcal{Z} \to \mathbb{R}^q$, and $G, H: \mathcal{Z} \to \Sn$ are continuously differentiable. SDCMPCCs are matrix analogues of MPCCs and include them as a special case. The complementarity constraint is described by the set
\begin{equation*}
\Omega:= \operatorname{gph}N_{\Sp}=\bigl\{(Y,Z)\in\Sp\times\Sm:\ YZ=0\bigr\}.
\end{equation*}
Thus the complementarity condition in \eqref{P} is $(G(w),H(w))\in\Omega$. Even for vector MPCCs, complementarity causes standard constraint qualifications to fail, and classical KKT conditions need not hold at local minimizers \cite{LuoPangRalph1996,ScheelScholtes2000}. SDCMPCCs combine this degeneracy with the nonpolyhedral geometry of $\Sp$ and the nonconvexity of $\Omega$. Stationarity analysis and algorithm design therefore require the geometry of the semidefinite complementarity set.

On the theoretical side, Ding et al.~\cite{DingSunYe2014} showed that the Robinson constraint qualification fails in the direct conic formulation of SDCMPCCs. Their proximal and limiting normal-cone formulas for $\Omega$ give explicit spectral characterizations of S- and M-stationarity, while a projection reformulation yields C-stationarity conditions. Wu et al.~\cite{WuZhangZhang2014} established a tangent-cone formula and characterized B-stationarity. They also studied its relationships with S-, M-, and C-stationarity under SDCMPCC-tailored constraint qualifications, together with related second-order sufficient conditions. Madariaga and Ram{\'i}rez~\cite{MadariagaRamirez2025} study stationarity for general cone complementarity constraints from a facial viewpoint. These results provide the geometric basis for the stationarity conditions used in SDCMPCC algorithms.

Algorithmic research on vector MPCCs includes regularization and relaxation methods \cite{Scholtes2001,KadraniDussaultBenkachroun2009,KanzowSchwartz2013}, penalty and exact-penalty methods \cite{RalphWright2004,HuRalph2004,LinFukushima2003}, and augmented Lagrangian methods \cite{IzmailovSolodovUskov2012,AndreaniSecchinSilva2019,GuoDeng2022}. In particular, Guo and Deng~\cite{GuoDeng2022} retain vector complementarity explicitly in the augmented Lagrangian subproblems and solve them by a nonmonotone projected gradient method. This inner method uses an explicit projection onto the vector complementarity set, which separates into independent scalar problems.

For the matrix case, Yan and Fukushima~\cite{YanFukushima2011} proposed a smoothing method for symmetric-cone complementarity constraints, including the semidefinite case, based on a nonsmooth projection reformulation. Under suitable assumptions, limits of stationary points of the smoothed problems are C-stationary. For general SDCMPCCs, Wu and Zhang~\cite{WuZhang2015} studied a bilinear penalty framework and proved S-stationarity of feasible accumulation points under their assumptions. The resulting subproblems are nonlinear semidefinite programs, but their numerical solution is outside the scope of that analysis. For an inverse linear SDP, Wu et al.~\cite{WuZhangZhangLu2016} developed a sequential convex programming method that exploits the DC structure of a nonsmooth partial penalty reformulation.

More recently, Santos et al.~\cite{SantosCostaKrulikovskiGeremia2026} analyze an augmented Lagrangian method based on a slack reformulation that retains semidefinite complementarity explicitly. They impose an approximate KKT condition on the separate semidefinite and scalar complementarity constraints. Their analysis gives AW-stationarity at feasible accumulation points and AC-stationarity under a multiplier-dependent accuracy requirement; a subsequent C-stationarity result uses an SDCMPCC-Robinson condition. These stationarity notions are defined in their formulation.

For second-order cone complementarity, Liang et al.~\cite{LiangZhuZhangLinZhu2026} study an augmented Lagrangian method for a direct cone-and-scalar formulation. Its K-stationarity target is the KKT condition of that formulation. Ye and Zhou~\cite[Proposition~7.1 and Corollary~7.4]{YeZhou2016} show that K-stationarity implies S-stationarity, while the converse can fail when a cone has dimension at least three. Thus a finite-multiplier target may impose an additional representation requirement; our M-stationarity guarantee is not stronger than K-stationarity. Conic subproblems and positive-parameter relaxations are different objects from the unrelaxed scalar-orthogonality formulation.

For the semidefinite problem, the choice of a subproblem stationarity criterion depends on how complementarity is represented. Replacing $\Omega$ with two semidefinite constraints and one scalar orthogonality constraint preserves the feasible set, but the resulting finite multipliers need not represent all regular normals to $\Omega$. This representation also has a nonzero abnormal multiplier at every feasible point, so convergence results requiring Robinson's constraint qualification cannot be invoked directly. Section~\ref{subsec:spectral_geometry} derives these properties, and Example~\ref{ex:scalar_normal_gap} exhibits a missing coupled normal even at a strictly complementary point. We focus on M-stationarity, which requires the complementarity multiplier pair to lie in the limiting normal cone to $\Omega$. Outer semicontinuity of this mapping preserves the normal-cone relation along convergent sequences of slack pairs and their associated normals \cite[equation~(2.1)]{JiaKanzowMehlitzWachsmuth2022}.

For approaches that retain complementarity explicitly, the matrix projection is a central computational difficulty. The general augmented Lagrangian framework of Jia et al.~\cite{JiaKanzowMehlitzWachsmuth2022} retains a closed set and penalizes the remaining smooth constraints. Its projected-gradient inner method requires projections onto the retained set. For semidefinite complementarity, a fixed common orthogonal basis reduces the projection problem to scalar problems with explicit solutions. Finding an optimal basis remains, in general, a coupled nonconvex optimization problem. Thus the inexpensive projection available for vector MPCCs does not directly extend to SDCMPCCs.

Building on this retained-set framework and the asymptotic stationarity analysis in \cite{AndreaniHaeserSecchinSilva2019,Mehlitz2020Asymptotic,JiaKanzowMehlitzWachsmuth2022}, we develop SPALM for SDCMPCCs. Its slack candidates satisfy complementarity, and the outer method controls the smooth linking errors. Inexact proximal alternation produces candidates that must pass both descent and same-point limiting-normal tests. The stationarity target includes the coupled normals omitted by the scalar representation. Vanishing absolute errors yield AM-stationarity at feasible accumulation points, and AM-regularity permits passage to M-stationarity without assuming bounded raw multipliers for the degenerate scalar reformulation. Generating and certifying matrix candidates remains a computational cost. Section~\ref{sec:ALM} compares the resulting accuracy requirement with the separate-constraint tests of \cite{SantosCostaKrulikovskiGeremia2026}; no general cost or runtime ordering follows from the different interfaces.

The main contributions are as follows.
\begin{enumerate}[label=(\roman*)]
\item
We specialize the asymptotic stationarity framework to SPALM and establish AM-stationarity of feasible accumulation points of infinite outer sequences when the absolute subproblem errors vanish (Theorem~\ref{thm:outer_AM_stationarity}). We characterize active-set SDCMPCC-LICQ and prove that it implies AM-regularity, which yields M-stationarity. Feasibility of accumulation points is treated separately from this stationarity implication.

\item
To solve the augmented Lagrangian subproblem at each outer iteration, we develop an inexact proximal alternating framework with verifiable block descent and stationarity conditions. Boundedness of the initial sublevel set, two-sided bounds on the proximal parameters, and vanishing block errors yield convergence of the full residual to zero along infinite accepted sequences (Theorem~\ref{thm:main}). Square-summable block errors give an $O(T^{-1/2})$ best-residual bound over $T$ accepted iterations (Corollary~\ref{cor:inner_complexity}).

\item
We construct Huber--Stiefel and fixed-basis procedures for generating feasible slack candidates, together with computable limiting-normal certificates for their verification. Acceptance requires both true descent and the prescribed normal-residual bound. We also derive spectral objective-gap bounds for the fixed-basis candidates.

\item
We construct a compact rank lifting that preserves the projected feasible set and global optimal value of the rank-constrained correlation problem and satisfies active-set SDCMPCC-LICQ at every feasible point. The M-stationarity guarantees apply to this lifted formulation (Remark~\ref{rem:familyIprime_lifting_scope}). For this lifting and a single-equality observation-fitting inverse SDP, we verify active-set SDCMPCC-LICQ and boundedness of the augmented Lagrangian sublevel sets under the stated model assumptions. For each fixed instance, comparing the inner initial point with an explicit feasible reference point gives a common upper bound on the outer AL values. With completed inner solves and vanishing absolute errors, this yields bounded primal iterates and M-stationary accumulation points (Lemma~\ref{lem:reference_AL_bound} and Corollary~\ref{cor:numerical_outer_boundedness}).
\end{enumerate}

The numerical study evaluates SPALM on $24$ parameter groups with five seeds per group. All $120$ runs satisfy the prescribed numerical success tests, combining absolute feasibility and stationarity tests with a normalized witness-based residual target. Retained joint spectral data and independent numerical verification support the residual evaluations. The tables and figures summarize final residuals and solver wall-clock times. The finite-precision interpretation is given in Remark~\ref{rem:implementation_theory}. In the inverse SDP experiments, we calibrate the stiffness matrix of a composite plate using a nominal stiffness matrix and a buckling mode computed on an independent finite-element mesh.

The remainder of the paper follows the analysis and computation outlined above. Section~\ref{sec:pre} develops the slack reformulation and stationarity framework. Section~\ref{sec:ALM} presents SPALM and its outer convergence results. Section~\ref{sec:inner_solver} analyzes accepted proximal alternating sequences, and Section~\ref{sec:slack_solver} develops slack-candidate procedures and their verification. Section~\ref{sec:numerical_experiments} presents the models and numerical results. The final section concludes the paper.

\paragraph{Notation.}
We denote by $\mathbb S^n$ the space of $n\times n$ real symmetric matrices, by $\mathbb S_+^n$ and $\mathbb S_-^n$ the cones of positive- and negative-semidefinite matrices, and by $\mathbb O^n$ the set of $n\times n$ orthogonal matrices.
For $A,B\in\mathbb{R}^{n\times n}$, the Frobenius inner product is defined by
$\langle A,B\rangle := \operatorname{tr}(A^{\top} B)$, and the corresponding Frobenius norm is
$\|A\|_F := \sqrt{\langle A,A\rangle}$. The notation $\|\cdot\|$ denotes the norm induced by the given inner product on a Euclidean space; in particular, it is the Euclidean norm for vectors. For explicit matrix arguments, we write $\|\cdot\|_F$ for the Frobenius norm and $\|\cdot\|_2$ for the spectral norm.
Product spaces are equipped with the canonical product norms induced by their component norms. In particular, for $u=(Y,Z)\in\Sn\times\Sn$, $\|u\|^2=\|Y\|_F^2+\|Z\|_F^2$.
For a differentiable mapping $F$, $F'(w)$ and $F'(w)^*$
denote its derivative and adjoint derivative, respectively.
For a scalar $s\in\mathbb{R}$, define
$s_+ := \max\{s,0\}\ \text{and} \ s_- := \max\{-s,0\}$. These operations are applied componentwise to vectors.
Moreover, $\operatorname{Diag}(x)$ denotes the diagonal matrix generated by $x$, whereas $\operatorname{diag}(A)$ denotes the vector of diagonal entries of $A$.
The eigenvalues $\lambda_i(A)$ follow the order in the stated spectral decomposition; $\lambda_j^\downarrow(A)$ denotes the $j$-th largest eigenvalue. Any spectral reordering also reorders the corresponding eigenvector columns. For a nonempty closed convex set $D$, $\Pi_D$ denotes its Euclidean projection, and
\begin{equation*}
N_D(x):=\{v:\langle v,z-x\rangle\leq0\ \text{for every }z\in D\}
\quad(x\in D),\qquad N_D(x):=\varnothing\quad(x\notin D).
\end{equation*}
Regular and limiting normals to closed nonconvex sets are defined below.

\section{Problem reformulation and stationarity framework}\label{sec:pre}
We separate the smooth constraints from semidefinite complementarity, establish M-stationarity equivalence under the unique slack lift, and introduce the asymptotic stationarity framework used by SPALM. Active-set SDCMPCC-LICQ and the absence of singular M-multipliers provide sufficient conditions for the required AM-regularity.

\subsection{Slack reformulation and normal-cone geometry}
\label{subsec:spectral_geometry}
To fit problem \eqref{P} into the structured geometric constraint framework of \cite{JiaKanzowMehlitzWachsmuth2022}, we introduce $u=(Y,Z)\in\mathcal U:=\Sn\times\Sn$ and $x=(w,u)\in\mathcal X:=\mathcal Z\times\mathcal U$. Problem \eqref{P} can be equivalently reformulated as
\begin{equation}\label{RE}
\min_{w,Y,Z} f(w)\quad\text{s.t.}\quad
g(w)\leq0,\quad h(w)=0,\quad G(w)=Y,\quad H(w)=Z,\quad (Y,Z)\in\Omega.
\end{equation}
Set $\varphi(x):=f(w)$ and
\begin{equation*}
\Psi(x):=(g(w),h(w),G(w)-Y,H(w)-Z),\qquad
C:=\mathbb R_-^p\times\{0\}^q\times\{0\}^{\Sn}\times\{0\}^{\Sn}.
\end{equation*}
Thus \eqref{RE} minimizes $\varphi(x)$ subject to $\Psi(x)\in C$ and $x\in\mathcal Z\times\Omega$. 
Write
\begin{equation*}
\mathcal F:=\{w:g(w)\leq0,\ h(w)=0,\ (G(w),H(w))\in\Omega\},\qquad
\widetilde{\mathcal F}:=\{x\in\mathcal Z\times\Omega:\Psi(x)\in C\}.
\end{equation*}
The map $w\mapsto(w,(G(w),H(w)))$ is a bijection from $\mathcal F$ to $\widetilde{\mathcal F}$ and preserves objective values. In particular, the optimal values of the two formulations coincide.
For any $u=(Y,Z)\in\Omega$, the linking error also bounds the distance of the matrix image to complementarity:
\begin{equation}\label{eq:link_geometric_distance}
\operatorname{dist}((G(w),H(w)),\Omega)
\leq\bigl(\|G(w)-Y\|_F^2+\|H(w)-Z\|_F^2\bigr)^{1/2}.
\end{equation}

We recall the following standard notions from variational analysis \cite{RockafellarWets1998}. For a closed subset $S$ of a finite-dimensional Euclidean space, $x^k\to_S\bar x$ means $x^k\in S$ and $x^k\to\bar x$. The outer limit is defined by
\[
\operatorname*{Limsup}_{x\to\bar x}\mathcal T(x):=\{v:\exists\,x^k\to\bar x,\ v^k\to v,\ v^k\in\mathcal T(x^k)\text{ for all }k\}.
\]
\begin{definition}[Regular and limiting normal cones]\label{def:normals}
For $\bar x\in S$, define
\begin{equation*}
\Nhat_S(\bar x):=\left\{v:\limsup_{\substack{x\to_S\bar x\\x\ne\bar x}}
\frac{\ip{v}{x-\bar x}}{\|x-\bar x\|}\leq0\right\},
\qquad
\Nlim_S(\bar x):=\operatorname*{Limsup}_{x\to_S\bar x}\Nhat_S(x).
\end{equation*}
Both cones are empty outside $S$. We have $\Nhat_S\subseteq\Nlim_S$, and the limiting normal-cone mapping $\Nlim_S$ is outer semicontinuous.
\end{definition}
For $u=(Y,Z)\in\Omega$, consider the common spectral representation \cite[Sections~2.3 and~3]{DingSunYe2014}:
\begin{equation}\label{decomposition}
\mathcal E(u):=Y+Z=P\operatorname{Diag}(\Lambda_{\alpha\alpha},0,\Lambda_{\gamma\gamma})P^\top,
\qquad \Lambda_{\alpha\alpha}\succ0,\quad\Lambda_{\gamma\gamma}\prec0.
\end{equation}
Here $\alpha,\beta,\gamma$ are the positive, zero, and negative spectral index sets; $P=(P_\alpha,P_\beta,P_\gamma)\in\On$,
$Y=P_\alpha\Lambda_{\alpha\alpha}P_\alpha^\top$, and $Z=P_\gamma\Lambda_{\gamma\gamma}P_\gamma^\top$.
\begin{equation}\label{eq:Sigma_def}
\Sigma_{ij}:=\frac{\lambda_i(\mathcal E(u))}{\lambda_i(\mathcal E(u))-\lambda_j(\mathcal E(u))}
\quad(i\in\alpha,\ j\in\gamma),\qquad
\widetilde\Gamma^Y:=P^\top\Gamma^YP,\quad\widetilde\Gamma^Z:=P^\top\Gamma^ZP.
\end{equation}
Let $E_{IJ}$ denote the all-ones matrix on $I\times J$,
and let $\circ$ denote the Hadamard product.
The common multiplier space is
\[
\begin{aligned}
\mathcal K_W(u):=
\Bigl\{(\Gamma^Y,\Gamma^Z)\in\Sn\times\Sn:\;&
\widetilde\Gamma^Y_{\alpha\alpha}=
\widetilde\Gamma^Y_{\alpha\beta}=
\widetilde\Gamma^Y_{\beta\alpha}=0,\quad 
\widetilde\Gamma^Z_{\gamma\gamma}=
\widetilde\Gamma^Z_{\beta\gamma}=
\widetilde\Gamma^Z_{\gamma\beta}=0,\\
&
\Sigma_{\alpha\gamma}\circ\widetilde\Gamma^Y_{\alpha\gamma}
+
(E_{\alpha\gamma}-\Sigma_{\alpha\gamma})
\circ\widetilde\Gamma^Z_{\alpha\gamma}=0
\Bigr\}.
\end{aligned}
\]

Multiplication in the basis $P$ gives the basis-independent identity
\begin{equation}\label{eq:KW_algebraic}
(\Gamma^Y,\Gamma^Z)\in\mathcal K_W(u)
\quad\Longleftrightarrow\quad Y\Gamma^Y-\Gamma^ZZ=0.
\end{equation}
Indeed, the $\alpha\alpha$ and $\alpha\beta$ blocks of $\widetilde\Gamma^Y$ and the $\beta\gamma$ and $\gamma\gamma$ blocks of $\widetilde\Gamma^Z$ vanish by invertibility of the nonzero spectral blocks. For $i\in\alpha$ and $j\in\gamma$, the remaining relation is $\lambda_i(\mathcal E(u))\widetilde\Gamma^Y_{ij}-\lambda_j(\mathcal E(u))\widetilde\Gamma^Z_{ij}=0$, which gives the $\alpha\gamma$ condition by \eqref{eq:Sigma_def}. A simultaneous rotation of the positive and negative eigenspaces gives a geometric interpretation of this coupling in dimension two.

\begin{example}[Coupled normals at a strictly complementary pair]\label{ex:scalar_normal_gap}
Let $n=2$ and set
\[
Y_0=\begin{pmatrix}1&0\\0&0\end{pmatrix},\qquad
Z_0=-\begin{pmatrix}0&0\\0&1\end{pmatrix},\qquad
W=\begin{pmatrix}0&1\\1&0\end{pmatrix}.
\]
The pair $u_0=(Y_0,Z_0)$ is strictly complementary since $Y_0-Z_0=I_2$. With $p_\theta=(\cos\theta,\sin\theta)^\top$ and $q_\theta=(-\sin\theta,\cos\theta)^\top$, every feasible pair near $u_0$ has the form $(a p_\theta p_\theta^\top,-b q_\theta q_\theta^\top)$ for $a,b>0$. Differentiating this local parameterization at $(a,b,\theta)=(1,1,0)$ gives
\[
T_\Omega(u_0)=\operatorname{span}\{(Y_0,0),(0,Z_0),(W,W)\}.
\]
Orthogonality to the rotation direction $(W,W)$ imposes $\Gamma^Y_{12}+\Gamma^Z_{12}=0$ on each normal pair. Taking the orthogonal complement of the tangent space of this local manifold gives
\begin{equation}\label{eq:coupled_psd_normal}
\Nhat_\Omega(u_0)=\Nlim_\Omega(u_0)
=\left\{\left(
\begin{pmatrix}0&r\\r&s\end{pmatrix},
\begin{pmatrix}t&-r\\-r&0\end{pmatrix}
\right):r,s,t\in\mathbb R\right\}.
\end{equation}

To compute the finite multiplier image at the same pair, write the scalar representation in its general form $\Omega=F^{-1}(D)$, where $F(Y,Z)=(Y,Z,\langle Y,Z\rangle)$ and $D=\Sp\times\Sm\times\mathbb R_+$. On $\Sp\times\Sm$, the inner product $\langle Y,Z\rangle$ is nonpositive and vanishes exactly when $YZ=0$.\footnote{At every feasible $u=(Y,Z)$, the tuple $(Z,Y,-1)$ satisfies \begin{equation}\label{eq:scalar_abnormal_multiplier}(Z,Y,-1)\in N_D(F(u))\setminus\{0\},\qquad F'(u)^*(Z,Y,-1)=0.\end{equation} After composition with $(G,H)$, this direction certifies the failure of Robinson's constraint qualification for the scalar representation of \eqref{P} \cite[Proposition~4.1]{DingSunYe2014}. Adding its nonnegative multiples to any KKT multiplier gives an unbounded ray of multipliers.} At $u=(Y,Z)\in\Omega$, the finite multiplier image is
\begin{equation}\label{eq:scalar_normal_image}
\begin{aligned}
\mathcal N^{\mathrm{img}}(Y,Z)
&:=F'(u)^*N_D(F(u))\\
&=\bigl(N_{\Sp}(Y)\times N_{\Sm}(Z)\bigr)+\mathbb R_-(Z,Y)
\subseteq\Nhat_\Omega(u).
\end{aligned}
\end{equation}
For $u+(\Delta Y,\Delta Z)\in\Omega$, orthogonality at both feasible pairs gives $\langle Z,\Delta Y\rangle+\langle Y,\Delta Z\rangle=-\langle\Delta Y,\Delta Z\rangle=O(\|(\Delta Y,\Delta Z)\|^2)$. Together with the two cone-normal inequalities, this proves the inclusion in \eqref{eq:scalar_normal_image}.

At $u_0$, the two cone normals are $N_{\Sp}(Y_0)=\mathbb R_+Z_0$ and $N_{\Sm}(Z_0)=\mathbb R_+Y_0$. Hence the finite multiplier image consists of the diagonal pairs
\[
\mathcal N^{\mathrm{img}}(Y_0,Z_0)=\{(sZ_0,tY_0):s,t\in\mathbb R\}.
\]
The scalar equality $\langle Y,Z\rangle=0$ gives the same image at $u_0$. The off-diagonal pair $\xi_0=(W,-W)$ belongs to the normal space in \eqref{eq:coupled_psd_normal} and is orthogonal to this multiplier image. Since $\|\xi_0\|=2$, the two distances are
\begin{equation}\label{eq:normal_residual_gap}
\operatorname{dist}\bigl(\xi_0,\Nlim_\Omega(u_0)\bigr)=0,\qquad
\operatorname{dist}\bigl(\xi_0,\mathcal N^{\mathrm{img}}(Y_0,Z_0)\bigr)=2.
\end{equation}

A path with derivative $\xi_0$ illustrates the local scaling of the scalar orthogonality residual. Consider the positive and negative semidefinite matrices
\[
Y(t)=\begin{pmatrix}1&t\\t&t^2\end{pmatrix},\qquad
Z(t)=-\begin{pmatrix}t^2&t\\t&1\end{pmatrix}.
\]
Along this path, $|\langle Y(t),Z(t)\rangle|=4t^2$. Since $(Y'(0),Z'(0))=\xi_0$ is normal to the local manifold at $u_0$ and has norm $2$, the distance to $\Omega$ satisfies $\operatorname{dist}((Y(t),Z(t)),\Omega)=2|t|+o(|t|)$ as $t\to0$.
\end{example}

In Example~\ref{ex:scalar_normal_gap}, the common kernel is trivial and the $\alpha\gamma$ coupling is already present. When $Y$ and $Z$ have a nontrivial common kernel, the conditions on the $\beta\beta$ blocks distinguish the following spectral multiplier sets \cite{DingSunYe2014}.

\begin{definition}[Spectral multiplier sets]\label{def:spectral_sets}
For $u\in\Omega$, define
\[
\begin{aligned}
\mathcal K_S(u)&:=\{(\Gamma^Y,\Gamma^Z)\in\mathcal K_W(u):
\widetilde\Gamma^Y_{\beta\beta}\preceq0,\ \widetilde\Gamma^Z_{\beta\beta}\succeq0\},\\
\mathcal K_M(u)&:=\Nlim_\Omega(u),\\
\mathcal K_C(u)&:=\{(\Gamma^Y,\Gamma^Z)\in\mathcal K_W(u):
\langle\widetilde\Gamma^Y_{\beta\beta},\widetilde\Gamma^Z_{\beta\beta}\rangle\leq0\}.
\end{aligned}
\]
Appendix~\ref{app:KM-auxiliary} recalls the full spectral characterization of $\mathcal K_M(u)$.
\end{definition}

These sets satisfy
\[
\mathcal K_S(u)=\Nhat_\Omega(u)\subseteq\mathcal K_M(u)\subseteq\mathcal K_C(u)\subseteq\mathcal K_W(u),\qquad
\operatorname{span}\Nhat_\Omega(u)=\operatorname{span}\Nlim_\Omega(u)=\mathcal K_W(u).
\]
The span identities follow from the inclusion chain and the fact that every symmetric $\beta\beta$ block is a difference of two positive semidefinite matrices. The remaining block relations are homogeneous and independent of the $\beta\beta$ blocks.

\subsection{Stationarity and asymptotic stationarity}\label{subsec:stationarity}
Using the spectral multiplier sets introduced above and the stationarity framework of \cite{DingSunYe2014}, we define the S-, M-, C-, and W-stationarity conditions and describe their relationships.
\begin{definition}[M-stationarity for the SDCMPCC]\label{def:stationarity_P}
A feasible point $\bar w\in\mathcal F$, with $\bar u=(G(\bar w),H(\bar w))$, is M-stationary if there exist $\lambda\in\mathbb R^p$, $\mu\in\mathbb R^q$, and $\Gamma^Y,\Gamma^Z\in\Sn$ such that
\begin{equation}\label{eq:stationarity_P}
0=\nabla f(\bar w)+g'(\bar w)^*\lambda+h'(\bar w)^*\mu
+G'(\bar w)^*\Gamma^Y+H'(\bar w)^*\Gamma^Z,
\end{equation}
\begin{equation}\label{eq:ineq_multiplier_P}
\lambda\geq0,\qquad\langle\lambda,g(\bar w)\rangle=0,
\end{equation}
and
\begin{equation}\label{eq:spectral_multiplier_P}
(\Gamma^Y,\Gamma^Z)\in\mathcal K_M(\bar u).
\end{equation}
\end{definition}
For \eqref{RE}, M-stationarity at $(\bar w,\bar u)\in\widetilde{\mathcal F}$ involves additional linking multipliers $M^G,M^H\in\Sn$ satisfying
\begin{equation}\label{eq:RE_stat_w}
0=\nabla f(\bar w)+g'(\bar w)^*\lambda+h'(\bar w)^*\mu
+G'(\bar w)^*M^G+H'(\bar w)^*M^H,
\end{equation}
\begin{equation}\label{eq:RE_stat_u}
M^G=\Gamma^Y,\qquad M^H=\Gamma^Z,
\end{equation}
together with \eqref{eq:ineq_multiplier_P}--\eqref{eq:spectral_multiplier_P}.

The following equivalence allows stationarity results for the slack reformulation to be transferred to the original problem.
\begin{proposition}[Equivalence of M-stationarity]\label{prop:Mstat_equiv}
For $\bar w\in\mathcal F$ and $\bar u=(G(\bar w),H(\bar w))$, $\bar w$ is M-stationary for \eqref{P} if and only if $(\bar w,\bar u)$ is M-stationary for \eqref{RE}.
\end{proposition}
\begin{proof}
Eliminating $M^G,M^H$ from \eqref{eq:RE_stat_w}--\eqref{eq:RE_stat_u} yields \eqref{eq:stationarity_P}. Conversely, setting $M^G=\Gamma^Y$ and $M^H=\Gamma^Z$ gives the M-stationarity conditions for \eqref{RE}.
\end{proof}
Replacing $\mathcal K_M$ in \eqref{eq:spectral_multiplier_P}
by $\mathcal K_S$, $\mathcal K_C$, or $\mathcal K_W$
defines S-, C-, or W-stationarity, respectively.
The inclusions among these multiplier sets yield
\[
\text{S-stationarity}
\ \Longrightarrow\
\text{M-stationarity}
\ \Longrightarrow\
\text{C-stationarity}
\ \Longrightarrow\
\text{W-stationarity}
\]

We next specialize the asymptotic stationarity and regularity framework of \cite[Section~2.2]{JiaKanzowMehlitzWachsmuth2022} to the slack reformulation \eqref{RE}. Define the multiplier image
\begin{equation}\label{def:multiplier_image}
\mathcal M(w,u,z):=\Psi'(w,u)^*N_C(\Psi(w,u)-z)
+(\{0\}\times\Nlim_\Omega(u)),
\end{equation}
where $w\in\mathcal Z$, $u\in\mathcal U$, and $z\in\mathbb R^p\times\mathbb R^q\times\Sn\times\Sn$.
Since $\nabla\varphi(w,u)=(\nabla f(w),0,0)$ and $\Nlim_{\mathcal Z\times\Omega}(w,u)=\{0\}\times\Nlim_\Omega(u)$, exact M-stationarity is
\begin{equation*}
0\in\nabla\varphi(\bar w,\bar u)+\mathcal M(\bar w,\bar u,0).
\end{equation*}

\begin{definition}[AM-stationarity and AM-regularity]\label{def:AM_stationarity}\label{def:AM_regular}
Let $(\bar w,\bar u)\in\widetilde{\mathcal F}$.
\begin{enumerate}[label=(\roman*)]
\item The point is asymptotically M-stationary (AM-stationary) if there exist $(w^k,u^k)\to(\bar w,\bar u)$, $z^k\to0$, and $e^k\to0$ with
$e^k\in\nabla\varphi(w^k,u^k)+\mathcal M(w^k,u^k,z^k)$ for every $k$.
\item Problem \eqref{RE} is AM-regular at this point if
\begin{equation*}
\operatorname*{Limsup}_{(w,u,z)\to(\bar w,\bar u,0)}\mathcal M(w,u,z)
\subseteq\mathcal M(\bar w,\bar u,0).
\end{equation*}
\end{enumerate}
\end{definition}
Condition (i) requires $u^k\in\Omega$ and $\Psi(w^k,u^k)-z^k\in C$, while $(w^k,u^k)$ need not be feasible for \eqref{RE}. AM-regularity ensures the outer semicontinuity of the multiplier image, under which AM-stationarity implies M-stationarity.

\subsection{Active-set SDCMPCC-LICQ}\label{sec:active_LICQ_AM}
To establish AM-regularity for the slack reformulation, we adapt the spectral SDCMPCC-LICQ condition in \cite[Theorem~6.1(ii)]{DingSunYe2014} to the active inequality constraints and derive an equivalent local uniform estimate.

For $\Upsilon=(\lambda,\mu,M^G,M^H)$ and $\Gamma=(\Gamma^Y,\Gamma^Z)$, define
\begin{equation}\label{eq:full_linearized_operator}
\mathcal J_w(\Upsilon,\Gamma):=\Psi'(w,u)^*\Upsilon+(0,\Gamma)
=\begin{pmatrix}
g'(w)^*\lambda+h'(w)^*\mu+G'(w)^*M^G+H'(w)^*M^H\\
-M^G+\Gamma^Y\\-M^H+\Gamma^Z
\end{pmatrix}.
\end{equation}
Call $(w,u,z)$ admissible if $u\in\Omega$ and $s:=\Psi(w,u)-z\in C$. For such a triple, set $\mathcal I(s_g):=\{i:(s_g)_i=0\}$ and
\begin{equation*}
\mathfrak L(w,u,z):=\operatorname{span}N_C(s)\times\mathcal K_W(u).
\end{equation*}

\begin{lemma}[Closedness of the linearized multiplier space]\label{lem:L_equiv}
For an admissible triple with $u=(Y,Z)$,
\begin{equation}\label{eq:L_algebraic}
\mathfrak L(w,u,z)=\{(\lambda,\mu,M^G,M^H,\Gamma^Y,\Gamma^Z):
\operatorname{Diag}(s_g)\lambda=0,\quad Y\Gamma^Y-\Gamma^ZZ=0\}.
\end{equation}
The mapping $\mathfrak L$ has a closed graph on admissible triples.
\end{lemma}
\begin{proof}
The activity condition is exactly $\operatorname{Diag}(s_g)\lambda=0$, and \eqref{eq:KW_algebraic} gives the spectral condition. These equalities are preserved under limits, as is admissibility by continuity of $\Psi$ and closedness of $C$ and $\Omega$.
\end{proof}

\begin{definition}[Active-set SDCMPCC-LICQ]\label{def:SDCMPCC_LICQ}
Active-set SDCMPCC-LICQ holds at $(\bar w,\bar u)\in\widetilde{\mathcal F}$ if
$\ker\mathcal J_{\bar w}\cap\mathfrak L(\bar w,\bar u,0)=\{0\}$.
\end{definition}

\begin{proposition}[Reduced and uniform characterizations of LICQ]
\label{prop:reduced_LICQ}
For $(\bar w,\bar u)\in\widetilde{\mathcal F}$, with $\bar u=(\bar Y,\bar Z)$, the following are equivalent:
\begin{enumerate}[label=(\roman*)]
\item Active-set SDCMPCC-LICQ holds at $(\bar w,\bar u)$.
\item The only tuple $(\lambda,\mu,\Gamma^Y,\Gamma^Z)$ satisfying
\begin{equation*}
\begin{aligned}
g'(\bar w)^*\lambda+h'(\bar w)^*\mu+G'(\bar w)^*\Gamma^Y+H'(\bar w)^*\Gamma^Z&=0,\\
\lambda_i=0\quad(i\notin\mathcal I(g(\bar w))),\qquad
\bar Y\Gamma^Y-\Gamma^Z\bar Z&=0
\end{aligned}
\end{equation*}
is zero. 
\item There exist a neighborhood $\mathcal N$ of $(\bar w,\bar u,0)$ and $c>0$ such that
\begin{equation}\label{eq:uniform_LICQ}
\|\mathcal J_w t\|\geq c\|t\|
\quad\text{for every admissible }(w,u,z)\in\mathcal N,\quad t\in\mathfrak L(w,u,z).
\end{equation}
\end{enumerate}
\end{proposition}
\begin{proof}
The last two equations of $\mathcal J_{\bar w}t=0$ give $M^G=\Gamma^Y$ and $M^H=\Gamma^Z$; elimination proves (i)$\Leftrightarrow$(ii). If (iii) fails, choose admissible triples within $1/k$ of $(\bar w,\bar u,0)$ and nonzero tuples with $\|\mathcal J_{w^k}t^k\|<\|t^k\|/k$. Normalize them and pass to a convergent subsequence on the unit sphere. Lemma~\ref{lem:L_equiv} and continuity of $\mathcal J_w$ give a unit tuple in $\ker\mathcal J_{\bar w}\cap\mathfrak L(\bar w,\bar u,0)$, contradicting (i). Conversely, (iii) applied at the base point implies (i).
\end{proof}

\begin{proposition}
[AM-regularity under active-set SDCMPCC-LICQ]\label{thm:LICQ_AMregular}
If active-set SDCMPCC-LICQ holds at $(\bar w,\bar u)\in\widetilde{\mathcal F}$, then \eqref{RE} is AM-regular there.
\end{proposition}
\begin{proof}
Let admissible triples $(w^k,u^k,z^k)\to(\bar w,\bar u,0)$ and $v^k\in\mathcal M(w^k,u^k,z^k)$ satisfy $v^k\to\bar v$. Choose representations $v^k=\mathcal J_{w^k}t^k$, where
\begin{equation*}
t^k=(\Upsilon^k,\Gamma^k)\in
N_C(\Psi(w^k,u^k)-z^k)\times\Nlim_\Omega(u^k)\subseteq\mathfrak L(w^k,u^k,z^k).
\end{equation*}
The uniform estimate \eqref{eq:uniform_LICQ} gives $c\|t^k\|\leq\|v^k\|$ for all sufficiently large \(k\). Thus, after passing to a subsequence, \(t^k\to(\bar\Upsilon,\bar\Gamma)\). By the closed-graph properties of the normal-cone mappings, $(\bar\Upsilon,\bar\Gamma)\in N_C(\Psi(\bar w,\bar u))\times N_\Omega^{\lim}(\bar u)$. Continuity then gives $\bar v=\mathcal J_{\bar w}(\bar\Upsilon,\bar\Gamma)\in\mathcal M(\bar w,\bar u,0)$.
\end{proof}

\begin{remark}[NNAMCQ and AM-regularity]
A weaker sufficient condition for AM-regularity is the absence of singular M-multipliers. Following \cite[Definition~6.1]{DingSunYe2014}, a nonzero tuple $(\lambda,\mu,\Gamma^Y,\Gamma^Z)$ is called a \emph{singular M-multiplier} for \eqref{P} at $\bar w\in\mathcal F$ if
\[
\begin{aligned}
g'(\bar w)^*\lambda+h'(\bar w)^*\mu +G'(\bar w)^*\Gamma^Y+H'(\bar w)^*\Gamma^Z&=0,\\
\lambda\in N_{\mathbb R_-^p}(g(\bar w)),\qquad (\Gamma^Y,\Gamma^Z)&\in\Nlim_\Omega(\bar u),
\end{aligned}
\]
where $\bar u=(G(\bar w),H(\bar w))$. The absence of such multipliers is equivalent to the no nonzero abnormal multiplier constraint qualification (NNAMCQ) for \eqref{RE}, since the linking multipliers satisfy $M^G=\Gamma^Y$ and $M^H=\Gamma^Z$.

Since
\[
N_C(\Psi(\bar w,\bar u)) \times\Nlim_\Omega(\bar u)
\subseteq\mathfrak L(\bar w,\bar u,0),
\]
active-set SDCMPCC-LICQ implies NNAMCQ, which in turn implies AM-regularity by \cite[Section~2.2, Eq.~(2.2)] {JiaKanzowMehlitzWachsmuth2022}.

The first implication is strict. Let $n=1$, $G(w)=w$, and $H(w)=-w$, with no other constraints. At the unique feasible point $w=0$, the homogeneous stationarity equation is $\Gamma^Y=\Gamma^Z$. While $\mathcal K_W(0,0)=\mathbb R^2$ admits the nonzero solution $(1,1)$, we have
\[ \{(\Gamma^Y,\Gamma^Z)\in\Nlim_\Omega(0,0): \Gamma^Y=\Gamma^Z\}=\{(0,0)\}.
\]
Hence, active-set SDCMPCC-LICQ fails, whereas NNAMCQ holds.
\end{remark}

\begin{remark}[Scalar complementarity specialization]
For $\Omega_{\mathrm{cc}}=\operatorname{gph}N_{\mathbb R_+^m}$, write $y=G(w)$ and $\zeta=-H(w)$.
The usual nonnegative pair $(y,\zeta)$ defines
$I^{+0},I^{0+},I^{00}$, and
\begin{equation*}
\mathcal K_W^{\mathrm{cc}}(y,-\zeta)
=\{(\Gamma^Y,\Gamma^Z):
\Gamma_i^Y=0\ (i\in I^{+0}),\
\Gamma_i^Z=0\ (i\in I^{0+})\}.
\end{equation*}
In this scalar setting, the active-set SDCMPCC-LICQ
reduces to the classical MPCC-LICQ
\cite[Definition~2]{AnitescuTsengWright2007},
with $\Gamma^Z$ corresponding to $H=-\zeta$.
However, diagonal-valued maps into the full matrix cone
need not yield the same LICQ, since off-diagonal
multipliers remain admissible.
\end{remark}

\section{Structure-preserving augmented Lagrangian method}\label{sec:ALM}

SPALM penalizes the smooth constraints of \eqref{RE} and retains $\Omega$ explicitly. We separate feasibility of accumulation points from their AM-stationarity; AM-regularity then gives M-stationarity.

\subsection{Augmented Lagrangian function}\label{sec:outer_AL}

Following the partial-penalization strategy in \cite{GuoDeng2022,JiaKanzowMehlitzWachsmuth2022}, we incorporate $\Psi(x)\in C$ into the objective and retain $x\in\mathcal Z\times\Omega$.

For $\rho>0$ and
$\Upsilon=(\lambda,\mu,M^G,M^H)$, define the augmented Lagrangian function by
\begin{equation*}
  \mathcal L_\rho(x;\Upsilon)
  :=
  \varphi(x)
  +\frac{\rho}{2}
  \operatorname{dist}^2
  \left(
    \Psi(x)+\frac{\Upsilon}{\rho},C
  \right).\end{equation*}
Using the product structure of $C$, we obtain the following expression in $x=(w,Y,Z)$:
\begin{equation}\label{eq:AL}
\begin{aligned}
  \mathcal L_\rho(w,Y,Z;\lambda,\mu,M^G,M^H)
  ={}&
  f(w)
  +\frac{\rho}{2}
  \left\|
    \left(g(w)+\frac{\lambda}{\rho}\right)_+
  \right\|^2
  +\frac{\rho}{2}
  \left\|
    h(w)+\frac{\mu}{\rho}
  \right\|^2
  \\
  &+
  \frac{\rho}{2}
  \left\|
    G(w)-Y+\frac{M^G}{\rho}
  \right\|_F^2
  +
  \frac{\rho}{2}
  \left\|
    H(w)-Z+\frac{M^H}{\rho}
  \right\|_F^2.
\end{aligned}
\end{equation}

At the $k$-th outer iteration, let $\bar\Upsilon^k := (\bar\lambda^k,\bar\mu^k,\bar M^{G,k},\bar M^{H,k})$ be the safeguarded multiplier vector, and define
\begin{equation*}
  L_k(x):=\mathcal L_{\rho_k}(x;\bar\Upsilon^k).
\end{equation*}
The corresponding augmented Lagrangian subproblem is
\begin{equation}\label{eq:outer_subproblem}
  \min_{x\in\mathcal Z\times\Omega}\ L_k(x).
\end{equation}
Step~3 below requires approximate M-stationarity for this subproblem.

\subsection{SPALM for SDCMPCCs}\label{sec:outer_algorithm}

The penalty update and stopping test use a common feasibility measure.  For $x=(w,u)$, with $u=(Y,Z)$, $\lambda\in\mathbb R^p$, and $\rho>0$, define
\begin{equation}\label{eq:feasibility_measure}
  \mathcal V_\rho(x;\lambda):=\max\Bigg\{
  \left\| g(w) - \Pi_{\mathbb R_-^p}
    \left( g(w)+\frac{\lambda}{\rho}\right)
  \right\|,
  \|h(w)\|,\|G(w)-Y\|_F, \|H(w)-Z\|_F
  \Bigg\}.
\end{equation}
This projection-based measure bounds the perturbation in the smooth constraints and in the normal-cone relation for the updated multipliers, as shown below. The first term vanishes when the inequality constraints and the corresponding multiplier conditions are satisfied.

\paragraph{Feasible-reference comparison.}
When a fixed feasible point $x^{\mathrm{ref}}=(w^{\mathrm{ref}},u^{\mathrm{ref}})$ of \eqref{RE} is available, we use the following optional rule in Step~3:
\begin{equation}\label{eq:reference_initialization}
x^{k,0}\in\argmin_{x\in\{x^k,x^{\mathrm{ref}}\}}L_k(x),
\qquad L_k(x^{k+1})\leq L_k(x^{k,0}).
\end{equation}

The reference point remains fixed across outer iterations, and both AL values are evaluated with the current parameters $(\rho_k,\bar\Upsilon^k)$. Accepted inner descent ensures the second inequality in \eqref{eq:reference_initialization}.

\begin{algorithm}[H]
\caption{Structure-preserving augmented Lagrangian method (SPALM)}
\label{alg:sALcc}
\footnotesize

\begin{algorithmic}
\algrenewcommand\algorithmicrequire{\textbf{Input:}}
\algrenewcommand\algorithmicensure{\textbf{Output:}}

\Require
An initial point $x^0=(w^0,u^0)\in\mathcal Z\times\Omega$, with $u^0=(Y^0,Z^0)$;
initial multipliers $\Upsilon^0=(\lambda^0,\mu^0,M^{G,0},M^{H,0})$;
a penalty parameter $\rho_0>0$;
parameters $\beta_\rho>1$ and $\eta\in(0,1)$;
safeguards $\lambda_{\max}>0$, $\mu_{\min}<\mu_{\max}$, $M_{\max}^G>0$, and $M_{\max}^H>0$;
a feasibility tolerance $\varepsilon_{\mathrm{feas}}>0$;
a stationarity tolerance $\varepsilon_{\mathrm{stat}}>0$;
and a nonincreasing positive sequence $\{\nu_k\}$ satisfying $\nu_k\to0$.

\Ensure
A computed iterate $x^k=(w^k,u^k)$, with
$u^k=(Y^k,Z^k)$, its updated multipliers $\Upsilon^k$, and a normal certificate for the slack stationarity test.

\State Set $k=0$ and
$V_0:=\mathcal V_{\rho_0}(x^0;0)$.

\Statex \textbf{Step 1: feasibility and stationarity termination test.}
\State If $k\geq1$ and
\begin{equation*}
  V_k\leq\varepsilon_{\mathrm{feas}},\qquad \left\|\nabla_wL_{k-1}(x^k)\right\|
  \leq\varepsilon_{\mathrm{stat}},
 \qquad \operatorname{dist}\left(
    -\nabla_uL_{k-1}(x^k),
    N_\Omega^{\lim}(u^k)
  \right)
  \leq\varepsilon_{\mathrm{stat}},
\end{equation*}
\State output $x^k$ and $\Upsilon^k$, retain a normal certificate for the slack test, and terminate.
\Statex \textbf{Step 2: safeguard the multipliers.}
\State Project the vector multipliers onto their prescribed boxes and the matrix multipliers onto Frobenius balls:
\begin{equation*}
  \bar\lambda^k := \Pi_{[0,\lambda_{\max}]^p}(\lambda^k),
  \qquad
  \bar\mu^k:=\Pi_{[\mu_{\min},\mu_{\max}]^q}(\mu^k).
\end{equation*}
\begin{equation*}
  \bar M^{G,k}:=\frac{M^{G,k}}{\max\{1,\|M^{G,k}\|_F/M^G_{\max}\}},
  \qquad
  \bar M^{H,k}:=\frac{M^{H,k}}{\max\{1,\|M^{H,k}\|_F/M^H_{\max}\}}.
\end{equation*}
\State Set
$\bar\Upsilon^k
  :=
  (\bar\lambda^k,\bar\mu^k,\bar M^{G,k},\bar M^{H,k})$.

\Statex \textbf{Step 3: approximately solve the AL subproblem.}
\State If $x^{\mathrm{ref}}$ is supplied, use \eqref{eq:reference_initialization}.
\State Compute $x^{k+1}=(w^{k+1},u^{k+1})\in\mathcal Z\times\Omega$, with
$u^{k+1}=(Y^{k+1},Z^{k+1})$ such that

\begin{equation*}
  \left\|\nabla_wL_k(x^{k+1})\right\|
  \leq\nu_{k+1},\qquad
  \operatorname{dist}\left(
    -\nabla_uL_k(x^{k+1}),
    N_\Omega^{\lim}(u^{k+1})
  \right)
  \leq\nu_{k+1}.
\end{equation*}

\Statex \textbf{Step 4: update the multipliers.}
\State Set
\begin{equation*}
  \lambda^{k+1}
  :=
  \left(
    \bar\lambda^k+\rho_k g(w^{k+1})
  \right)_+,
  \qquad
  \mu^{k+1}
  :=
  \bar\mu^k+\rho_k h(w^{k+1}),
\end{equation*}
\begin{equation*}
  M^{G,k+1}
  :=
  \bar M^{G,k}
  +\rho_k\bigl(G(w^{k+1})-Y^{k+1}\bigr),
  \qquad
  M^{H,k+1}
  :=
  \bar M^{H,k}
  +\rho_k\bigl(H(w^{k+1})-Z^{k+1}\bigr).
\end{equation*}

\Statex \textbf{Step 5: compute the new feasibility measure.}
\State Set
\begin{equation*}
  V_{k+1}
  :=
  \mathcal V_{\rho_k}(x^{k+1};\bar\lambda^k).
\end{equation*}

\Statex \textbf{Step 6: update the penalty parameter.}
\State Set
\begin{equation*}
  \rho_{k+1}
  :=
  \begin{cases}
    \rho_k,
      &\text{if }V_{k+1}\leq\eta V_k,\\[1mm]
    \beta_\rho\rho_k,
      &\text{otherwise}.
  \end{cases}
\end{equation*}

\State Set $k\leftarrow k+1$ and return to Step~1.
\end{algorithmic}
\end{algorithm}
\noindent\textit{Multiplier relations.} 
The safeguards give
$M_{\mathrm{saf}}:=\sup_k\|\bar\Upsilon^k\|<\infty$,
while the updated multipliers $\Upsilon^{k+1}$ may be unbounded. To express the multiplier update and the subproblem gradient, define
\begin{equation}\label{eq:outer_perturbation}
z^{k+1}:=\Psi(x^{k+1})
-\Pi_C\left(\Psi(x^{k+1})+\frac{\bar\Upsilon^k}{\rho_k}\right).
\end{equation}
Projection optimality and the multiplier update give

\begin{equation}\label{eq:outer_multiplier_normal}
\Upsilon^{k+1}=\bar\Upsilon^k+\rho_k z^{k+1}
\in N_C\bigl(\Psi(x^{k+1})-z^{k+1}\bigr).
\end{equation}
Differentiating the squared-distance expression for $L_k$ gives
\begin{equation}\label{eq:outer_gradient_identity}
\nabla L_k(x^{k+1})
=\nabla\varphi(x^{k+1})+\Psi'(x^{k+1})^*\Upsilon^{k+1}.
\end{equation}
In particular, $\nabla_uL_k(x^{k+1})=-(M^{G,k+1},M^{H,k+1})$, so the slack test in Step~3 measures the distance of the updated linking multipliers to $\Nlim_\Omega(u^{k+1})$.

\noindent\textit{Perturbation bounds.}
The four blocks of $z^{k+1}$ are the residuals in \eqref{eq:feasibility_measure}; hence

\begin{equation}\label{eq:outer_residual_bounds}
V_{k+1}\leq\|z^{k+1}\|\leq2V_{k+1},\qquad
\operatorname{dist}\bigl(\Psi(x^{k+1}),C\bigr)\leq\|z^{k+1}\|.
\end{equation}
Since projection onto $C$ does not increase distances,
we also obtain

\begin{equation}\label{eq:outer_shift_bound}
\|z^{k+1}\|
\leq\operatorname{dist}\bigl(\Psi(x^{k+1}),C\bigr)
+\frac{\|\bar\Upsilon^k\|}{\rho_k}
\leq\operatorname{dist}\bigl(\Psi(x^{k+1}),C\bigr)
+\frac{M_{\mathrm{saf}}}{\rho_k}.
\end{equation}

\noindent\textit{Termination and error bounds.}
Suppose that Step~1 terminates at $x^k$.
Choose $\Gamma^k\in\Nlim_\Omega(u^k)$ satisfying
$\|\nabla_uL_{k-1}(x^k)+\Gamma^k\|
\leq\varepsilon_{\mathrm{stat}}$;
such a normal exists because $\Nlim_\Omega(u^k)$ is nonempty
and closed.
A retained inner normal with residual at most $\nu_k$ can be
reused when $\nu_k\leq\varepsilon_{\mathrm{stat}}$; cf.\
\cite[equations~(4.4)--(4.5)]{JiaKanzowMehlitzWachsmuth2022}.
With $z^k$ defined by \eqref{eq:outer_perturbation} and
$e^k:=\nabla L_{k-1}(x^k)+(0,\Gamma^k)$,
the projection and gradient identities, together with the
stopping test, give
\begin{equation}\label{eq:outer_finite_certificate}
\begin{gathered}
x^k\in\mathcal Z\times\Omega,\qquad
\Psi(x^k)-z^k\in C,\\
\Upsilon^k\in N_C\bigl(\Psi(x^k)-z^k\bigr),\qquad
\Gamma^k\in\Nlim_\Omega(u^k),\\
e^k=\nabla\varphi(x^k)
+\Psi'(x^k)^*\Upsilon^k+(0,\Gamma^k),\\
\|z^k\|\leq2\varepsilon_{\mathrm{feas}},\qquad
\|e^k\|\leq\sqrt2\,\varepsilon_{\mathrm{stat}}.
\end{gathered}
\end{equation}
Thus the returned point satisfies the perturbed M-stationarity
system \eqref{eq:outer_finite_certificate} with the stated
residual bounds.
For the inequality block $z_g^k$, the normal relation yields
$\lambda_i^k g_i(w^k)=\lambda_i^k z_{g,i}^k$
for $i=1,\ldots,p$, so bounds on the unperturbed complementarity
products also depend on the magnitudes of the updated multipliers.

\subsection{Global convergence analysis}

We first establish sufficient conditions for feasibility of accumulation points and then analyze their AM- and M-stationarity. For the asymptotic analysis, we consider Algorithm~\ref{alg:sALcc} with Step~1 omitted and assume that every Step~3 solve is completed. 
For the slack reformulation \eqref{RE}, the objective is $\varphi(x)=f(w)$, the smooth constraint mapping is $\Psi$, and the retained set is $\mathcal Z\times\Omega$.

Jia et al.~\cite[Proposition~4.1]{JiaKanzowMehlitzWachsmuth2022} gave sufficient conditions for feasibility of accumulation points. Our feasibility measure \eqref{eq:feasibility_measure} uses the maximum of the four residual block norms. In view of the norm equivalence in \eqref{eq:outer_residual_bounds}, their result yields the following proposition.
\begin{proposition}[Feasibility of accumulation points]\label{prop:feasibility}
Suppose that Algorithm~\ref{alg:sALcc} generates an infinite sequence $\{x^{k+1}\}$ with Step 1 omitted. If either of the following conditions holds, then every accumulation point of this sequence is feasible for \eqref{RE}:
\begin{enumerate}[label=(\roman*)]
  \item The sequence $\{\rho_k\}$ is bounded.
  \item There exists $\bar L\in\mathbb R$ such that
  \begin{equation*}
    L_k(x^{k+1})\leq\bar L\qquad\text{for all }k.
  \end{equation*}
\end{enumerate}
\end{proposition}

The next lemma ensures condition \emph{(ii)} by comparison with a fixed feasible reference point.
\begin{lemma}[Uniform AL bound from a feasible reference point]
\label{lem:reference_AL_bound}
Let $x^{\mathrm{ref}}$ be feasible for \eqref{RE}. Suppose that $\rho_k\geq\rho_0>0$, $\|\bar\Upsilon^k\|\leq M_{\mathrm{saf}}$, and the returned points satisfy \eqref{eq:reference_initialization}. Then
\begin{equation}\label{eq:reference_AL_bound}
L_k(x^{k+1})\leq L_k(x^{\mathrm{ref}})
\leq f(w^{\mathrm{ref}})+\frac{M_{\mathrm{saf}}^2}{2\rho_0}
=:\bar L_{\mathrm{ref}}\qquad\text{for every }k.
\end{equation}
\end{lemma}
\begin{proof}
Feasibility gives $\Psi(x^{\mathrm{ref}})\in C$, so
\begin{equation*}
\operatorname{dist}\left(\Psi(x^{\mathrm{ref}})+\frac{\bar\Upsilon^k}{\rho_k},C\right)
\leq\frac{\|\bar\Upsilon^k\|}{\rho_k}.
\end{equation*}
Substitution into the squared-distance definition of $L_k$ gives $L_k(x^{\mathrm{ref}})\leq f(w^{\mathrm{ref}})+\|\bar\Upsilon^k\|^2/(2\rho_k)$. The initialization comparison and descent in \eqref{eq:reference_initialization} give $L_k(x^{k+1})\leq L_k(x^{\mathrm{ref}})$. The bounds on $\rho_k$ and $\bar\Upsilon^k$ yield \eqref{eq:reference_AL_bound}.
\end{proof}

Jia et al.~\cite[Theorem~4.3 and Corollary~4.4]{JiaKanzowMehlitzWachsmuth2022} established AM-stationarity of feasible accumulation points and M-stationarity under AM-regularity. In our algorithm, the two stationarity tests in Step~3 yield an overall stationarity error bounded by $\sqrt{2}\,\nu_{k+1}$. Applying their results, together with Proposition~\ref{prop:feasibility} and the regularity and stationarity equivalence results in Section~\ref{sec:pre}, gives the following theorem.

\begin{theorem}[AM- and M-stationarity of accumulation points]\label{thm:outer_AM_stationarity}
Consider an infinite sequence $\{x^{k+1}\}$ generated by Algorithm~\ref{alg:sALcc} with Step~1 omitted. Suppose that every Step~3 solve is completed and that $\nu_k\to0$. Then the following statements hold:
\begin{enumerate}[label=(\roman*)]
  \item Every feasible accumulation point of this sequence is AM-stationary for \eqref{RE}.
  \item Furthermore, if either condition \emph{(i)} or \emph{(ii)} in Proposition~\ref{prop:feasibility} holds, then every accumulation point of this sequence is feasible and AM-stationary for \eqref{RE}.
\end{enumerate}
For any feasible accumulation point $\bar x=(\bar w,\bar u)$, if \eqref{RE} is AM-regular at $\bar x$, then $\bar x$ is M-stationary for \eqref{RE}, and $\bar w$ is M-stationary for \eqref{P}.
\end{theorem}

By Section~\ref{sec:active_LICQ_AM}, the AM-regularity required in the last assertion of Theorem~\ref{thm:outer_AM_stationarity} holds if \eqref{P} admits no singular M-multiplier at $\bar w$, in particular if active-set SDCMPCC-LICQ holds at $\bar x$.

Suppose that each Step~3 solve is completed and $\nu_k\to0$. For a run continued with Step~1 omitted, the argument of \cite[proof of Proposition~4.1(a) and discussion following Theorem~4.3]{JiaKanzowMehlitzWachsmuth2022}, together with \eqref{eq:outer_residual_bounds}, gives $V_{k+1}\to0$ if $\{\rho_k\}$ is bounded and $z^{k+1}\to_{\mathcal K}0$ along every subsequence $x^{k+1}\to_{\mathcal K}\bar x\in\widetilde{\mathcal F}$. Since $V_{k+1}\leq\|z^{k+1}\|$ and both stationarity residuals at $x^{k+1}$ are at most $\nu_{k+1}$, the test in Step~1 is satisfied after finitely many iterations for any positive stopping tolerances whenever $\{\rho_k\}$ is bounded or the iterates have a feasible accumulation point.

\begin{corollary}[Local multiplier convergence]\label{cor:outer_multiplier_convergence}
Let $x^j=(w^j,u^j)\to\bar x=(\bar w,\bar u)\in\widetilde{\mathcal F}$, where active-set SDCMPCC-LICQ holds, and let $z^j,e^j\to0$. If $t^j\in N_C(\Psi(x^j)-z^j)\times\Nlim_\Omega(u^j)$ satisfies $\mathcal J_{w^j}t^j=e^j-\nabla\varphi(x^j)$, then $t^j$ converges to the unique M-stationarity multiplier tuple at $\bar x$.
\end{corollary}
\begin{proof}
The normal relations give $t^j\in\mathfrak L(w^j,u^j,z^j)$, so \eqref{eq:uniform_LICQ} yields $c\|t^j\|\leq\|e^j\|+\|\nabla f(w^j)\|$ for all sufficiently large $j$. By continuity and closedness of the normal-cone graphs, every cluster point is an M-stationarity multiplier tuple at $\bar x$. Any two such tuples differ by an element of $\ker\mathcal J_{\bar w}\cap\mathfrak L(\bar w,\bar u,0)=\{0\}$, proving convergence.
\end{proof}

For the run considered in
Theorem~\ref{thm:outer_AM_stationarity}, along any subsequence converging to a feasible point satisfying active-set SDCMPCC-LICQ, Corollary~\ref{cor:outer_multiplier_convergence} applies to the updated multipliers and any Step~3 normal witness. Indeed, \eqref{eq:outer_multiplier_normal} and \eqref{eq:outer_gradient_identity} give the required relations for $t^{k+1}=(\Upsilon^{k+1},\Gamma^{k+1})$ with $e^{k+1}=\nabla L_k(x^{k+1})+(0,\Gamma^{k+1})$ and $\|e^{k+1}\|\leq\sqrt{2}\,\nu_{k+1}$.

\begin{remark}[Comparison with Santos et al.~\cite{SantosCostaKrulikovskiGeremia2026}]
Both methods retain semidefinite complementarity and penalize smooth linking constraints. Santos et al. impose approximate KKT conditions for the scalar representation of semidefinite complementarity. Their Theorem~6.1 gives AW-stationarity at feasible limits and AC-stationarity under the additional coupling $\delta_k\leq\epsilon_k/(1+\gamma_k)$ in their Assumption~6.1. SPALM instead tests the limiting-normal distance directly in Step~3.
\end{remark}
The next two sections give verifiable block conditions and procedures for generating slack candidates.

\section{An inexact proximal alternating framework with verifiable block conditions}\label{sec:inner_solver}

We solve the augmented Lagrangian subproblem at each outer iteration by proximal alternating updates. The analysis concerns accepted iterates satisfying the prescribed descent and stationarity tests.

\subsection{Proximal block subproblems and the shifted projection formula}

Throughout this section, consider the subproblem at the $k$-th outer iteration and write $(\rho,\Upsilon):=(\rho_k,\bar\Upsilon^k)$. These parameters remain unchanged during the inner solve. Set
\begin{equation*}
L_\rho(w,u):=\mathcal L_\rho(w,u;\Upsilon),
\end{equation*}
where the augmented Lagrangian is given in \eqref{eq:AL}. The AL subproblem is
\begin{equation}\label{eq:inner}
\min_{w\in\mathcal Z,\;u\in\Omega}L_\rho(w,u).
\end{equation}

For fixed $w$, define the unregularized slack center
\begin{equation*}
  c_\rho(w)
  :=
  \left(
    G(w)+\frac{M^G}{\rho},
    H(w)+\frac{M^H}{\rho}
  \right)
  \in\mathcal U.
\end{equation*}
Up to a term independent of $u$, the slack-dependent part of $L_\rho(w,u)$ is $\frac\rho2\|u-c_\rho(w)\|^2$. Thus minimizing it over $\Omega$ is equivalent to projecting $c_\rho(w)$ onto $\Omega$.

To stabilize both block updates and obtain sufficient decrease, we
introduce proximal regularization. For proximal parameters
$a_t>0$ and $b_t>0$, define
\begin{equation}\label{eq:theta}
  \Theta^t(w)
  :=
  L_\rho(w,u^t)
  +
  \frac{a_t}{2}\norm{w-w^t}^2
\end{equation}
and
\begin{equation}\label{eq:phi}
  \Phi^t(u)
  :=
  L_\rho(w^{t+1},u)
  +
  \frac{b_t}{2}\norm{u-u^t}^2,
  \qquad
  u\in\Omega.
\end{equation}

Completing the square shows that $\Phi^t(u)$ differs from
$\frac{\rho+b_t}{2}\|u-\widehat c^t\|^2$ by a term independent of $u$, where
\begin{equation*}
\widehat c^t:=\frac{\rho c_\rho(w^{t+1})+b_tu^t}{\rho+b_t}.
\end{equation*}
Hence an exact proximal slack update satisfies
\begin{equation}\label{eq:shiftproj}
u_\star^{t+1}\in\operatorname*{argmin}_{u\in\Omega}\|u-\widehat c^t\|^2.
\end{equation}

Since $\Omega$ is nonempty and closed, \eqref{eq:shiftproj} has a solution. Every exact projection gives $\Phi^t(u_\star^{t+1})\leq\Phi^t(u^t)$ and
\begin{equation*}
-\nabla_u\Phi^t(u_\star^{t+1})=(\rho+b_t)(\widehat c^t-u_\star^{t+1})
\in\widehat N_\Omega(u_\star^{t+1})\subset\Nlim_\Omega(u_\star^{t+1}).
\end{equation*}

\begin{lemma}[Lipschitz dependence of $\nabla_wL_\rho$ on slack variables]\label{lem:wlip}
Assume $f,g,h,G,H\in C^1$. For every bounded set $\mathcal B\subset\mathcal Z$, there exists $L_{wu}<\infty$ such that
\begin{equation}\label{eq:wlip}
\|\nabla_wL_\rho(w,u)-\nabla_wL_\rho(w,u')\|
\leq L_{wu}\|u-u'\|
\qquad\forall w\in\mathcal B,\quad\forall u,u'\in\mathcal U.
\end{equation}
Here, for a linear operator $\mathcal A$ between finite-dimensional Euclidean spaces, $\norm{\mathcal A}_{\mathrm{op}}$ denotes its induced operator norm. One may take
\begin{equation*}
  L_{wu}:=\sqrt2\rho\sup_{w\in\mathcal B}
  \max\{\norm{G'(w)}_{\mathrm{op}},\norm{H'(w)}_{\mathrm{op}}\}.
\end{equation*}
\end{lemma}

\begin{proof}
Differentiating \eqref{eq:AL} and subtracting at $u=(Y,Z)$ and $u'=(Y',Z')$ gives
\begin{equation}\label{eq:wgrad}
\begin{aligned}
&\nabla_wL_\rho(w,u)-\nabla_wL_\rho(w,u')\\
&\qquad=-\rho G'(w)^*(Y-Y')-\rho H'(w)^*(Z-Z').
\end{aligned}
\end{equation}
The triangle inequality, $\|Y-Y'\|_F+\|Z-Z'\|_F\leq\sqrt2\|u-u'\|$, and boundedness of the continuous derivatives on $\overline{\mathcal B}$ prove the claim.

\end{proof}

For later use, define the pointwise slack Lipschitz factor
\begin{equation*}
  L_{wu}(w) :=
  \sqrt{2}\rho
  \max\left\{ \norm{G'(w)}_{\mathrm{op}}, \norm{H'(w)}_{\mathrm{op}}
  \right\}.
\end{equation*}
Then, at the $(t+1)$-st accepted iterate, we write $L_{wu}^{t+1}:=L_{wu}(w^{t+1})$.

For \eqref{eq:num_familyIprime_GH} and \eqref{eq:num_familyIIprime_GH}, $G(w)=X$ and $H(w)=-S$ act on separate components. The product norm in \eqref{eq:wgrad} therefore permits the smaller cross-Lipschitz constant $\rho$ in the subsequent residual bounds, in place of the generic factor $L_{wu}(w)$ defined above.

\subsection{Verifiable block updates and accepted iterates}

A trial point becomes an accepted iterate only after both block tests below pass.

\begin{algorithm}[H]
\caption{Inexact proximal alternating framework with verifiable block conditions}
\label{alg:verified_inner}
\footnotesize

\begin{algorithmic}
\algrenewcommand\algorithmicrequire{\textbf{Input:}}
\algrenewcommand\algorithmicensure{\textbf{Output:}}

\Require An initial accepted point
$x^0=(w^0,u^0)$ with $u^0\in\Omega$;
positive proximal-parameter sequences $\{a_t\}$ and $\{b_t\}$;
nonnegative tolerance sequences
$\{\delta_w^t\}$ and $\{\delta_u^t\}$; a prescribed target $\nu>0$.

\Ensure Accepted iterates $\{x^t\}$ satisfying the prescribed block conditions, together with normals $\Gamma^{t+1}\in\Nlim_\Omega(u^{t+1})$ certifying the slack residuals.

\State Set $t=0$.

\Statex \textbf{Step 1: inexact proximal $w$-update.}

\State Invoke a $w$-block oracle to generate trial points
for the proximal objective $\Theta^t$ in \eqref{eq:theta}.

\State Accept a trial point as $w^{t+1}$ only if
\begin{equation}\label{eq:wdescent}
  \Theta^t(w^{t+1})
  \leq
  \Theta^t(w^t)
\end{equation}
\State and
\begin{equation}\label{eq:wres}
  \norm{
    \nabla_wL_\rho(w^{t+1},u^t)
    +
    a_t(w^{t+1}-w^t)
  }
  \leq
  \delta_w^t.
\end{equation}

\Statex \textbf{Step 2: inexact proximal slack update.}

\State Invoke a $u$-block oracle to generate feasible trial points
for the proximal objective $\Phi^t$ in \eqref{eq:phi}.

\State Accept a trial point as $u^{t+1}\in\Omega$ only if
\begin{equation}\label{eq:udescent}
  \Phi^t(u^{t+1})
  \leq
  \Phi^t(u^t)
\end{equation}
\State and
\begin{equation}\label{eq:ures}
  \dist\left(
    -\nabla_u\Phi^t(u^{t+1}),
    \Nlim_\Omega(u^{t+1})
  \right)
  \leq
  \delta_u^t.
\end{equation}

\State Retain a normal $\Gamma^{t+1}\in\Nlim_\Omega(u^{t+1})$ such that
$\|\nabla_u\Phi^t(u^{t+1})+\Gamma^{t+1}\|\leq\delta_u^t$.
The computable sufficient test \eqref{eq:slack_certificate_test} supplies such a normal.

\Statex \textbf{Step 3: iterate update.}
\State Set $x^{t+1}:=(w^{t+1},u^{t+1})$. If $\widehat{\mathcal R}_\rho^{\,t+1}\leq\nu$, return this point and its normal. Otherwise increase $t$ and repeat.

\end{algorithmic}
\end{algorithm}
The slack residual condition \eqref{eq:ures} is equivalently
\begin{equation}\label{eq:ures_equiv}
\operatorname{dist}\left(
-\nabla_uL_\rho(w^{t+1},u^{t+1})-b_t(u^{t+1}-u^t),
N_\Omega^{\lim}(u^{t+1})\right)\leq\delta_u^t.
\end{equation}

For example, L-BFGS with line search \cite{LiuNocedal1989} can generate $w$-trials, and Section~\ref{sec:slack_solver} supplies feasible slack trials. Budget exhaustion before acceptance is reported as \emph{target not verified}.

For target $\nu>0$, the full-residual bound \eqref{eq:computable_full_residual} is checked after both updates. Returning $x^{t+1}$ and its retained normal when $\widehat{\mathcal R}_\rho^{\,t+1}\leq\nu$ implies Step~3 of Algorithm~\ref{alg:sALcc}. Conversely, its two componentwise tests only imply a combined residual at most $\sqrt2\nu$. Omitting the stopping test defines the infinite accepted sequences below.

Both the existence argument and the convergence analysis use the initial sublevel set

\begin{equation}\label{eq:inner_sublevel}
\mathcal S_\rho(x^0)
:=\{x\in\mathcal Z\times\Omega:L_\rho(x)\leq L_\rho(x^0)\}.
\end{equation}
If this set is bounded, it is compact, since $L_\rho$ is continuous and $\mathcal Z\times\Omega$ is closed. In particular, $L_\rho$ attains a finite minimum on it.

\begin{proposition}[Existence of accepted block updates]
\label{prop:inner_block_existence}
Assume that $f,g,h,G,H\in C^1$ and that $\mathcal S_\rho(x^0)$ is bounded. At every accepted iterate $x^t$, and for any $a_t,b_t>0$, exact minimizers
\begin{equation*}
w_\star^{t+1}\in\argmin_{w\in\mathcal Z}\Theta^t(w),
\qquad u_\star^{t+1}\in\argmin_{u\in\Omega}\Phi^t(u)
\end{equation*}
exist, where $w^{t+1}=w_\star^{t+1}$ is used to define $\Phi^t$. These minimizers satisfy the two block acceptance conditions with zero residuals.
\end{proposition}

\begin{proof}
The accepted descent conditions give $L_\rho(x^t)\leq L_\rho(x^0)$. The set $\{w:\Theta^t(w)\leq\Theta^t(w^t)\}$ is nonempty and closed. Every point in this set satisfies
\begin{equation*}
L_\rho(w,u^t)\leq\Theta^t(w)\leq L_\rho(x^t)\leq L_\rho(x^0),
\end{equation*}
so the set is bounded by the assumption on $\mathcal S_\rho(x^0)$. Thus $\Theta^t$ attains a global minimum, and its gradient vanishes there. This gives \eqref{eq:wdescent} and \eqref{eq:wres} with zero residual. The projection representation \eqref{eq:shiftproj} provides an exact slack update satisfying \eqref{eq:udescent} and \eqref{eq:ures} with zero residual. Induction therefore gives accepted updates at every iteration when exact block minimizers are available.
\end{proof}

\subsection{Convergence and finite-accuracy guarantees for accepted sequences}
We establish convergence and finite-accuracy guarantees using sufficient decrease and estimates of the full stationarity residual.

\begin{lemma}[Sufficient decrease and vanishing steps]
\label{lem:suffdec}
Suppose that Algorithm~\ref{alg:verified_inner} generates a finite or infinite accepted sequence $\{x^t=(w^t,u^t)\}$, with a constant $\cmin>0$ such that
\begin{equation*}
  a_t\geq\cmin,
  \qquad
  b_t\geq\cmin
  \qquad
  \forall t.
\end{equation*}
Then every accepted transition satisfies
\begin{equation}\label{eq:wblockdec}
  L_\rho(w^t,u^t)-L_\rho(w^{t+1},u^t)
  \geq
  \frac{a_t}{2}\norm{w^{t+1}-w^t}^{2},
\end{equation}
and
\begin{equation}\label{eq:ublockdec}
  L_\rho(w^{t+1},u^t)-L_\rho(w^{t+1},u^{t+1})
  \geq
  \frac{b_t}{2}\norm{u^{t+1}-u^t}^{2}.
\end{equation}
Consequently,
\begin{equation}\label{eq:jointdec}
  L_\rho(x^t)-L_\rho(x^{t+1})
  \geq
  \frac{\cmin}{2}\norm{x^{t+1}-x^t}^{2}.
\end{equation}

If $\mathcal S_\rho(x^0)$ is bounded and the sequence is infinite, then
\begin{equation*}
  \sum_{t=0}^{\infty}
  \norm{x^{t+1}-x^t}^{2}
  <\infty,
\end{equation*}
and hence
\begin{equation}\label{eq:vanishing_steps}
  \norm{x^{t+1}-x^t}\to0.
\end{equation}
In particular,
\begin{equation*}
  \norm{w^{t+1}-w^t}\to0,
  \qquad
  \norm{u^{t+1}-u^t}\to0.
\end{equation*}
\end{lemma}
\begin{proof}
Expanding \eqref{eq:wdescent} and \eqref{eq:udescent} gives \eqref{eq:wblockdec} and \eqref{eq:ublockdec}; adding them proves \eqref{eq:jointdec}. All accepted points remain in $\mathcal S_\rho(x^0)$. If this set is bounded, continuity gives a finite lower bound $m_\rho$ for $L_\rho$ there. Summation yields
\begin{equation*}
\frac{\cmin}{2}\sum_{t=0}^{T}\|x^{t+1}-x^t\|^2
\leq L_\rho(x^0)-L_\rho(x^{T+1})\leq L_\rho(x^0)-m_\rho,
\end{equation*}
which proves square summability and \eqref{eq:vanishing_steps}.

\end{proof}

For $x=(w,u)\in\mathcal Z\times\Omega$, define the full M-stationarity residual of the inner subproblem by
\begin{equation}\label{eq:fullresidual}
  \mathcal R_\rho(x):=
  \left(\norm{\nabla_wL_\rho(w,u)}^2+
    \dist^2\left(-\nabla_uL_\rho(w,u),\Nlim_\Omega(u)\right)
  \right)^{1/2}.
\end{equation}

\begin{lemma}[One-step residual bound and normal-certificate transfer]
\label{lem:block_full_residual}
Assume that $f,g,h,G,H\in C^1$. Let $x^t=(w^t,u^t)$ and $x^{t+1}=(w^{t+1},u^{t+1})$ be consecutive accepted iterates of Algorithm~\ref{alg:verified_inner}, with $a_t,b_t>0$. Write
\begin{equation*}
\Delta w^t:=w^{t+1}-w^t,\qquad
\Delta u^t:=u^{t+1}-u^t=(\Delta Y^t,\Delta Z^t).
\end{equation*}
There exists $\Gamma^{t+1}\in\Nlim_\Omega(u^{t+1})$ such that the block errors
\begin{equation*}
e_w^t:=\nabla_wL_\rho(w^{t+1},u^t)+a_t\Delta w^t,
\qquad
e_u^t:=\nabla_u\Phi^t(u^{t+1})+\Gamma^{t+1}
\end{equation*}
satisfy $\|e_w^t\|\leq\delta_w^t$ and $\|e_u^t\|\leq\delta_u^t$. For any such normal,
\begin{equation}\label{eq:inner_error_transfer}
\begin{aligned}
\nabla_wL_\rho(x^{t+1})
&=e_w^t-a_t\Delta w^t
-\rho G'(w^{t+1})^*\Delta Y^t
-\rho H'(w^{t+1})^*\Delta Z^t,\\
\nabla_uL_\rho(x^{t+1})+\Gamma^{t+1}
&=e_u^t-b_t\Delta u^t.
\end{aligned}
\end{equation}
Consequently,
\begin{equation}\label{eq:fullwres}
\|\nabla_wL_\rho(x^{t+1})\|
\leq L_{wu}^{t+1}\|\Delta u^t\|+a_t\|\Delta w^t\|+\delta_w^t
\end{equation}
and
\begin{equation}\label{eq:fullures}
\begin{aligned}
\dist\bigl(-\nabla_uL_\rho(x^{t+1}),\Nlim_\Omega(u^{t+1})\bigr)
&\leq\|\nabla_uL_\rho(x^{t+1})+\Gamma^{t+1}\|\\
&\leq b_t\|\Delta u^t\|+\delta_u^t.
\end{aligned}
\end{equation}
Define the computable upper bound
\begin{equation}\label{eq:computable_full_residual}
\widehat{\mathcal R}_\rho^{\,t+1}
:=\left[
\bigl(L_{wu}^{t+1}\|\Delta u^t\|+a_t\|\Delta w^t\|+\delta_w^t\bigr)^2
+\bigl(b_t\|\Delta u^t\|+\delta_u^t\bigr)^2
\right]^{1/2}
\end{equation}
and the full error vector
\begin{equation*}
E^{t+1}:=\bigl(\nabla_wL_\rho(x^{t+1}),
\nabla_uL_\rho(x^{t+1})+\Gamma^{t+1}\bigr).
\end{equation*}
Then
\begin{equation}\label{eq:full_residual_upper_bound}
\mathcal R_\rho(x^{t+1})\leq\|E^{t+1}\|
\leq\widehat{\mathcal R}_\rho^{\,t+1}.
\end{equation}
\end{lemma}

\begin{proof}
The distance to the nonempty closed set $\Nlim_\Omega(u^{t+1})$ is attained; a normal retained from \eqref{eq:slack_certificate_test} also suffices. The definitions of $e_w^t,e_u^t$ and \eqref{eq:wgrad} give \eqref{eq:inner_error_transfer}. Taking norms and applying Lemma~\ref{lem:wlip} proves the bounds. The same normal is used throughout; the proximal term $-b_t\Delta u^t$ is part of the error.

\end{proof}

\begin{theorem}[Residual convergence and M-stationarity]
\label{thm:main}
Assume that $f,g,h,G,H\in C^1$, and that there exist constants $0<\cmin\leq\cmax<\infty$ such that
\begin{equation*}
\cmin\leq a_t,b_t\leq\cmax\qquad\text{for all }t.
\end{equation*}
Suppose that Algorithm~\ref{alg:verified_inner} generates an infinite accepted sequence $\{x^t=(w^t,u^t)\}$, that $\mathcal S_\rho(x^0)$ is bounded, and that $\delta_w^t,\delta_u^t\to0$. Then:
\begin{enumerate}
\item[(i)] The accepted sequence $\{x^t\}$ is bounded.
\item[(ii)] The computable upper bound and the full residual satisfy
\begin{equation}\label{eq:fullresidual_zero}
\widehat{\mathcal R}_\rho^{\,t+1}\to0,
\qquad \mathcal R_\rho(x^{t+1})\to0.
\end{equation}
\item[(iii)] Every accumulation point $\bar x=(\bar w,\bar u)$ of $\{x^t\}$ is M-stationary for the AL subproblem \eqref{eq:inner}; that is,
\begin{equation*}
\nabla_wL_\rho(\bar w,\bar u)=0,
\qquad -\nabla_uL_\rho(\bar w,\bar u)\in\Nlim_\Omega(\bar u).
\end{equation*}
\end{enumerate}
\end{theorem}

\begin{proof}
Lemma~\ref{lem:suffdec} places all accepted iterates in the compact set $\mathcal S_\rho(x^0)$ and gives $\|x^{t+1}-x^t\|\to0$. This proves (i). Continuity of $G'$ and $H'$ implies boundedness of $\{L_{wu}^{t+1}\}$. Together with the parameter bounds and vanishing block errors, \eqref{eq:computable_full_residual} gives $\widehat{\mathcal R}_\rho^{\,t+1}\to0$. Lemma~\ref{lem:block_full_residual} then proves (ii).

Let $x^{t_j}\to\bar x$. Vanishing steps give $x^{t_j+1}\to\bar x$, and closedness of $\Omega$ gives $\bar u\in\Omega$. Use the normals $\Gamma^{t_j+1}$ from Lemma~\ref{lem:block_full_residual}. Since $\|E^{t_j+1}\|\to0$, continuity of the two gradients yields
\begin{equation*}
\nabla_wL_\rho(\bar x)=0,
\qquad \Gamma^{t_j+1}\to-\nabla_uL_\rho(\bar x).
\end{equation*}
The closed graph of the limiting normal-cone mapping implies $-\nabla_uL_\rho(\bar x)\in\Nlim_\Omega(\bar u)$, proving (iii).
\end{proof}

Propositions~\ref{prop:familyIprime_level_bounded} and~\ref{prop:familyIIprime_level_bounded} verify the sublevel-set boundedness assumption for our numerical models.

\begin{corollary}[Finite-accuracy output]
\label{cor:inner_finite_accuracy}\label{cor:fixed_accuracy_inner}
Assume $C^1$ data, boundedness of $\mathcal S_\rho(x^0)$, $0<\cmin\leq a_t,b_t\leq\cmax<\infty$, and an infinite accepted sequence. Fix $\nu>0$. If either
\emph{(i)} $\delta_w^t,\delta_u^t\to0$, or
\emph{(ii)} $0\leq\delta_w^t\leq c_w\nu$, $0\leq\delta_u^t\leq c_u\nu$ with $\kappa:=\sqrt{c_w^2+c_u^2}<1$, then eventually
\begin{equation*}
\mathcal R_\rho(x^{t+1})\leq\|E^{t+1}\|\leq\widehat{\mathcal R}_\rho^{\,t+1}<\nu.
\end{equation*}
For the AL subproblem at outer iteration $k$ with parameters $(\rho_k,\bar\Upsilon^k)$, choose $x^{k,0}\in\mathcal Z\times\Omega$, and impose the preceding assumptions on its initial sublevel set and accepted sequence. Denote its accepted iterates, retained normals and bounds by $x^{k,t}$, $\Gamma^{k,t+1}$ and $\widehat{\mathcal R}_k^{\,t+1}$. The initial point need not equal $x^k$, and the sublevel-set bound need not be uniform in $k$. If a feasible reference point is supplied, choose $x^{k,0}$ by \eqref{eq:reference_initialization}. Lemma~\ref{lem:suffdec} then gives $L_k(x^{k,t})\leq L_k(x^{k,0})$ for every accepted iterate, so the returned point also satisfies the value comparison. Under either regime with target $\nu_{k+1}$, the index
\begin{equation}\label{eq:inner_stopping_index}
t_k:=\min\{t\in\mathbb N_0:\widehat{\mathcal R}_k^{\,t+1}\leq\nu_{k+1}\}
\end{equation}
is finite. Returning $x^{k+1}=x^{k,t_k+1}$ and $\Gamma^{k+1}=\Gamma^{k,t_k+1}\in\Nlim_\Omega(u^{k+1})$ gives
\begin{equation*}
\|\bigl(\nabla_wL_k(x^{k+1}),\nabla_uL_k(x^{k+1})+\Gamma^{k+1}\bigr)\|\leq\nu_{k+1},
\end{equation*}
and hence both conditions in Step~3 of Algorithm~\ref{alg:sALcc}.
\end{corollary}
\begin{proof}
In case \emph{(i)}, apply Theorem~\ref{thm:main}. In case \emph{(ii)}, Lemma~\ref{lem:suffdec} gives vanishing steps, while $L_{wu}^{t+1},a_t,b_t$ are bounded. Thus \eqref{eq:computable_full_residual} gives $\limsup_t\widehat{\mathcal R}_\rho^{\,t+1}\leq\kappa\nu<\nu$. The returned bound follows from Lemma~\ref{lem:block_full_residual} at that same point and normal.
\end{proof}
For example, $c_w=c_u=1/2$ is valid. The retained normal can also be used in \eqref{eq:outer_finite_certificate} when $\nu_{k+1}\leq\varepsilon_{\mathrm{stat}}$, or when its slack residual is directly verified against $\varepsilon_{\mathrm{stat}}$.

\begin{corollary}[Explicit residual bounds in accepted iterations]
\label{cor:inner_complexity}
Assume that $f,g,h,G,H\in C^1$, that $\mathcal S_\rho(x^0)$ is bounded, and that $0<\cmin\leq a_t,b_t\leq\cmax<\infty$ for all accepted transitions. Let $\mathcal B_\rho$ be the projection of $\mathcal S_\rho(x^0)$ onto the $w$-space, and define
\begin{equation*}
\begin{aligned}
m_\rho&:=\min_{x\in\mathcal S_\rho(x^0)}L_\rho(x),
&\Delta_\rho&:=L_\rho(x^0)-m_\rho,\\
\overline L_{wu}&:=\sup_{w\in\mathcal B_\rho}L_{wu}(w),
&K_\rho&:=\sqrt{2\cmax^{\,2}+\overline L_{wu}^{\,2}}.
\end{aligned}
\end{equation*}
These constants are finite. For any $N\geq1$ accepted transitions, put
\begin{equation*}
\epsilon_t:=\sqrt{(\delta_w^t)^2+(\delta_u^t)^2},
\qquad E_N:=\sum_{t=0}^{N-1}\epsilon_t^2.
\end{equation*}
Then
\begin{equation}\label{eq:inner_finite_horizon_bound}
\min_{0\leq t<N}\bigl(\widehat{\mathcal R}_\rho^{\,t+1}\bigr)^2
\leq\frac{4K_\rho^2\Delta_\rho/\cmin+2E_N}{N}.
\end{equation}
In particular, along an infinite accepted sequence with $E_\infty:=\sum_{t=0}^\infty\epsilon_t^2<\infty$,
\begin{equation*}
\min_{0\leq t<N}\widehat{\mathcal R}_\rho^{\,t+1}
\leq\sqrt{\frac{4K_\rho^2\Delta_\rho/\cmin+2E_\infty}{N}}
=O(N^{-1/2}).
\end{equation*}
If instead $\epsilon_t\leq\kappa\nu$ for a fixed $\nu>0$ and $0\leq\kappa<1$, then
\begin{equation}\label{eq:inner_fixed_accuracy_bound}
\min_{0\leq t<N}\widehat{\mathcal R}_\rho^{\,t+1}
\leq\kappa\nu+K_\rho\sqrt{\frac{2\Delta_\rho}{\cmin N}}.
\end{equation}
The same upper bounds hold with $\mathcal R_\rho(x^{t+1})$ in place of $\widehat{\mathcal R}_\rho^{\,t+1}$.
\end{corollary}

\begin{proof}
Compactness of $\mathcal S_\rho(x^0)$ and continuity of the data and their derivatives ensure finiteness of the constants. Let $s_t:=\|x^{t+1}-x^t\|$. By the product-norm triangle inequality in \eqref{eq:computable_full_residual} and
\begin{equation*}
\left\|\begin{pmatrix}a_t&L_{wu}^{t+1}\\0&b_t\end{pmatrix}\right\|_{\mathrm{op}}
\leq
\left\|\begin{pmatrix}a_t&L_{wu}^{t+1}\\0&b_t\end{pmatrix}\right\|_F
\leq K_\rho,
\end{equation*}
we obtain $\widehat{\mathcal R}_\rho^{\,t+1}\leq K_\rho s_t+\epsilon_t$. Summing the finite-transition decrease inequality \eqref{eq:jointdec} gives
\begin{equation*}
\sum_{t=0}^{N-1}s_t^2
\leq\frac{2\bigl(L_\rho(x^0)-L_\rho(x^N)\bigr)}{\cmin}
\leq\frac{2\Delta_\rho}{\cmin}.
\end{equation*}
Squaring the pointwise residual bound, using $(p+q)^2\leq2p^2+2q^2$, and averaging over $t=0,\ldots,N-1$ proves \eqref{eq:inner_finite_horizon_bound}. The square-summable case follows by replacing $E_N$ with $E_\infty$. If $\epsilon_t\leq\kappa\nu$, choose an index minimizing $s_t$ to obtain
\begin{equation*}
\min_{0\leq t<N}\widehat{\mathcal R}_\rho^{\,t+1}
\leq\kappa\nu+K_\rho\min_{0\leq t<N}s_t
\leq\kappa\nu+K_\rho\sqrt{\frac{2\Delta_\rho}{\cmin N}}.
\end{equation*}
Finally, Lemma~\ref{lem:block_full_residual} gives the bounds for $\mathcal R_\rho$.
\end{proof}

Under the fixed-accuracy allocation in Corollary~\ref{cor:fixed_accuracy_inner}, \eqref{eq:inner_fixed_accuracy_bound} guarantees a point with $\widehat{\mathcal R}_\rho^{\,t+1}<\nu$ among the first $N$ accepted transitions whenever
\begin{equation*}
N>B_\nu:=\frac{2K_\rho^2\Delta_\rho}{\cmin(1-\kappa)^2\nu^2}.
\end{equation*}
Thus $\lfloor B_\nu\rfloor+1$ is a sufficient number of accepted transitions for the first tolerance hit. 

A known lower bound for $L_\rho$ on $\mathcal S_\rho(x^0)$ may replace $m_\rho$ in these estimates. Usable bounds for the remaining constants then give a numerical upper bound on the number of accepted transitions required to reach the target.

For a fixed AL subproblem and initial sublevel set, the bound is $O(\nu^{-2})$ when $\cmin>0$, $\cmax$, and $\kappa<1$ are independent of $\nu$. It counts accepted transitions; the cost of generating and verifying trial points is additional, and the constants may vary across outer iterations.

In the next section, we use the shifted projection representation to generate feasible slack candidates, which are then checked against the descent and stationarity tests.

\section{The slack \texorpdfstring{$u$}{u}-subproblem: feasible-candidate generation}\label{sec:slack_solver}

At the $t$-th inner iteration, write
\begin{equation*}
  \widehat c^t=(A^t,B^t),
  \qquad A^t,B^t\in\Sn.
\end{equation*}
Since the projection routines developed below apply to a generic center, we suppress the iteration index and write $(A,B)=(A^t,B^t)$. The exact proximal slack update then reduces to
\begin{equation}\label{eq:projAB}
  \min_{(Y,Z)\in\Omega}
  \ \norm{Y-A}_F^{2}+\norm{Z-B}_F^{2}.
\end{equation}

\begin{remark}[Complexity of exact projection]\label{rem:projection-hardness} When $n$ is part of the input, deciding whether the optimal value of \eqref{eq:projAB} is at most a given rational threshold is NP-hard for rational $A,B$. A reduction from Simple Max-Cut \cite{GareyJohnsonStockmeyer1976} constructs centers $A=D+\varepsilon C$ and $B=-D+\varepsilon C$, with $D$ positive diagonal and $C$ encoding the graph. A uniform second-order estimate allows the maximum cut size to be recovered from the projection value. The rational construction is computable in polynomial time, yielding the claimed reduction. The complete proof is given in Appendix~\ref{app:projection-hardness}. \end{remark}

An exact solution of \eqref{eq:projAB} is computationally expensive. We therefore use a common diagonal reconstruction to generate feasible candidates through either a Huber--Stiefel search or a fixed eigenbasis of $A$ or $B$.

For the proximal slack update, set $q_t(u):=\|u-\widehat c^t\|^2$. Completing the square in \eqref{eq:phi} gives
\begin{equation}\label{eq:slack_true_descent_bridge}
\Phi^t(u)=C_t+\frac{\rho+b_t}{2}q_t(u),
\qquad
\Phi^t(u_{\mathrm{cand}})\leq\Phi^t(u^t)
\ \Longleftrightarrow\ q_t(u_{\mathrm{cand}})\leq q_t(u^t),
\end{equation}
where $C_t$ is independent of $u$. Thus candidate generation is followed by a comparison of the true squared distances and a separate normal-residual test at the same candidate.

\subsection{Diagonal decoupling for a fixed orthogonal matrix}\label{subsec:diagonal_decoupling}

Every $(Y,Z)\in\Omega$ has a joint representation $Y=P\Diag(y)P^\top$, $Z=P\Diag(z)P^\top$, with $P\in\On$, $y,z\in\mathbb{R}^n$, $y\geq0$, $z\leq0$ and $y_iz_i=0$; conversely, these conditions imply membership in $\Omega$.

For such a matrix $P$, define the diagonal vectors
\begin{equation*}d_A(P):=\diag(P^{\top}AP), \qquad d_B(P):=\diag(P^{\top}BP),\end{equation*}
whose $i$-th components are denoted by $d_{A,i}(P)$ and $d_{B,i}(P)$, respectively.

For fixed $P$, orthogonal invariance of the Frobenius norm
implies that the off-diagonal terms are independent of $y$ and $z$.
Hence, the minimization over $y$ and $z$ separates into the following
scalar complementarity problems, one for each $i=1,\ldots,n$, as in \cite[Section~4, equations~(27)--(29)]{GuoDeng2022}:

\begin{equation*}
  \min_{y_i,z_i}\ (y_i-d_{A,i}(P))^2+(z_i-d_{B,i}(P))^2
  \quad\text{s.t.}\quad
  y_i\geq0,
  \ z_i\leq0,
  \ y_iz_i=0.
\end{equation*}
The minimizers over the two complementarity branches are
\begin{equation*}
  (y_i,z_i)=((d_{A,i}(P))_+,0)
  \quad\text{and}\quad
  (y_i,z_i)=(0,-(d_{B,i}(P))_-),
\end{equation*} respectively. Select the better branch by
\begin{equation*}
  \chi_i(P):=\begin{cases}
  1, & (d_{A,i}(P))_+^2\geq(d_{B,i}(P))_-^2,\\
  0, & (d_{A,i}(P))_+^2<(d_{B,i}(P))_-^2.
  \end{cases}\end{equation*}
Ties are resolved by choosing $\chi_i(P)=1$. Define the scalar reconstruction by
\begin{equation*}
y_i(P):=\chi_i(P)(d_{A,i}(P))_+,
\qquad z_i(P):=-(1-\chi_i(P))(d_{B,i}(P))_-.
\end{equation*}
Collect these coefficients into the vectors
\[
y(P):=\bigl(y_1(P),\ldots,y_n(P)\bigr)^\top,
\qquad
z(P):=\bigl(z_1(P),\ldots,z_n(P)\bigr)^\top.
\]
The corresponding matrices are
\begin{equation}\label{eq:reconY}
Y(P):=P\Diag(y(P))P^{\top},
\qquad
Z(P):=P\Diag(z(P))P^{\top}.
\end{equation}
Thus $u(P):=(Y(P),Z(P))\in\Omega$.

For each fixed $P$, the reconstructed candidate minimizes the
diagonal subproblem. Its objective value is
\begin{equation}\label{eq:minvalP}
\begin{aligned}
V(P)
&:=\|Y(P)-A\|_F^2+\|Z(P)-B\|_F^2\\
&=\|A\|_F^2+\|B\|_F^2-F(P),
\end{aligned}
\end{equation}
where
\begin{equation*}
\begin{aligned}
F(P)
&:=\sum_{i=1}^n
\max\{(d_{A,i}(P))_+^2,(d_{B,i}(P))_-^2\}\\
&=\sum_{i=1}^n
\bigl(\max\{d_{A,i}(P),-d_{B,i}(P),0\}\bigr)^2.
\end{aligned}
\end{equation*}
Thus, the search over $P$ reduces to the nonsmooth
maximization problem $\displaystyle\max_{P\in\On} F(P)$. 

\subsection{Huber--Stiefel feasible-candidate generation}

To obtain a continuously differentiable approximation of $F$,
we smooth the maximum of $d_{A,i}(P)$ and $-d_{B,i}(P)$
in each summand. For a smoothing parameter $\delta>0$, we use the following Huber-type approximation of the absolute value \cite{Huber1964}:
\begin{equation*}
  H_\delta(s):=\begin{cases}
  \dfrac{s^2}{2\delta}+\dfrac{\delta}{2}, & |s|\leq\delta,\\[0.7em]
  |s|, & |s|>\delta,
  \end{cases}\end{equation*}
with derivative
\begin{equation*}H_\delta'(s):=\begin{cases}
  s/\delta, & |s|\leq\delta,\\
  \operatorname{sign}(s), & |s|>\delta.
  \end{cases}\end{equation*}

A continuously differentiable approximation of $\max\{a,b\}$ is
\begin{equation*}M_\delta(a,b):=\frac{a+b+H_\delta(a-b)}{2}.\end{equation*}

Since $0\leq H_\delta(s)-|s|\leq\frac{\delta}{2}$, the smoothed maximum satisfies
\begin{equation}\label{eq:Mdelta_scalar_bound}
  0\leq M_\delta(a,b)-\max\{a,b\} \leq \frac{\delta}{4}.
\end{equation}
Thus $M_\delta(a,b)$ converges uniformly to $\max\{a,b\}$ as $\delta\downarrow0$.

Set
\begin{equation*}m_i(P):=M_\delta(d_{A,i}(P),-d_{B,i}(P))
  =\frac{d_{A,i}(P)-d_{B,i}(P)+H_\delta(d_{A,i}(P)+d_{B,i}(P))}{2}.\end{equation*}
The smoothed Stiefel objective is
\begin{equation}\label{eq:Fdelta}
  F_\delta(P):= \sum_{i=1}^n
\bigl(m_i(P)\bigr)_+^2,
  \qquad P\in\On.
\end{equation}
Since the scalar function $s\mapsto (s_+)^2$ is continuously differentiable and $M_\delta$ is $C^1$, the function $F_\delta$ is continuously differentiable on $\On$.

The following result quantifies the approximation error uniformly
over $\mathbb{O}^n$.
\begin{proposition}[Uniform approximation of the Stiefel objective]
\label{prop:Fdelta_uniform}
Let
\begin{equation*}R_0:=\max\{\|A\|_2,\|B\|_2\}. \end{equation*}
Then, for every $\delta>0$ and every $P\in\On$,
\begin{equation}\label{eq:Fdelta_uniform_bound}
  0\leq F_\delta(P)-F(P) \leq
  n\left( \frac{R_0\delta}{2}
+\frac{\delta^2}{16} \right).
\end{equation}
Hence $F_\delta\to F$ uniformly on $\On$ as $\delta\downarrow0$.

\end{proposition}

\begin{proof}
For $P\in\On$, set $q_i(P):=\max\{d_{A,i}(P),-d_{B,i}(P)\}$. By \eqref{eq:Mdelta_scalar_bound}, $0\leq m_i(P)-q_i(P)\leq\delta/4$. Since the mapping $s\mapsto s_+$ is nondecreasing and
$1$-Lipschitz,
\begin{equation}\label{eq:positive_part_uniform_bound}
0\leq(m_i(P))_+-(q_i(P))_+\leq\delta/4.
\end{equation}
Since $0\leq(q_i(P))_+\leq R_0$, the difference of squares is at most $(\delta/4)(2R_0+\delta/4)$. Summing over $i$ proves \eqref{eq:Fdelta_uniform_bound}.

\end{proof}

To carry out the Stiefel search, we compute the Euclidean gradient of $F_\delta$. Let
\begin{equation*}
  \zeta_i:=H_\delta'(d_{A,i}(P)+d_{B,i}(P)),
  \qquad
  \gamma_i:=2(m_i(P))_+,
\end{equation*}
and define
\begin{equation*}
  D_\gamma:=\Diag(\gamma_1,\ldots,\gamma_n),
  \qquad
  D_{\gamma\zeta}:=\Diag(\gamma_1\zeta_1,\ldots,\gamma_n\zeta_n).
\end{equation*}
By direct differentiation and the chain rule, the Euclidean gradient of $F_\delta$ is
\begin{equation}\label{eq:eucgrad}
  \nabla_{P}F_\delta(P)=(A-B)P D_\gamma+(A+B)P D_{\gamma\zeta}.
\end{equation}

Under the Frobenius metric, projection onto $T_P\On=\{\eta:P^\top\eta+\eta^\top P=0\}$ gives the Riemannian gradient \cite{EdelmanAriasSmith1998,AbsilMahonySepulchre2008}:

\begin{equation}\label{eq:riemgrad}
  \grad F_\delta(P)=\nabla_P F_\delta(P)
  -P\sym\!\left(P^{\top}\nabla_P F_\delta(P)\right).
\end{equation}
Here $\sym(C):=(C+C^\top)/2$ for any square matrix $C$.

For the manifold search below, we use the standard QR-based retraction \cite{AbsilMahonySepulchre2008},
\begin{equation*}
  R_P(\eta):=  \operatorname{qf}(P+\eta),
\end{equation*}
where $\operatorname{qf}(X)$ denotes the orthogonal factor $Q$ in the QR factorization $X=QR$, with the diagonal entries of $R$ chosen to be positive.

Since $F_\delta$ is continuously differentiable and $R_P$ is a retraction, whenever $\grad F_\delta(P)\neq0$, the ascent direction $\eta=\grad F_\delta(P)$ satisfies the standard ascent Armijo condition for all sufficiently small $\tau>0$. Hence, the
backtracking procedure in the algorithm below terminates after finitely many reductions.

\begin{algorithm}[H]
\caption{Huber--Stiefel feasible-candidate generator}
\label{alg:stiefel_candidate}
\footnotesize

\begin{algorithmic}
\algrenewcommand\algorithmicrequire{\textbf{Input:}}
\algrenewcommand\algorithmicensure{\textbf{Output:}}

\Require
$A,B\in\Sn$; initial $P^0\in\On$; smoothing parameter
$\delta>0$; Armijo parameters $\alpha>0$,
$\beta_{\mathrm{ls}}\in(0,1)$, $\sigma\in(0,1)$; stopping tolerances
$\eta_{\mathrm{st}}>0$, $\eta_F>0$; maximum number of iterations $K_{\max}$.

\Ensure
A feasible candidate
$u_{\mathrm{cand}}=(Y_{\mathrm{cand}},Z_{\mathrm{cand}})\in\Omega$.

\State Set $s=0$.
  \Statex \textbf{Step 1: Gradient and ascent direction.}
  \State Compute $\grad F_\delta(P^s)$ from
  \eqref{eq:eucgrad}--\eqref{eq:riemgrad} and set
  $\eta^s:=\grad F_\delta(P^s)$.

\Statex \textbf{Step 2: Armijo backtracking and update.}
\State Starting from $\tau=\alpha$, reduce
$\tau\leftarrow\beta_{\mathrm{ls}}\tau$ until, with
\begin{equation*}
  P_{\mathrm{cand}}:=R_{P^s}(\tau\eta^s),
\end{equation*}
the Armijo condition
\begin{equation}\label{eq:armijo}
  F_\delta(P_{\mathrm{cand}})-F_\delta(P^s)
  \geq
  \sigma\tau
  \ip{\grad F_\delta(P^s)}{\eta^s}
\end{equation}
holds. Set $P^{s+1}:=P_{\mathrm{cand}}$.

  \Statex \textbf{Step 3: Stiefel stopping tests.}
  \State If
\begin{equation*}
\|\operatorname{grad}F_\delta(P^{s+1})\|_F<\eta_{\mathrm{st}},
  \qquad\text{or}\qquad
  |F_\delta(P^{s+1})-F_\delta(P^s)|<\eta_F,
\end{equation*}
or if $s+1\geq K_{\max}$, proceed to Step~4.
Otherwise, set $s\leftarrow s+1$ and return to Step~1.

\Statex \textbf{Step 4: Feasible reconstruction.}
\State Set
\begin{equation*}
  u_{\mathrm{cand}}
  :=
  \bigl(
    Y(P^{s+1}),
    Z(P^{s+1})
  \bigr)
\end{equation*}
using \eqref{eq:reconY}, and return $u_{\mathrm{cand}}\in\Omega$.

\end{algorithmic}
\end{algorithm}

\paragraph{Post-verification and acceptance.}

By \eqref{eq:minvalP}, $q_t(u(P))=V(P)$. If $P^0$ is a common
eigenbasis of the current feasible point $u^t$, scalar minimization
yields $V(P^0)\leq q_t(u^t)$. Increasing $F_\delta$, however,
need not preserve this bound.

Define
\begin{equation*}
E_\delta:=n\left(\frac{R_0\delta}{2}+\frac{\delta^2}{16}\right).
\end{equation*}
By \eqref{eq:Fdelta_uniform_bound}, for any $P_0,P_1\in\On$,
\begin{equation}\label{eq:smoothed_true_gain}
F(P_1)-F(P_0)
\geq F_\delta(P_1)-F_\delta(P_0)-E_\delta.
\end{equation}
Thus, under the above initialization, a smoothed gain of at least
$E_\delta$ from $P^0$ ensures
$q_t(u(P_1))=V(P_1)\leq V(P^0)\leq q_t(u^t)$.

The gain condition is sufficient but not necessary for descent. Acceptance is determined by the true-value comparison in \eqref{eq:slack_true_descent_bridge} and the certificate test \eqref{eq:slack_certificate_test}. The following example explains why the normal-residual test is needed in addition to the Stiefel stopping criteria.

\begin{remark}[Stiefel stopping versus slack verification]
A zero Stiefel residual need not give a small slack residual. Take
\begin{equation*}
A=\Diag(2,0),\quad B=0,\quad
P=2^{-1/2}\begin{pmatrix}1&-1\\1&1\end{pmatrix},\quad 0<\delta<1.
\end{equation*}
Then $d_A(P)=(1,1)^\top$, $d_B(P)=0$, $m_i=1$, $\gamma_i=2$ and $\zeta_i=1$. Thus $\nabla_PF_\delta(P)=4AP$ and $\grad F_\delta(P)=0$, but $u(P)=(I,0)$ and
\begin{equation*}
\operatorname{dist}\bigl((\rho+b_t)(A-I,B),N_\Omega^{\lim}(I,0)\bigr)=(\rho+b_t)\|A-I\|_F=(\rho+b_t)\sqrt2.
\end{equation*}
Here $N_\Omega^{\lim}(I,0)=\{0\}\times\mathbb S^2$, so the Stiefel stopping tests cannot replace the normal-residual test.
\end{remark}

\subsection{Fixed-basis feasible-candidate generation}
Feasible candidates can also be generated using a fixed eigenbasis of either $A$ or $B$. We describe this construction below and derive computable objective-gap bounds in the next subsection.
Let $A=P_A\Lambda_AP_A^\top$ and $B=P_B\Lambda_BP_B^\top$
be orthogonal eigendecompositions, where $\Lambda_A = \operatorname{Diag}(\lambda(A)), \Lambda_B = \operatorname{Diag}(\lambda(B))$.
For each $i=1,\ldots,n$, the scalar branch rule in
Subsection~\ref{subsec:diagonal_decoupling} gives
\[
\begin{aligned}
y_i(P_A)
&=\chi_i(P_A)(\lambda_i(A))_+,
&
z_i(P_A)
&=-(1-\chi_i(P_A))
  \bigl((P_A^\top BP_A)_{ii}\bigr)_-,
\\
y_i(P_B)
&=\chi_i(P_B)
  \bigl((P_B^\top AP_B)_{ii}\bigr)_+,
&
z_i(P_B)
&=-(1-\chi_i(P_B))(\lambda_i(B))_-.
\end{aligned}
\]
Reconstructing the matrices in their respective bases
according to \eqref{eq:reconY}, we obtain
\begin{align}
Y_A&:=P_A\operatorname{Diag}(y(P_A))P_A^\top,\qquad
Z_A:=P_A\operatorname{Diag}(z(P_A))P_A^\top,\qquad
V_A:=V(P_A),\label{eq:fixedA_Y}\\
Y_B&:=P_B\operatorname{Diag}(y(P_B))P_B^\top,\qquad
Z_B:=P_B\operatorname{Diag}(z(P_B))P_B^\top,\qquad
V_B:=V(P_B).\label{eq:fixedB_Y}
\end{align}

Both candidates lie in $\Omega$, so $V_A$ and $V_B$ are upper
bounds on the optimal value of the slack projection problem.
They may be tried before or alongside the Stiefel search,
subject to the same descent and normal-residual tests.

\subsection{Objective-gap bounds for fixed-basis candidates}

Let
\begin{equation}\label{eq:proj_opt_value}
v_\star
:=
\min_{(Y,Z)\in\Omega}
\left\{
\|Y-A\|_F^2+\|Z-B\|_F^2
\right\}
\end{equation}
denote the optimal value of \eqref{eq:projAB}.

Keeping the eigenpair ordering of $A$ fixed, define
\begin{equation}\label{eq:spectral_upper}
U_{\mathrm{spec}}
:=
\max_{\pi}
\sum_{i=1}^n
\max\left\{
(\lambda_i(A))_+^2,
(\lambda_{\pi(i)}(B))_-^2
\right\},
\end{equation}
where $\pi$ ranges over all permutations of $\{1,\ldots,n\}$.
Equivalently, $U_{\mathrm{spec}}$ is the sum of the $n$ largest
of the $2n$ values $(\lambda_i(A))_+^2$ and
$(\lambda_i(B))_-^2$, $i=1,\ldots,n$.
Set
\begin{equation}\label{eq:spectral_lower}
L_{\mathrm{spec}}
:=
\|A\|_F^2+\|B\|_F^2-U_{\mathrm{spec}}.
\end{equation}

To attain the maximum in \eqref{eq:spectral_upper}, pair the values $(\lambda_i(A))_+^2$ and $(\lambda_j(B))_-^2$ in opposite sorted orders. Keep the eigenpair ordering of $A$ fixed and reorder the eigenvalues of $B$ and the corresponding columns of $P_B$ accordingly. Relabelling the reordered eigenvalues gives
\begin{equation}\label{eq:optimal_spectral_pairing}
U_{\mathrm{spec}}
=
\sum_{i=1}^n
\max\left\{
(\lambda_i(A))_+^2,
(\lambda_i(B))_-^2
\right\}.
\end{equation}
Using the resulting eigenvector matrices, define
\begin{equation*}
T:=(P_A^\top P_B)\circ(P_A^\top P_B),
\end{equation*}
where $\circ$ denotes the Hadamard product. The following bound relates the fixed-$A$ objective gap to the alignment of the paired eigenbases.

\begin{proposition}[Objective-gap bound for the fixed-\(A\) candidate]
\label{prop:fixedA_gap}
Under the pairing \eqref{eq:optimal_spectral_pairing},
with $R_0=\max\{\|A\|_2,\|B\|_2\}$, we have
\begin{equation}\label{eq:fixedA_optimality_bound}
0
\leq V_A-v_\star
\leq V_A-L_{\mathrm{spec}}
\leq
2\sqrt n\,R_0\,\|I-T\|_2\,\|B\|_F.
\end{equation}
\end{proposition}

The proof, including the pairing argument, is given in
Appendix~\ref{app:fixed-basis-proof}.

By symmetry, the fixed-$B$ candidate satisfies
\begin{equation*}
0
\leq V_B-v_\star
\leq V_B-L_{\mathrm{spec}}
\leq
2\sqrt n\,R_0\,\|I-T\|_2\,\|A\|_F,
\end{equation*}
since $\|I-T^\top\|_2=\|I-T\|_2$. 

The factor $\|I-T\|_2$ depends on the chosen eigenbases and
the pairing \eqref{eq:optimal_spectral_pairing}, including
choices within repeated eigenspaces. $AB=BA$
guarantees neither $T=I$ under this pairing nor
$v_\star=L_{\mathrm{spec}}$.
Under the same pairing condition, $T=I$ implies
$V_A=V_B=v_\star=L_{\mathrm{spec}}$.

Thus, with $V_{\min}:=\min\{V_A,V_B\}$,
\begin{equation*}
0\leq V_{\min}-v_\star
\leq V_{\min}-L_{\mathrm{spec}}.
\end{equation*}
$V_{\min}-L_{\mathrm{spec}}$ is a computable upper bound on the objective gap of the better fixed-basis candidate. A zero value certifies global optimality, whereas a positive value leaves optimality undecided. Acceptance is evaluated separately through the descent and normal-residual tests.

\subsection{Limiting-normal certificates and slack verification}
\label{subsec:slack_verification}

We now turn to the normal-residual condition required for
slack acceptance.
For $u=(Y,Z)\in\Omega$ and a test vector
$\xi=(\xi^Y,\xi^Z)\in\mathcal U$, we construct computable bounds
on $\operatorname{dist}(\xi,N_\Omega^{\lim}(u))$.
Throughout this subsection, $u\in\Omega$ and $\alpha,\beta,\gamma$ are the exact positive, zero, and negative spectral index sets of $Y+Z$. We use the joint spectral representation from Subsection~\ref{subsec:spectral_geometry} to express the test vector in spectral coordinates:
\begin{equation*}
\widetilde\xi^Y:=P^\top\xi^YP,
\qquad
\widetilde\xi^Z:=P^\top\xi^ZP,
\qquad
A_\beta:=\widetilde\xi^Y_{\beta\beta},
\quad
B_\beta:=\widetilde\xi^Z_{\beta\beta}.
\end{equation*}

For a symmetric matrix $C=U\Diag(c)U^\top$, define its spectral positive and negative parts by
$C_\pm:=U\Diag(c_\pm)U^\top$. We construct certificate branches and compute their projections, obtaining both distance bounds and explicit normals for slack verification.

\paragraph{Certificate sets.}
To make limiting-normal certificates easier to compute,
we construct a finite family of branches by imposing additional
conditions on the zero spectral block of $\mathcal K_W(u)$.
We begin with the two pure branches
\begin{equation*}
\begin{aligned}
\mathcal B_Y(u)
&:=\{(\Gamma^Y,\Gamma^Z)\in\mathcal K_W(u):
\widetilde\Gamma^Y_{\beta\beta}=0\},\\
\mathcal B_Z(u)
&:=\{(\Gamma^Y,\Gamma^Z)\in\mathcal K_W(u):
\widetilde\Gamma^Z_{\beta\beta}=0\},
\end{aligned}
\end{equation*}
and let
\begin{equation*}
\mathcal N_{\mathrm{base}}(u)
:=
\mathcal K_S(u)\cup\mathcal B_Y(u)\cup\mathcal B_Z(u).
\end{equation*}

For $\beta\ne\varnothing$, let $\mathcal P(\beta)$ be the set of all ordered partitions $\pi=(I_+,I_0,I_-)$ of $\{1,\ldots,|\beta|\}$, with empty parts allowed.
Choose finite sets
$\mathcal Q(u;\xi)\subseteq\mathbb O^{|\beta|}$ and
$\Theta_{\mathrm{cert}}\subseteq[0,1]$.
The selected bases may depend on both $u$ and $\xi$.
For $\pi\in\mathcal P(\beta)$, $Q\in\mathcal Q(u;\xi)$, and
$\theta\in\Theta_{\mathrm{cert}}$, write
\begin{equation*}
G_Y:=Q^\top\widetilde\Gamma^Y_{\beta\beta}Q,
\qquad
G_Z:=Q^\top\widetilde\Gamma^Z_{\beta\beta}Q,
\end{equation*}
and define
\begin{equation*}
\begin{aligned}
\mathcal B_{\pi,Q,\theta}(u)
:=\bigl\{
(\Gamma^Y,\Gamma^Z)\in\mathcal K_W(u):\;&
(G_Y)_{++}=(G_Y)_{+0}=0,\\
&(G_Y)_{00}\preceq0,\quad(G_Z)_{00}\succeq0,\\
&(G_Z)_{0-}=(G_Z)_{--}=0,\\
&\theta(G_Y)_{+-}+(1-\theta)(G_Z)_{+-}=0
\bigr\}.
\end{aligned}
\end{equation*}
Here the block subscripts refer to the partition $\pi$, and the transposed block conditions follow from symmetry.
Set
\begin{equation*}
\mathcal N_{\mathrm{cert}}(u;\xi)
:=
\mathcal N_{\mathrm{base}}(u)
\cup
\bigcup_{\substack{
\pi\in\mathcal P(\beta),\ Q\in\mathcal Q(u;\xi)\\
\theta\in\Theta_{\mathrm{cert}}}}
\mathcal B_{\pi,Q,\theta}(u).
\end{equation*}

If $\beta=\varnothing$ or no mixed branches are selected, set $\mathcal N_{\mathrm{cert}}(u;\xi):=\mathcal N_{\mathrm{base}}(u)$.

\begin{proposition}[Limiting-normal certificates]
\label{prop:normal_certificates}
For every $u\in\Omega$ and $\xi\in\mathcal U$,
\begin{equation}\label{eq:certificate_inclusions}
\mathcal N_{\mathrm{cert}}(u;\xi)
\subseteq N_\Omega^{\lim}(u)
=\mathcal K_M(u)
\subseteq\mathcal K_W(u).
\end{equation}
Consequently,
\begin{equation}\label{eq:normal_distance_bounds}
\operatorname{dist}(\xi,\mathcal K_W(u))
\leq\operatorname{dist}(\xi,N_\Omega^{\lim}(u))
\leq\operatorname{dist}
\bigl(\xi,\mathcal N_{\mathrm{cert}}(u;\xi)\bigr).
\end{equation}
If $\beta=\varnothing$, all sets in
\eqref{eq:certificate_inclusions} coincide, and both inequalities
in \eqref{eq:normal_distance_bounds} are equalities.
\end{proposition}

\begin{proof}
We use the spectral characterization in \cite[Proposition~3.3 and Theorem~3.1]{DingSunYe2014}, recalled in Appendix~\ref{app:KM-auxiliary}. Fix $\pi=(I_+,I_0,I_-)$ and $Q$. Absorb a permutation of the partition into $Q$ so that the positive, zero, and negative coordinates occur in this order. For $s\downarrow0$, assign a common value $a_+(s)$ to the coordinates in $I_+$, zero to those in $I_0$, and $a_-(s)$ to those in $I_-$. Choose
\begin{equation*}
\begin{array}{c|ccc}
 & a_+(s) & a_-(s) & a_+(s)/(a_+(s)-a_-(s))\\ \hline
0<\theta<1 & s\theta & -s(1-\theta) & \theta\\
\theta=0 & s^2 & -s & s/(1+s)\\
\theta=1 & s & -s^2 & 1/(1+s)
\end{array}
\end{equation*}
The resulting vectors $a(s)$ are ordered and converge to zero. Their divided-difference matrices converge to an element of $\mathbb U_{|\beta|}$ whose $+-$ block is the constant $\theta$. The other blocks have exactly the values specified in Appendix~\ref{app:KM-auxiliary}. Hence the degenerate-block conditions defining $\mathcal B_{\pi,Q,\theta}(u)$ satisfy that characterization, including the two endpoint values of $\theta$.

The all-positive, all-negative, and all-zero partitions give $\mathcal B_Y(u)$, $\mathcal B_Z(u)$, and $\mathcal K_S(u)$, respectively. Empty parts impose no conditions. This proves the first inclusion in \eqref{eq:certificate_inclusions}; the equality and the last inclusion follow from the same characterization. Distances reverse set inclusion, giving \eqref{eq:normal_distance_bounds}. When $\beta=\varnothing$, no degenerate-block restrictions remain, so all the sets coincide.
\end{proof}

For fixed $(u,\xi)$ and selected bases, each branch is a nonempty closed convex cone. Comparing the branch projections therefore yields a closest normal in $\mathcal N_{\mathrm{cert}}(u;\xi)$. The number of branches is at most $3+|\mathcal Q(u;\xi)|\,|\Theta_{\mathrm{cert}}|\,3^{|\beta|}$.

\paragraph{Explicit distance formulas.}
To evaluate the resulting certificate distances, we first compute
the contribution from the common linear block constraints.
Write $d_W:=\operatorname{dist}(\xi,\mathcal K_W(u))$.
These constraints give
\begin{equation}\label{eq:normal_W_distance}
\begin{aligned}
d_W^2
&=\|\widetilde\xi^Y_{\alpha\alpha}\|_F^2
+2\|\widetilde\xi^Y_{\alpha\beta}\|_F^2+\|\widetilde\xi^Z_{\gamma\gamma}\|_F^2
+2\|\widetilde\xi^Z_{\beta\gamma}\|_F^2\\
&\quad+2\sum_{\substack{i\in\alpha\\j\in\gamma}}
\frac{
\bigl(
\Sigma_{ij}\widetilde\xi^Y_{ij}
+(1-\Sigma_{ij})\widetilde\xi^Z_{ij}
\bigr)^2
}{
\Sigma_{ij}^2+(1-\Sigma_{ij})^2
}.
\end{aligned}
\end{equation}

The basic certificate distance is
\begin{equation}\label{eq:base_certificate_distance}
\begin{aligned}
&\operatorname{dist}^2
\bigl(\xi,\mathcal N_{\mathrm{base}}(u)\bigr)\\
&\quad=d_W^2+\min\bigl\{
\|(A_\beta)_+\|_F^2+\|(B_\beta)_-\|_F^2,
\ \|A_\beta\|_F^2,\ \|B_\beta\|_F^2
\bigr\}.
\end{aligned}
\end{equation}

For a mixed branch, set
$\widehat A:=Q^\top A_\beta Q$ and
$\widehat B:=Q^\top B_\beta Q$.
Then
\begin{equation}\label{eq:mixed_certificate_distance}
\begin{aligned}
&\operatorname{dist}^2
\bigl(\xi,\mathcal B_{\pi,Q,\theta}(u)\bigr)\\
&\quad=d_W^2
+\|\widehat A_{++}\|_F^2
+2\|\widehat A_{+0}\|_F^2
+\|(\widehat A_{00})_+\|_F^2\\
&\qquad+\|\widehat B_{--}\|_F^2
+2\|\widehat B_{0-}\|_F^2
+\|(\widehat B_{00})_-\|_F^2\\
&\qquad+
\frac{
2\|\theta\widehat A_{+-}
+(1-\theta)\widehat B_{+-}\|_F^2
}{
\theta^2+(1-\theta)^2
}.
\end{aligned}
\end{equation}
The squared distance to $\mathcal N_{\mathrm{cert}}(u;\xi)$ is the minimum of the squared distances to its included branches.

These formulas follow from blockwise orthogonal projections,
with the factors $2$ accounting for transposed off-diagonal blocks in the Frobenius norm. The same projections also provide an explicit normal witness for the slack verification below. For a pair of scalar entries or equally sized blocks $(U_0,V_0)$ subject to $aU+bV=0$, with $a^2+b^2>0$, set
\begin{equation}\label{eq:certificate_cross_projection}
E:=aU_0+bV_0,
\qquad
U_*:=U_0-\frac{aE}{a^2+b^2},
\qquad
V_*:=V_0-\frac{bE}{a^2+b^2}.
\end{equation}
For the common $\alpha\gamma$ block, apply this formula entrywise with $(a,b)=(\Sigma_{ij},1-\Sigma_{ij})$. For a mixed $+-$ block, apply it with $(a,b)=(\theta,1-\theta)$ to $(\widehat A_{+-},\widehat B_{+-})$, and fill the transpose blocks by symmetry. Set all blocks required to vanish equal to zero and retain all free blocks.

For the basic branches, the projected $\beta\beta$ pair is $(-(A_\beta)_-,(B_\beta)_+)$, $(0,B_\beta)$, or $(A_\beta,0)$ for $\mathcal K_S$, $\mathcal B_Y$, or $\mathcal B_Z$, respectively. For a mixed branch, use
\begin{equation*}
(G_Y)_{00}=-(\widehat A_{00})_- ,
\qquad
(G_Z)_{00}=(\widehat B_{00})_+.
\end{equation*}
Rotate the mixed degenerate block back by $Q$, then rotate the full pair back by $P$. These orthogonal transformations preserve the product Frobenius norm. Return the minimizing normal $\Gamma=(\Gamma^Y,\Gamma^Z)$, its branch type, the spectral data $P$ and $(\alpha,\beta,\gamma)$, and, where applicable, $Q,\pi,\theta$, together with $r=\|\xi-\Gamma\|$. The returned branch data permit an independent check of membership, while \eqref{eq:base_certificate_distance} or \eqref{eq:mixed_certificate_distance} gives the squared residual.

\paragraph{Slack verification.}
To use these certificates in the slack update, consider a feasible
trial point $u^{t+1}\in\Omega$ and set
\begin{equation*}
\xi^{t+1}
:=-\nabla_u\Phi^t(u^{t+1})
=(\rho+b_t)(\widehat c^t-u^{t+1}).
\end{equation*}

The computable test is
\begin{equation}\label{eq:slack_certificate_test}
\operatorname{dist}\bigl(
\xi^{t+1},
\mathcal N_{\mathrm{cert}}(u^{t+1};\xi^{t+1})
\bigr)
\leq\delta_u^t.
\end{equation}
Let $\Gamma^{t+1}$ be a closest certificate normal returned by the preceding block projections, attaining the distance in this test. When the test passes, Proposition~\ref{prop:normal_certificates} gives
\[
\Gamma^{t+1}\in N_\Omega^{\lim}(u^{t+1}),
\qquad
\|\nabla_u\Phi^t(u^{t+1})+\Gamma^{t+1}\|
\leq\delta_u^t.
\]
Thus the slack residual condition \eqref{eq:ures} holds, and the candidate is accepted if it also satisfies the descent condition \eqref{eq:udescent}. Retain this normal and its branch data for evaluating the full inner residual in Lemma~\ref{lem:block_full_residual}. At inner termination, the retained normal may also be used in the outer slack test whenever $\|\nabla_u L_k(x^{k+1})+\Gamma^{k+1}\|\leq\varepsilon_{\mathrm{stat}}$.

The candidate-generation routines and certificate choices used in the numerical experiments will be specified in Section~\ref{sec:numerical_experiments}.

\section{Numerical experiments}\label{sec:numerical_experiments}
We evaluate SPALM on two test families, denoted by Families~$\mathrm{I}$ and~$\mathrm{II}$. For the runs reported in Tables~1 and~2, we use two relative noise levels, $\eta_{\mathrm{rel}}\in\{0.01,0.05\}$, and five seeds for each parameter group. Both families fit the SDCMPCC model \eqref{P} and its slack reformulation \eqref{RE}. We describe the implementation, evaluation criteria, models, and data generation below.

\subsection{Implementation and evaluation criteria}
\label{subsec:numerical_protocol}
We evaluate finite-budget runs using the final complete outer tuple, its updated multipliers, and normal witnesses constructed from retained joint spectral data. The formulas in Section~\ref{subsec:slack_verification} define the certificate calculations; Remark~\ref{rem:implementation_theory} specifies their interpretation in floating-point arithmetic. All matrix blocks use the Frobenius metric and all other blocks the Euclidean metric, with the canonical product norm. An isometric vectorization of symmetric matrices weights off-diagonal upper-triangular entries by $\sqrt2$.

\paragraph{Reported residuals.}
At a complete tuple $(w,Y,Z,\lambda,\mu,M^G,M^H)$, define
\begin{equation*}
\begin{aligned}
r_{\mathrm{pri}}:=\max\{&\|g(w)_+\|_2,\|h(w)\|,\|G(w)-Y\|_F,\|H(w)-Z\|_F,\\
&\operatorname{dist}(Y,\mathbb S_+^n),\operatorname{dist}(Z,\mathbb S_-^n),\|YZ\|_F\},\\
r_\lambda:=&\frac{\|(-\lambda)_+\|_2}{1+\|\lambda\|_2},\qquad
r_{\lambda g}:=\frac{\|\lambda\circ g(w)\|_2}{1+\|\lambda\|_2+\|g(w)\|_2},\\
s_w:=&\nabla f(w)+g'(w)^*\lambda+h'(w)^*\mu+G'(w)^*M^G+H'(w)^*M^H,\\
r_{\mathrm{stat}}:=&\frac{\|s_w\|}{1+\|\nabla f(w)\|+\|\lambda\|_2+\|\mu\|+\|M^G\|_F+\|M^H\|_F}.
\end{aligned}
\end{equation*}
The inequality residuals are zero for Family~$\mathrm{II}$. The primal residual vanishes exactly at feasible points, while $V_k$ in \eqref{eq:feasibility_measure} instead controls the penalty update. Model-specific expressions for $s_w$ are given below.

To evaluate the slack-stationarity component, we use the retained
joint spectral representation of $u=(Y,Z)$.
Following Section~\ref{subsec:slack_verification}, the signs and structural zeros of the
retained coefficients determine the spectral partition used
for normal-witness evaluation.
With $\xi=(M^G,M^H)$, let $A_\beta$ and $B_\beta$ denote
its zero-block spectral coordinates in this representation.

The selected certificate family contains the three basic branches.
For $1\leq|\beta|\leq6$, it also includes mixed branches constructed
using
\begin{equation*}
\mathcal Q=
\{I,Q(A_\beta),Q(B_\beta),Q(A_\beta+B_\beta),Q(A_\beta-B_\beta)\},
\qquad
\Theta_{\mathrm{cert}}=\{0,1\},
\end{equation*}
where $Q(C)$ is a selected orthogonal eigenbasis and duplicate bases
are discarded. All ordered partitions of the zero block are tested.
For $|\beta|=0$ or $|\beta|>6$, the three basic branches are used.

At the same slack point $u=(Y,Z)$, the reporting witness
$\Gamma_{\mathrm{rep}}$ minimizes $\|\xi-\Gamma\|$
over the retained inner witness and the projections onto
the selected certificate branches.
Define the normalization factor and the reporting residual by
\begin{equation*}
s_N:=1+\sqrt{\|M^G\|_F^2+\|M^H\|_F^2},
\qquad
r_{M,\Gamma}:=\frac{\|\xi-\Gamma_{\mathrm{rep}}\|}{s_N},
\end{equation*}
and
\begin{equation*}
R_{M,\Gamma}:=
\max\{r_{\mathrm{pri}},r_\lambda,r_{\lambda g},
      r_{\mathrm{stat}},r_{M,\Gamma}\}.
\end{equation*}

A run is counted as successful when the absolute feasibility and
stationarity tests associated with Algorithm~\ref{alg:sALcc},
the reporting test $R_{M,\Gamma}\leq10^{-4}$, and the independent
numerical checks are satisfied. The absolute feasibility and
stationarity tolerances are both $10^{-4}$. The absolute slack test in Algorithm~\ref{alg:sALcc} uses the retained inner normal,
whereas $\Gamma_{\mathrm{rep}}$ is used for the reporting residual.

Tables report the final complete outer tuple and its updated multipliers. Here $N$ is the number of planned runs and $M$ is the number of successful runs. Each numerical column reports the median of its per-run values; $R_{M,\Gamma}$ is evaluated for each tuple
before aggregation. Solver wall-clock time includes candidate generation, witness evaluation, and checkpoint storage. Independent-verification time is recorded separately.

\paragraph{Inner schedule and finite budgets.}
At outer iteration $k=0,1,\ldots$, we choose the inner initial point $x^{k,0}$ according to \eqref{eq:reference_initialization}, using the fixed feasible reference specified below for each family. We set
\begin{equation*}
a_t=b_t=10^{-2},\qquad
\nu_{k+1}=\max\{10^{-6},10^{-1}(0.8)^{k+1}\},\qquad
\delta_w^t=\delta_u^t=\nu_{k+1}/4.
\end{equation*}
Each accepted block update is checked against the proximal descent condition and the corresponding absolute error test. The retained slack normal is used in the full residual evaluation. In both families, $G$ extracts one matrix block and $H$ extracts the negative of another, so the inner return criterion is
\begin{equation*}
\widehat{\mathcal R}_{\rho_k}^{\,t+1}
=
\left[
\bigl(\rho_k\|\Delta u^t\|
+a_t\|\Delta w^t\|+\delta_w^t\bigr)^2
+
\bigl(b_t\|\Delta u^t\|+\delta_u^t\bigr)^2
\right]^{1/2}
\leq\nu_{k+1}.
\end{equation*}
The outer and inner iteration limit is $400$. Both families use $\rho_0=0.2$ and a feasibility reduction factor of $0.8$. The penalty growth factors are $1.5$ and $1.15$, respectively. A requested penalty increase exceeding $30$ terminates the run. Iteration, search, and time budgets have separately recorded termination statuses.

L-BFGS and conjugate gradients use iteration limits of $160$ and
$500$, respectively, with a CG relative tolerance of $10^{-13}$.
The Huber--Stiefel search uses smoothing parameter $10^{-3}$,
gradient tolerance $10^{-10}$, and an iteration limit of $480$.
Fixed-support polishing uses at most two initial bases, an iteration
limit of $240$, and gradient tolerance $10^{-10}$.
Three seeds are executed concurrently with one BLAS thread per process.

\paragraph{Slack-candidate generation.}
Slack candidates are generated by fixed-basis projections,
Huber--Stiefel searches, and fixed-support polishing.
For a proximal center $(A,B)$, let $P=[P_+,P_-]\in\On$ have a fixed column partition. Compressed spectral projections give
\begin{equation*}
Y=P_+\Pi_{\mathbb S_+}(P_+^\top AP_+)P_+^\top,\qquad
Z=P_-\Pi_{\mathbb S_-}(P_-^\top BP_-)P_-^\top,
\end{equation*}
where the projections act on symmetric matrices of the appropriate
dimensions. The compressed eigendecompositions produce the joint
spectral basis retained for normal-witness evaluation.

To improve the candidate while keeping the column partition fixed, fixed-support polishing varies $P$ to maximize
\begin{equation*}
F_{\mathrm{sup}}(P)=
\|\Pi_{\mathbb S_+}(P_+^\top AP_+)\|_F^2+
\|\Pi_{\mathbb S_-}(P_-^\top BP_-)\|_F^2.
\end{equation*}
Each generated candidate is evaluated by the proximal descent
and absolute normal-residual tests.

\begin{remark}[Finite-precision implementation]
\label{rem:implementation_theory}
The implementation uses finite search budgets, an accuracy floor $\nu_{\min}=10^{-6}$, and an absolute descent-comparison tolerance of $10^{-11}$. Structural spectral representations, normal witnesses, and residuals are checked in floating-point arithmetic. The reported verification is at the numerical-witness level; rigorous roundoff certification requires validated error enclosures. The selected certificate family provides a sufficient test for normal-residual accuracy. The asymptotic results in Sections~\ref{sec:ALM}--\ref{sec:inner_solver} concern infinite sequences satisfying the stated vanishing-error assumptions.
\end{remark}

\subsection{Family \texorpdfstring{$\mathrm{I}$}{I}:
a compact low-rank correlation formulation}
\label{subsec:family-1-prime}

Numerical methods for nearest low-rank correlation matrix problems include \cite{LiQi2011,ZhangWu2003,BiHanPan2013}. Given an observed symmetric matrix $C_{\mathrm{obs}}\in\mathbb S^n$
and an integer $1\leq r\leq n$, we consider

\begin{equation}\label{eq:num_familyIprime_rank_problem}
\begin{aligned}
 \min_{X\in\Sn}\quad
 & \frac12\norm{X-C_{\mathrm{obs}}}_F^2\\
 \mathrm{s.t.}\quad
 & \diag(X)=\mathbf 1,\qquad X\succeq0,\qquad
   \operatorname{rank}(X)\leq r.
\end{aligned}
\end{equation}
The representation in \cite{DingSunYe2014} uses a trace equality and an additional semidefinite constraint on an auxiliary matrix. For a fixed lifting parameter $\vartheta\in(0,1)$, we use the following compact two-matrix SDCMPCC formulation:
\begin{equation}\label{eq:num_familyIprime_model}
\begin{aligned}
 \min_{X,S\in\Sn}\quad
 & \frac12\norm{X-C_{\mathrm{obs}}}_F^2\\
 \mathrm{s.t.}\quad
 & \diag(X)=\mathbf 1,\\
 & \norm{I-S}_F^2\leq r+\vartheta,\\
 & X\succeq0,\qquad S\succeq0,\qquad XS=0.
\end{aligned}
\end{equation}
Set
\begin{equation}\label{eq:num_familyIprime_GH}
w=(X,S),\quad g(w)=\|I-S\|_F^2-r-\vartheta,\quad h(w)=\diag(X)-\mathbf1,\quad G(w)=X,\ H(w)=-S.
\end{equation}

For this model, the reporting stationarity vector
$s_w=(s_X,s_S)$ has the explicit blocks
\begin{equation*}
\begin{aligned}
 s_X&=X-C_{\mathrm{obs}}+\Diag(\mu)+M^G,\\
 s_S&=2\lambda(S-I)-M^H,\qquad
 \|s_w\|^2=\|s_X\|_F^2+\|s_S\|_F^2.
\end{aligned}
\end{equation*}

The projection of the feasible set of \eqref{eq:num_familyIprime_model} onto
the $X$-space is exactly the feasible set of the rank-constrained correlation
problem \eqref{eq:num_familyIprime_rank_problem}.  Indeed, if
$k=\operatorname{rank}(X)$, then $XS=0$ implies
$\operatorname{range}(X)\subseteq\ker S$.  Hence $S$ has at least $k$ zero
eigenvalues and
\begin{equation*}
 k\leq\norm{I-S}_F^2\leq r+\vartheta<r+1.
\end{equation*}
Since $k$ is an integer, this gives $k\leq r$.  Conversely, for any
$X\succeq0$ with $\operatorname{rank}(X)=k\leq r$, let $P_X$ be the
orthogonal projector onto $\operatorname{range}(X)$ and set $S=I-P_X$.
Then $S\succeq0$, $XS=0$, and $\norm{I-S}_F^2=k\leq r < r + \vartheta$.

The lifted feasible set is also compact. Indeed, every feasible pair
satisfies
\begin{equation*}
 \|X\|_F\leq\operatorname{tr}(X)=n,\qquad
 \|S\|_F\leq\sqrt n+\sqrt{r+\vartheta}.
\end{equation*}
The constraints are closed, and the set is nonempty, as shown by $X=\mathbf1\mathbf1^\top$ and $S=I-X/n$. The continuous objective therefore attains its minimum.

\begin{remark}[Scope of the lifting]\label{rem:familyIprime_lifting_scope}
The two formulations have the same global optimal value, and the projection $(X,S)\mapsto X$ maps the global minimizers of \eqref{eq:num_familyIprime_model} onto those of \eqref{eq:num_familyIprime_rank_problem}.
For Family~$\mathrm{I}$, our M-stationarity results and numerical stopping criteria refer to the lifted formulation \eqref{eq:num_familyIprime_model}.
\end{remark}

\paragraph{Constraint qualification and inner level boundedness.}
\begingroup\emergencystretch=3em
We next verify the constraint qualification used in the outer
analysis and the level-boundedness assumption for the inner subproblems.
\begin{proposition}[Active-set SDCMPCC-LICQ for Family~$\mathrm{I}$]
\label{prop:familyIprime_LICQ}
The slack reformulation induced by \eqref{eq:num_familyIprime_model}
satisfies the active-set SDCMPCC-LICQ of
Definition~\ref{def:SDCMPCC_LICQ} at every feasible point. Equivalently, the
unique slack lift $((X,S),(X,-S))$ of every feasible $(X,S)$ satisfies this
condition. This conclusion does not require strict complementarity or
$\operatorname{rank}(X)=r$.
\end{proposition}
The proof is given in Appendix~\ref{app:familyIprime-assumptions}. The threshold $r+\vartheta\in(r,r+1)$ enforces the integer rank bound. At a feasible pair with $S$ an orthogonal projector, $\norm{I-S}_F^2=n-\operatorname{rank}(S)$ is an integer, so the scalar inequality is strict.

\par\endgroup

\begin{corollary}[Absence of singular M-multipliers for Family~$\mathrm{I}$]
\label{cor:familyI_no_singularM}
For the constraint representation \eqref{eq:num_familyIprime_GH} of
Family~$\mathrm{I}$, no feasible point of
\eqref{eq:num_familyIprime_model} admits a singular M-multiplier.
\end{corollary}
The proof is given in Appendix~\ref{app:familyIprime-assumptions}.

\begin{remark}[Dependence on the constraint representation]
\label{rem:familyI_formulation_comparison}
The formulation in \cite{DingSunYe2014} uses an auxiliary matrix $U$ satisfying
\begin{equation*}
 \operatorname{tr}(U)=r,\qquad U\succeq0,\qquad (X,U-I)\in\Omega,
\end{equation*}
in addition to $\operatorname{diag}(X)=\mathbf1$. If $\operatorname{rank}(X)=r$, complementarity implies that $U$ acts as the identity on $\operatorname{range}(X)$. Positive semidefiniteness and the trace equality then force $U=P_X$, where $P_X$ is the orthogonal projector onto $\operatorname{range}(X)$.
At this point, the pair $(X,U-I)$ is strictly complementary. Set the trace multiplier to $\tau=1$, the diagonal-equality multiplier to zero, and
\begin{equation*}
 W=-(I-P_X),\qquad (\Gamma^G,\Gamma^H)=(0,-P_X).
\end{equation*}
Then $W\in N_{\Sp}(U)$ and
\begin{equation*}
 (\Gamma^G,\Gamma^H)\in\Nhat_\Omega(X,U-I)
 \subseteq\Nlim_\Omega(X,U-I),
\end{equation*}
as follows from the normal-cone formula in \cite[Proposition~3.2]{DingSunYe2014}. Moreover,
\begin{equation*}
 \Gamma^G=0,\qquad \tau I+W+\Gamma^H=0.
\end{equation*}
Thus the earlier constraint representation admits a nonzero abnormal multiplier at every feasible point with $\operatorname{rank}(X)=r$. In contrast, Corollary~\ref{cor:familyI_no_singularM} excludes singular M-multipliers at every feasible point of \eqref{eq:num_familyIprime_model}.
\end{remark}

\begin{proposition}[Inner level boundedness for Family~$\mathrm{I}$]
\label{prop:familyIprime_level_bounded}
Fix $\rho>0$ and any safeguarded multiplier vector $\Upsilon$. For every
$x^0$ in the inner domain, the sublevel set $\mathcal S_\rho(x^0)$ in
\eqref{eq:inner_sublevel} is bounded.
\end{proposition}
The proof is given in Appendix~\ref{app:familyIprime-assumptions}.

\paragraph{A feasible reference point for Family $\mathrm{I}$.}
For $e=\mathbf1\in\mathbb R^n$, set
\begin{equation}\label{eq:familyI_reference}
X^{\mathrm{ref}}=ee^\top,\qquad
S^{\mathrm{ref}}=I-\frac{ee^\top}{n},\qquad
x^{\mathrm{ref}}=((X^{\mathrm{ref}},S^{\mathrm{ref}}),(X^{\mathrm{ref}},-S^{\mathrm{ref}})).
\end{equation}
This point is feasible: $\diag X^{\mathrm{ref}}=e$, $X^{\mathrm{ref}}S^{\mathrm{ref}}=0$, both matrices are positive semidefinite, and $\|I-S^{\mathrm{ref}}\|_F^2=1<r+\vartheta$ for $r\geq1$.
Lemma~\ref{lem:reference_AL_bound} therefore applies with
$f(w^{\mathrm{ref}})=\tfrac12\|ee^\top-C_{\mathrm{obs}}\|_F^2$.
For the inequality term, $g(w^{\mathrm{ref}})=1-r-\vartheta<0$ and $\bar\lambda^k\geq0$ give
$(g(w^{\mathrm{ref}})+\bar\lambda^k/\rho_k)_+\leq\bar\lambda^k/\rho_k$.
The reference point is determined by the dimension $n$.

\paragraph{Planted instances and exact relative noise.}
We first construct a planted feasible pair $(X^\star,S^\star)$ with prescribed ranks and common kernel dimension, then perturb $X^\star$ to obtain the observations.
Let $n\geq2$ and choose integers $r,d_0$ satisfying
$1\leq r\leq n$ and $0\leq d_0\leq r-1$.
Set $r_0:=r-d_0$, so $1\leq r_0\leq r$.
Generate $\widetilde F\in\mathbb R^{n\times r_0}$ with independent
standard normal entries and normalize every row,
\begin{equation*}
 F_{i:}:=\frac{\widetilde F_{i:}}{\norm{\widetilde F_{i:}}_2},
 \qquad X^\star:=FF^\top.
\end{equation*}
With probability one, the normalization is well defined and
preserves the full column rank of $\widetilde F$. Hence
$X^\star\succeq0$, $\diag(X^\star)=\mathbf 1$, and
$\operatorname{rank}(X^\star)=r_0$.
Let $Q_X\in\mathbb R^{n\times r_0}$ have orthonormal columns
spanning $\operatorname{range}(X^\star)$, and choose
$Q_0\in\mathbb R^{n\times d_0}$ with orthonormal columns in
$\ker X^\star$. Define
\begin{equation*}
 S^\star:=I-Q_XQ_X^\top-Q_0Q_0^\top.
\end{equation*}
Consequently,
\begin{equation*}
 \operatorname{rank}(S^\star)=n-r,\qquad
 X^\star S^\star=0,\qquad
 \norm{I-S^\star}_F^2=r,
\end{equation*}
and $\dim(\ker X^\star\cap\ker S^\star)=d_0$.
The scalar inequality is strictly satisfied, with
$g(X^\star,S^\star)=-\vartheta$.
The parameter $d_0$ jointly specifies the planted common kernel dimension and rank $r_0=r-d_0$. The spectral partition at termination is determined from the retained joint spectral representation.

To obtain observations at a prescribed relative noise level,
draw $\widetilde E\in\mathbb R^{n\times n}$ with independent
standard normal entries and set
\begin{equation}\label{eq:num_familyIprime_noise}
\begin{aligned}
 E_0
 &:=\frac{\widetilde E+\widetilde E^\top}{2}
 -\Diag\!\left(\diag\!\left(
 \frac{\widetilde E+\widetilde E^\top}{2}\right)\right),\\
 E&:=\frac{E_0}{\norm{E_0}_F},\\
 C_{\mathrm{obs}}
 &:=X^\star+\eta_{\mathrm{rel}}\norm{X^\star}_F E,
 \qquad \eta_{\mathrm{rel}}\in\{0.01,0.05\}.
\end{aligned}
\end{equation}
Since $n\geq2$, $E_0$ is nonzero with probability one, so its
normalization is well defined. It follows that
\begin{equation*}
 \diag(C_{\mathrm{obs}})=\mathbf 1,
 \qquad
 \frac{\norm{C_{\mathrm{obs}}-X^\star}_F}
      {\norm{X^\star}_F}
 =\eta_{\mathrm{rel}}.
\end{equation*}
For fixed $(n,r,d_0,\text{seed})$, the same normalized direction
$E$ is used at both relative noise levels.

All Family~$\mathrm{I}$ experiments use $\vartheta=1/2$. The test grid is
\begin{equation*}
 (n,r)\in\{(20,5),(30,8),(50,10),(100,20)\},\quad
 d_0\in\{0,2\},\quad
 \eta_{\mathrm{rel}}\in\{0.01,0.05\},
\end{equation*}
with five seeds for every parameter group.  The experiments use the implementation settings in Section~\ref{subsec:numerical_protocol}.
The $w$-block is treated by L-BFGS with an analytic blockwise fallback.
The per-run time limits are $5400$ seconds. Time limits are checked at algorithmic
checkpoints.

\begin{table}[htbp]
\centering
\caption{Finite-budget results for Family~$\mathrm{I}$.
Residuals and solver times are medians over five runs per
parameter group; $M/N$ denotes successful runs over planned runs. Time is the solver wall-clock time, including candidate generation, witness evaluation, and checkpoint storage.}
\label{tab:family-1-prime-results}
\normalsize
\setlength{\tabcolsep}{5pt}
\renewcommand{\arraystretch}{1.12}
\begin{tabular}{rrrrrrrrrr}
\toprule
$d_0$ & $n$ & $r$ & $\eta_{\mathrm{rel}}$ & $M/N$ &
$R_{M,\Gamma}$ & $r_{\mathrm{stat}}$ & $r_{M,\Gamma}$ &
$r_{\mathrm{pri}}$ & Time (s) \\
\midrule
0 & 20 & 5 & 0.01 & 5/5 & $6.74\!\times\!10^{-5}$ & $6.74\!\times\!10^{-5}$ & $2.75\!\times\!10^{-6}$ & $2.84\!\times\!10^{-5}$ & 1.18 \\
0 & 20 & 5 & 0.05 & 5/5 & $8.72\!\times\!10^{-5}$ & $1.27\!\times\!10^{-5}$ & $7.03\!\times\!10^{-7}$ & $8.72\!\times\!10^{-5}$ & 2.31 \\
0 & 30 & 8 & 0.01 & 5/5 & $7.11\!\times\!10^{-5}$ & $7.11\!\times\!10^{-5}$ & $4.83\!\times\!10^{-6}$ & $4.21\!\times\!10^{-5}$ & 1.02 \\
0 & 30 & 8 & 0.05 & 5/5 & $9.21\!\times\!10^{-5}$ & $2.39\!\times\!10^{-5}$ & $1.60\!\times\!10^{-6}$ & $9.21\!\times\!10^{-5}$ & 2.40 \\
0 & 50 & 10 & 0.01 & 5/5 & $7.05\!\times\!10^{-5}$ & $7.05\!\times\!10^{-5}$ & $3.36\!\times\!10^{-6}$ & $4.76\!\times\!10^{-5}$ & 1.83 \\
0 & 50 & 10 & 0.05 & 5/5 & $9.06\!\times\!10^{-5}$ & $5.92\!\times\!10^{-6}$ & $3.20\!\times\!10^{-7}$ & $9.06\!\times\!10^{-5}$ & 4.65 \\
0 & 100 & 20 & 0.01 & 5/5 & $5.25\!\times\!10^{-5}$ & $5.25\!\times\!10^{-5}$ & $2.92\!\times\!10^{-6}$ & $3.21\!\times\!10^{-5}$ & 5.91 \\
0 & 100 & 20 & 0.05 & 5/5 & $7.96\!\times\!10^{-5}$ & $4.22\!\times\!10^{-6}$ & $3.10\!\times\!10^{-7}$ & $7.96\!\times\!10^{-5}$ & 46.02 \\
2 & 20 & 5 & 0.01 & 5/5 & $7.43\!\times\!10^{-5}$ & $6.77\!\times\!10^{-5}$ & $6.44\!\times\!10^{-7}$ & $7.43\!\times\!10^{-5}$ & 1.74 \\
2 & 20 & 5 & 0.05 & 5/5 & $9.11\!\times\!10^{-5}$ & $7.85\!\times\!10^{-6}$ & $1.04\!\times\!10^{-7}$ & $9.11\!\times\!10^{-5}$ & 3.69 \\
2 & 30 & 8 & 0.01 & 5/5 & $6.78\!\times\!10^{-5}$ & $6.43\!\times\!10^{-5}$ & $1.47\!\times\!10^{-6}$ & $6.08\!\times\!10^{-5}$ & 1.69 \\
2 & 30 & 8 & 0.05 & 5/5 & $9.65\!\times\!10^{-5}$ & $8.69\!\times\!10^{-6}$ & $2.53\!\times\!10^{-7}$ & $9.65\!\times\!10^{-5}$ & 3.41 \\
2 & 50 & 10 & 0.01 & 5/5 & $8.06\!\times\!10^{-5}$ & $5.23\!\times\!10^{-5}$ & $1.15\!\times\!10^{-6}$ & $8.06\!\times\!10^{-5}$ & 2.57 \\
2 & 50 & 10 & 0.05 & 5/5 & $9.02\!\times\!10^{-5}$ & $5.64\!\times\!10^{-6}$ & $1.34\!\times\!10^{-7}$ & $9.02\!\times\!10^{-5}$ & 5.34 \\
2 & 100 & 20 & 0.01 & 5/5 & $8.60\!\times\!10^{-5}$ & $3.24\!\times\!10^{-5}$ & $1.62\!\times\!10^{-6}$ & $8.60\!\times\!10^{-5}$ & 6.90 \\
2 & 100 & 20 & 0.05 & 5/5 & $8.64\!\times\!10^{-5}$ & $4.10\!\times\!10^{-6}$ & $2.05\!\times\!10^{-7}$ & $8.64\!\times\!10^{-5}$ & 59.44 \\
\bottomrule
\end{tabular}
\end{table}

\subsection{Family \texorpdfstring{$\mathrm{II}$}{II}:
inverse linear SDP with a primal observation}
\label{subsec:family-2-prime}

Given an observed cost matrix $\widehat C\in\Sn$ and a primal
observation $X_{\mathrm{obs}}\in\Sn$, we fit a cost matrix
and a corresponding primal optimizer of the forward SDP $\min\{\langle C,U\rangle:\langle A,U\rangle=1,\ U\succeq0\}$. Following the minimum-distance approach to noisy observations \cite{AswaniShenSiddiq2018,ChanMahmoodZhu2025}, the objective penalizes deviations from $\widehat C$ and $X_{\mathrm{obs}}$. This fixes the constraint data, unlike left-hand-side recovery in \cite{GhatraniGhate2024}; related inverse SDCMPCC formulations also adjust right-hand sides \cite{WuZhangZhangLu2016}.

The dual of the forward SDP is
\begin{equation*}
\max_{y\in\mathbb R,\ S\in\Sn}\ y
\quad\text{s.t.}\quad C-yA=S\succeq0,
\end{equation*}
where $y$ is the scalar equality multiplier and $S$ is the dual slack. Write $A=\overline A/b_0$, with $b_0>0$ and $\overline A$ indefinite. The working dual coordinate $\widetilde y=y/b_0$ converts $yA+S=C$ into $C-\widetilde y\overline A-S=0$, giving

\begin{equation}\label{eq:num_familyIIprime_model}
\begin{aligned}
  \min_{C,X,S,\widetilde y}\quad
  & \frac12\norm{C-\widehat C}_F^2
    +\frac12\norm{X-X_{\mathrm{obs}}}_F^2\\
  \mathrm{s.t.}\quad
  & \langle A,X\rangle=1,\\
  & C-\widetilde y\,\overline A-S=0,\\
  & X\succeq0,\qquad S\succeq0,\qquad XS=0.
\end{aligned}
\end{equation}
To express \eqref{eq:num_familyIIprime_model} in the form of
\eqref{P}, set
\begin{equation*}
 w=(C,X,S,\widetilde y),\qquad
 h(w)=\bigl(\langle A,X\rangle-1,\,
 C-\widetilde y\,\overline A-S\bigr).
\end{equation*}
Consistently with the sign convention in \eqref{P}, define
\begin{equation}\label{eq:num_familyIIprime_GH}
 G(w)=X,\qquad H(w)=-S.
\end{equation}

\paragraph{Engineering interpretations.} Variable-ratio actuators allow a trade-off between motor torque and speed \cite{girard2017gear,mooney2013actuator}. For an ideal continuously variable transmission with inactive operating limits, fixed positive output power gives $\tau\omega=1$ in normalized coordinates $z=(\tau,\omega)^\top$. With $X=zz^\top$ and $A=\tfrac12\bigl(\begin{smallmatrix}0&1\\1&0\end{smallmatrix}\bigr)$, the power constraint becomes $\langle A,X\rangle=1$ and the quadratic operating cost is $\langle C,X\rangle=z^\top Cz$. Shear-buckling tests on panels \cite{nicassio2022shear} motivate a second interpretation. For fixed reference shear prestress, the reduced geometric stiffness matrix $A$ can be indefinite \cite{defaria2018buckling}. If the effective elastic stiffness $C$ is positive definite, then at the first positive buckling factor $y$ the tangent stiffness $S=C-yA\succeq0$ satisfies $Sq=0$ for a critical mode $q$ normalized by $q^\top A q=1$. Thus $X=qq^\top$ satisfies $XS=0$. In \eqref{eq:num_familyIIprime_model}, $\widehat C$ is the nominal cost or stiffness matrix, and $X_{\mathrm{obs}}$ is the outer product of an observed torque--speed vector or normalized buckling mode.

For the numerical experiments, we use the scaled slack $T_n=nS$ and the working variable $w=(C,X,T_n,\widetilde y)$. In these coordinates the constraint maps in \eqref{eq:num_familyIIprime_model} are
\begin{equation*}
 h(w)=\bigl(\langle A,X\rangle-1,\ C-\widetilde y\,\overline A-T_n/n\bigr),
 \qquad G(w)=X,\qquad H(w)=-T_n.
\end{equation*}
Writing the equality multiplier as $\mu=(\sigma,\Lambda)\in\mathbb R\times\Sn$, the stationarity blocks used in Section~\ref{subsec:numerical_protocol} are
\begin{equation*}
\begin{aligned}
 s_C&=C-\widehat C+\Lambda, &
 s_X&=X-X_{\mathrm{obs}}+\sigma A+M^G,\\
 s_{T_n}&=-\Lambda/n-M^H, &
 s_{\widetilde y}&=-\langle\overline A,\Lambda\rangle.
\end{aligned}
\end{equation*}
All stopping tests and reported residuals for Family~$\mathrm{II}$ use the product norm in these working coordinates. For fixed $n$, the invertible scaling $(S,Z_S)\mapsto(T_n,Z)=(nS,nZ_S)$ preserves complementarity and multiplies the second linking residual by $n$. The active-set LICQ and AL sublevel-boundedness results below therefore apply in the working coordinates as well.

The forward single-equality SDP is strictly feasible. To see this, choose
a unit eigenvector $v$ of $A$ with eigenvalue $\alpha>0$ and take
$\epsilon>0$ small enough that $1-\epsilon\operatorname{tr}(A)>0$. Then
\begin{equation*}
 U_\epsilon:=\epsilon I+
 \frac{1-\epsilon\operatorname{tr}(A)}{\alpha}vv^\top\succ0,
 \qquad \langle A,U_\epsilon\rangle=1.
\end{equation*}
By Slater's theorem, a finite optimal value of the forward problem is
attained by its dual. Every forward optimal solution can therefore be
paired with $y$ and $S=C-yA\succeq0$ satisfying $XS=0$.
Thus the KKT representation is complete for this family whenever
a forward optimum exists. Strict feasibility alone does not guarantee a finite
forward optimal value for every cost matrix $C$.

\paragraph{Constraint qualification and inner level boundedness.}
\begin{proposition}[Active-set SDCMPCC-LICQ for Family~$\mathrm{II}$]
\label{prop:familyIIprime_LICQ}
The unique slack lift $((C,X,S,\widetilde y),(X,-S))$ of every feasible
point of
\eqref{eq:num_familyIIprime_model} satisfies the active-set SDCMPCC-LICQ of
Definition~\ref{def:SDCMPCC_LICQ}.
\end{proposition}
The proof in Appendix~\ref{app:familyIIprime-assumptions} uses the
single equality $\langle A,X\rangle=1$ and does not require $A$ to be
indefinite. The following level-boundedness result uses indefiniteness.

\begin{proposition}[Inner level boundedness for Family~$\mathrm{II}$]
\label{prop:familyIIprime_level_bounded}
Suppose that $\overline A$ is indefinite. Fix $\rho>0$ and any safeguarded
multiplier vector $\Upsilon$. For every $x^0$ in the inner domain, the
sublevel set $\mathcal S_\rho(x^0)$ in \eqref{eq:inner_sublevel} is bounded.
\end{proposition}
The proof is given in Appendix~\ref{app:familyIIprime-assumptions}.

\paragraph{A feasible reference point for Family $\mathrm{II}$.}
Choose a unit eigenvector $v$ of $A$ with $Av=\alpha v$ and $\alpha>0$, and set
\begin{equation}\label{eq:familyII_reference}
X^{\mathrm{ref}}=\frac{vv^\top}{\alpha},\qquad
C^{\mathrm{ref}}=S^{\mathrm{ref}}=0,\qquad
\widetilde y^{\mathrm{ref}}=0.
\end{equation}
Then $\langle A,X^{\mathrm{ref}}\rangle=1$, the matrix equality holds, and the complementarity conditions hold. Pairing these variables with $(Y^{\mathrm{ref}},Z^{\mathrm{ref}})=(X^{\mathrm{ref}},0)$ gives a feasible reference for \eqref{RE}, with
$f(w^{\mathrm{ref}})=\tfrac12\|\widehat C\|_F^2+\tfrac12\|vv^\top/\alpha-X_{\mathrm{obs}}\|_F^2$.
A positive eigenvalue suffices for this construction. Indefiniteness is still used to bound the remaining variables in the model-level convergence result below.

The model bounds and the reference-point comparison give the following conclusion for both families.

\begin{corollary}[Outer boundedness with feasible-reference comparison]
\label{cor:numerical_outer_boundedness}
Fix an instance of Family~$\mathrm{I}$ or Family~$\mathrm{II}$, with $\overline A$ indefinite in the latter case. Consider an infinite run of Algorithm~\ref{alg:sALcc} with Step~1 omitted, completed Step~3 solves, and $\nu_k\to0$.
Supply the corresponding feasible reference point in \eqref{eq:familyI_reference} or \eqref{eq:familyII_reference} and enforce \eqref{eq:reference_initialization} at every outer iteration.
Then the primal outer sequence is bounded and has a nonempty set of accumulation points. Every accumulation point is M-stationary for the corresponding SDCMPCC and its slack reformulation.
If Step~1 is retained with positive stopping tolerances, the algorithm terminates after finitely many outer iterations, provided every inner solve completes.
\end{corollary}
The proof is given in Appendix~\ref{app:outer-model-bounds}. A common AL upper bound yields the stated convergence conclusions, and Lemma~\ref{lem:reference_AL_bound} provides this bound under the feasible-reference comparison. For Family~$\mathrm{I}$, the stationarity conclusion concerns the lifting, with the scope described in Remark~\ref{rem:familyIprime_lifting_scope}. In all reported runs, $L_k(x^k)<L_k(x^{\mathrm{ref}})$ at every outer iteration, so $x^{k,0}=x^k$. Accepted inner updates satisfy $L_k(x^{k+1})\leq L_k(x^{k,0})$.

Propositions~\ref{prop:familyIprime_level_bounded} and~\ref{prop:familyIIprime_level_bounded} also ensure that each AL subproblem attains a global minimum: its initial sublevel set is nonempty and compact, and points outside it have larger objective values. At a minimizer $x_*=(w_*,u_*)$, the Fermat rule gives
\begin{equation*}
\nabla_wL_k(x_*)=0,\qquad
-\nabla_uL_k(x_*)\in\Nhat_\Omega(u_*)\subseteq\Nlim_\Omega(u_*).
\end{equation*}
Thus a point satisfying Step~3 exists for every positive target. Proposition~\ref{prop:inner_block_existence} and Corollary~\ref{cor:inner_finite_accuracy} further give finitely many accepted transitions under ideal block oracles and the stated parameter and error conditions. These existence statements do not guarantee that the practical routines return a candidate passing both acceptance tests within the allotted budget.

\paragraph{Composite-plate instances and exact relative noise.}
The test instance is the simply supported square plate under pure shear studied in \cite[Table~1 and Figs.~5--6]{defaria2018buckling}. The plate has side length $0.40\,\mathrm{m}$ and consists of a single $0.15\,\mathrm{mm}$ T300--5208 ply oriented at $+45^\circ$; its reference elastic moduli are $p_0=(E_1,E_2,G_{12})=(154.0,11.13,6.98)\,\mathrm{GPa}$ and its Poisson ratio is $\nu_{12}=0.304$. We discretize the plate with the first $m$ normalized sine modes in each coordinate direction, which gives a Ritz basis of dimension $n=m^2$, and we take $m\in\{4,5,7,10\}$ throughout. Let $K_m(p)$ denote the reduced elastic stiffness matrix at moduli $p$ and $B_m$ the reduced geometric stiffness matrix for a unit shear load that induces compression along the diagonal $x=y$.
The modal observation is extracted from an independent reference solution. On a $16\times16$ mesh of bicubic Hermite elements we compute the mode associated with the smallest positive buckling load at $p_0$, project the resulting displacement field onto the Ritz basis in $L^2$, and normalize the coefficient vector to obtain $q_m\in\mathbb R^n$ with $\|q_m\|_2=1$. Writing $g_m=q_m^\top B_mq_m>0$, we set 
\begin{equation*} 
X_{\mathrm{obs}}=q_mq_m^\top,\qquad A=\frac{B_m}{g_m},\qquad \overline A=\frac{B_m}{\|B_m\|_F},\qquad b_0=\frac{g_m}{\|B_m\|_F}. \end{equation*}
Then $\langle A,X_{\mathrm{obs}}\rangle=1$ holds by construction, and $A=\overline A/b_0$. Note that $\overline A$ is indefinite at every tested order. Since $m$ fixes the basis, $X_{\mathrm{obs}}$ and the constraint matrices are shared by all seeds and noise levels at that order.
The nominal stiffness is generated by perturbing the relative magnitudes of the three elastic moduli, since a common stiffness factor cancels under the normalization below. For each seed we draw $z\sim\mathcal N(0,I_3)$, subtract the mean of its entries, and rescale the result to a unit direction $d$; this direction is reused across all Ritz orders and noise levels of that seed. Setting $\widehat p_i=(p_0)_i\exp(td_i)$ and $\widehat K_m=K_m(\widehat p)$, we choose $t>0$ so that 
\begin{equation}
\label{eq:num_familyIIprime_noise} 
\frac{\|\widehat K_m-K_m(p_0)\|_F}{\|K_m(p_0)\|_F}=\eta_{\mathrm{rel}},\qquad \eta_{\mathrm{rel}}\in\{0.01,0.05\}. 
\end{equation} 
We set $\kappa_m=\|\widehat K_m\|_F$ and $\widehat C=\widehat K_m/\kappa_m$. The inverse problem fits the effective stiffness matrix in the Ritz basis, represented by $C$ in normalized coordinates and by $\kappa_m C$ in physical units.
The iteration is initialized from the nominal forward problem. Let $\widetilde y^0$ be the smallest positive eigenvalue of $\widehat C v=\widetilde y\,\overline A v$, and normalize its eigenvector $v_0$ by $v_0^\top A v_0=1$. In the working coordinate $T_n=nS$ this yields the starting point 
\begin{equation*} 
C^0=\widehat C,\qquad X^0=v_0v_0^\top,\qquad T_n^0=n\bigl(\widehat C-\widetilde y^0\overline A\bigr), 
\end{equation*}
which is feasible by construction and satisfies $T_n^0v_0=0$ for the nominal mode. We normalize the tangent stiffness by $\kappa_m/n$, the root mean square of the entries of $\widehat K_m$, so that at each fixed noise level the initial residual of the observed mode, $\|T_n^0q_m\|_2$, remains of the same order of magnitude across Ritz discretizations.
For each $(m,\eta_{\mathrm{rel}})$ group we run five seeds. The quadratic $w$-block is solved by matrix-free conjugate gradients with an absolute block-residual test; the inner and outer iteration limits are both $400$, and the penalty growth factor is $1.15$. Solver timing begins after the nominal forward initialization, and each run is limited to $7200$ seconds.

\begin{table}[htbp]
\centering
\caption{Finite-budget results for Family~$\mathrm{II}$ on the composite-plate instances. The Ritz order is $m$ and the matrix dimension is $n=m^2$. Residuals and solver times are medians over five runs per parameter group. $M/N$ denotes successful runs over planned runs.}
\label{tab:family-2-prime-results}
\normalsize
\setlength{\tabcolsep}{5pt}
\renewcommand{\arraystretch}{1.12}
\begin{tabular}{rrrrrrrrr}
\toprule
$m$ & $n$ & $\eta_{\mathrm{rel}}$ & $M/N$ & $R_{M,\Gamma}$ & $r_{\mathrm{stat}}$ & $r_{M,\Gamma}$ & $r_{\mathrm{pri}}$ & Time (s) \\
\midrule
4  & 16  & 0.01 & 5/5 & $8.62\!\times\!10^{-5}$ & $8.62\!\times\!10^{-5}$ & $7.15\!\times\!10^{-7}$ & $5.81\!\times\!10^{-5}$ & 2.65 \\
4  & 16  & 0.05 & 5/5 & $8.76\!\times\!10^{-5}$ & $8.76\!\times\!10^{-5}$ & $9.15\!\times\!10^{-7}$ & $5.62\!\times\!10^{-5}$ & 2.80 \\
5  & 25  & 0.01 & 5/5 & $8.84\!\times\!10^{-5}$ & $1.49\!\times\!10^{-6}$ & $8.84\!\times\!10^{-5}$ & $2.90\!\times\!10^{-5}$ & 3.74 \\
5  & 25  & 0.05 & 5/5 & $8.29\!\times\!10^{-5}$ & $2.24\!\times\!10^{-6}$ & $5.47\!\times\!10^{-5}$ & $4.58\!\times\!10^{-5}$ & 3.79 \\
7  & 49  & 0.01 & 5/5 & $9.55\!\times\!10^{-5}$ & $4.52\!\times\!10^{-5}$ & $1.32\!\times\!10^{-6}$ & $9.55\!\times\!10^{-5}$ & 14.91 \\
7  & 49  & 0.05 & 5/5 & $8.67\!\times\!10^{-5}$ & $4.52\!\times\!10^{-5}$ & $1.62\!\times\!10^{-6}$ & $8.46\!\times\!10^{-5}$ & 15.12 \\
10 & 100 & 0.01 & 5/5 & $7.66\!\times\!10^{-5}$ & $5.61\!\times\!10^{-5}$ & $6.74\!\times\!10^{-6}$ & $7.66\!\times\!10^{-5}$ & 96.84 \\
10 & 100 & 0.05 & 5/5 & $9.20\!\times\!10^{-5}$ & $6.97\!\times\!10^{-5}$ & $7.99\!\times\!10^{-6}$ & $9.20\!\times\!10^{-5}$ & 87.44 \\
\bottomrule
\end{tabular}
\end{table}

\subsection{Discussion of the results}

All $120$ planned runs satisfied the success criteria in Section~\ref{subsec:numerical_protocol} within the prescribed budgets, comprising $80$ runs for Family~$\mathrm{I}$ and $40$ for Family~$\mathrm{II}$ (Tables~\ref{tab:family-1-prime-results} and~\ref{tab:family-2-prime-results}). Across parameter groups, the median $R_{M,\Gamma}$ ranges from $5.25\times10^{-5}$ to $9.65\times10^{-5}$ for Family~$\mathrm{I}$ and from $7.66\times10^{-5}$ to $9.55\times10^{-5}$ for Family~$\mathrm{II}$.

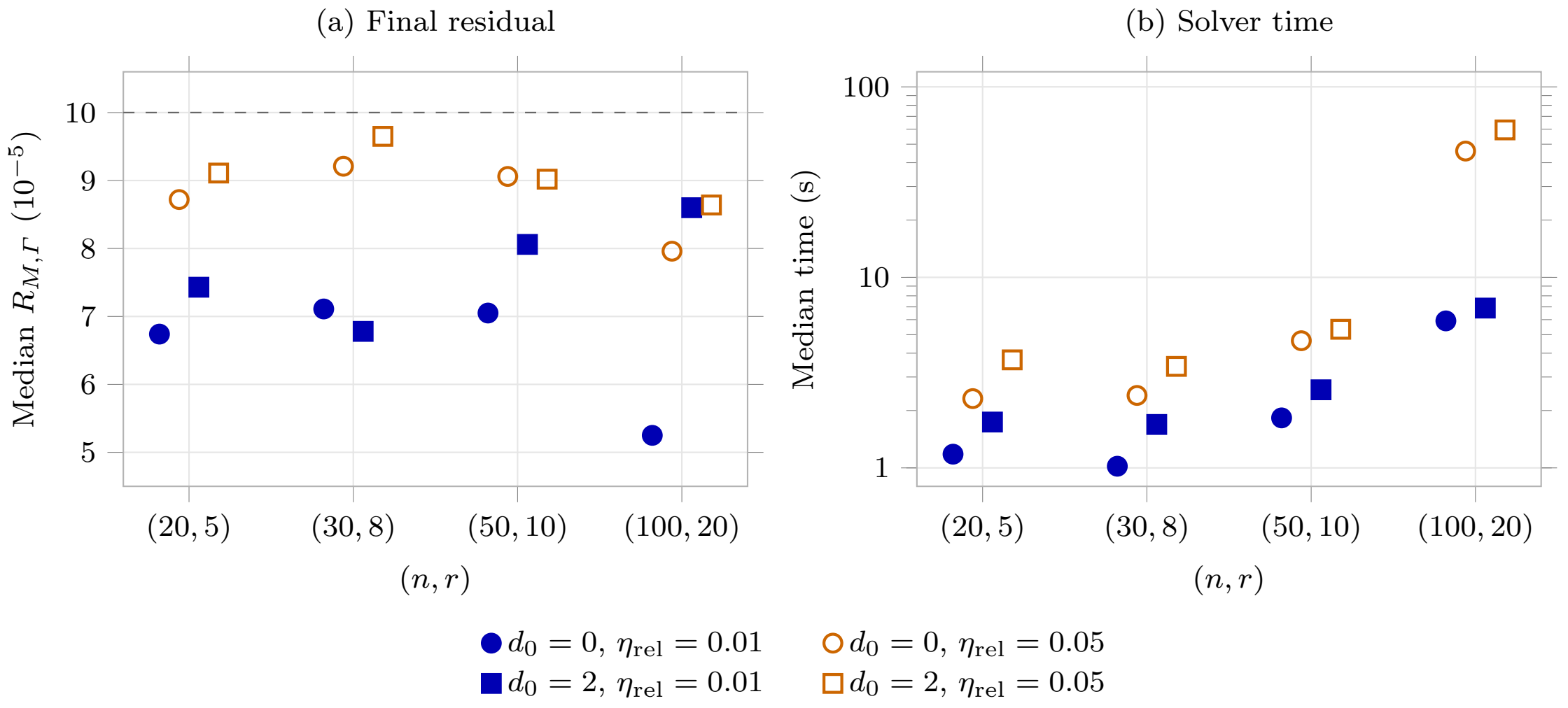
\begin{figure}[htbp]
  \centering
  \begin{tikzpicture}
\begin{groupplot}[
  group style={group size=2 by 1,horizontal sep=1.6cm},
  width=0.44\linewidth,height=5.5cm,
  xmin=-0.4,xmax=3.4,xtick={0,1,2,3},
  xticklabels={{$(20,5)$},{$(30,8)$},{$(50,10)$},{$(100,20)$}},
  xlabel={$(n,r)$},tick label style={font=\footnotesize},
  label style={font=\small},title style={font=\small},
  grid=major,grid style={gray!20},axis line style={gray!60},
  tick align=outside,legend cell align=left,
  legend style={draw=none,font=\footnotesize,at={(1.08,-0.31)},
    anchor=north,legend columns=2,/tikz/every even column/.append style={column sep=0.5cm}}
]
\nextgroupplot[title={(a) Final residual},
  ylabel={Median $R_{M,\Gamma}$ ($10^{-5}$)},
  ymin=4.5,ymax=10.6,ytick={5,6,7,8,9,10}]
\addplot[black!60,dashed,forget plot] coordinates {(-0.4,10) (3.4,10)};
\addplot[only marks,mark=*,mark size=2.4pt,line width=0.8pt,color=blue!70!black] coordinates {(-0.18,6.74) (0.82,7.11) (1.82,7.05) (2.82,5.25)};
\addlegendentry{$d_0=0$, $\eta_{\mathrm{rel}}=0.01$}
\addplot[only marks,mark=o,mark size=2.4pt,line width=0.8pt,color=orange!80!black] coordinates {(-0.06,8.72) (0.94,9.21) (1.94,9.06) (2.94,7.96)};
\addlegendentry{$d_0=0$, $\eta_{\mathrm{rel}}=0.05$}
\addplot[only marks,mark=square*,mark size=2.4pt,line width=0.8pt,color=blue!70!black] coordinates {(0.06,7.43) (1.06,6.78) (2.06,8.06) (3.06,8.60)};
\addlegendentry{$d_0=2$, $\eta_{\mathrm{rel}}=0.01$}
\addplot[only marks,mark=square,mark size=2.4pt,line width=0.8pt,color=orange!80!black] coordinates {(0.18,9.11) (1.18,9.65) (2.18,9.02) (3.18,8.64)};
\addlegendentry{$d_0=2$, $\eta_{\mathrm{rel}}=0.05$}
\nextgroupplot[title={(b) Solver time},ylabel={Median time (s)},
  ymode=log,ymin=0.8,ymax=120,ytick={1,10,100},
  yticklabels={$1$,$10$,$100$}]
\addplot[only marks,mark=*,mark size=2.4pt,line width=0.8pt,color=blue!70!black] coordinates {(-0.18,1.18) (0.82,1.02) (1.82,1.83) (2.82,5.91)};
\addplot[only marks,mark=o,mark size=2.4pt,line width=0.8pt,color=orange!80!black] coordinates {(-0.06,2.31) (0.94,2.40) (1.94,4.65) (2.94,46.02)};
\addplot[only marks,mark=square*,mark size=2.4pt,line width=0.8pt,color=blue!70!black] coordinates {(0.06,1.74) (1.06,1.69) (2.06,2.57) (3.06,6.90)};
\addplot[only marks,mark=square,mark size=2.4pt,line width=0.8pt,color=orange!80!black] coordinates {(0.18,3.69) (1.18,3.41) (2.18,5.34) (3.18,59.44)};
\end{groupplot}
\end{tikzpicture}
  \caption{Final residuals and solver times for Family~$\mathrm{I}$, using the group medians in Table~\ref{tab:family-1-prime-results}. Each marker summarizes five runs at the indicated $(n,r)$, noise level, and $d_0$. Residuals are displayed in units of $10^{-5}$; the dashed line corresponds to $R_{M,\Gamma}=10^{-4}$. The time axis is logarithmic.}
  \label{fig:family-i-endpoint-summary}
\end{figure}

\begin{figure}[htbp]
  \centering
  \begin{tikzpicture}
\begin{groupplot}[
  group style={group size=2 by 1,horizontal sep=1.6cm},
  width=0.44\linewidth,height=5.5cm,
  xmin=-0.4,xmax=3.4,xtick={0,1,2,3},
  xticklabels={{$(4,16)$},{$(5,25)$},{$(7,49)$},{$(10,100)$}},
  xlabel={$(m,n)$},tick label style={font=\footnotesize},
  label style={font=\small},title style={font=\small},
  grid=major,grid style={gray!20},axis line style={gray!60},
  tick align=outside,legend cell align=left,
  legend style={draw=none,font=\footnotesize,at={(1.08,-0.31)},
    anchor=north,legend columns=2,/tikz/every even column/.append style={column sep=0.5cm}}
]
\nextgroupplot[title={(a) Final residual},
  ylabel={Median $R_{M,\Gamma}$ ($10^{-5}$)},
  ymin=4.5,ymax=10.6,ytick={5,6,7,8,9,10}]
\addplot[black!60,dashed,forget plot] coordinates {(-0.4,10) (3.4,10)};
\addplot[only marks,mark=*,mark size=2.4pt,line width=0.8pt,color=blue!70!black] coordinates {(-0.07,8.62) (0.93,8.84) (1.93,9.55) (2.93,7.66)};
\addlegendentry{$\eta_{\mathrm{rel}}=0.01$}
\addplot[only marks,mark=o,mark size=2.4pt,line width=0.8pt,color=orange!80!black] coordinates {(0.07,8.76) (1.07,8.29) (2.07,8.67) (3.07,9.20)};
\addlegendentry{$\eta_{\mathrm{rel}}=0.05$}
\nextgroupplot[title={(b) Solver time},ylabel={Median time (s)},
  ymode=log,ymin=0.8,ymax=120,ytick={1,10,100},
  yticklabels={$1$,$10$,$100$}]
\addplot[only marks,mark=*,mark size=2.4pt,line width=0.8pt,color=blue!70!black] coordinates {(-0.07,2.65) (0.93,3.74) (1.93,14.91) (2.93,96.84)};
\addplot[only marks,mark=o,mark size=2.4pt,line width=0.8pt,color=orange!80!black] coordinates {(0.07,2.80) (1.07,3.79) (2.07,15.12) (3.07,87.44)};
\end{groupplot}
\end{tikzpicture}
  \caption{Final residuals and solver times for Family~$\mathrm{II}$, using the group medians in Table~\ref{tab:family-2-prime-results}. Each marker summarizes five runs at the indicated Ritz order $m$, matrix dimension $n=m^2$, and noise level. The residual scale, reporting threshold, and logarithmic time axis are as in Figure~\ref{fig:family-i-endpoint-summary}.}
  \label{fig:family-ii-endpoint-summary}
\end{figure}
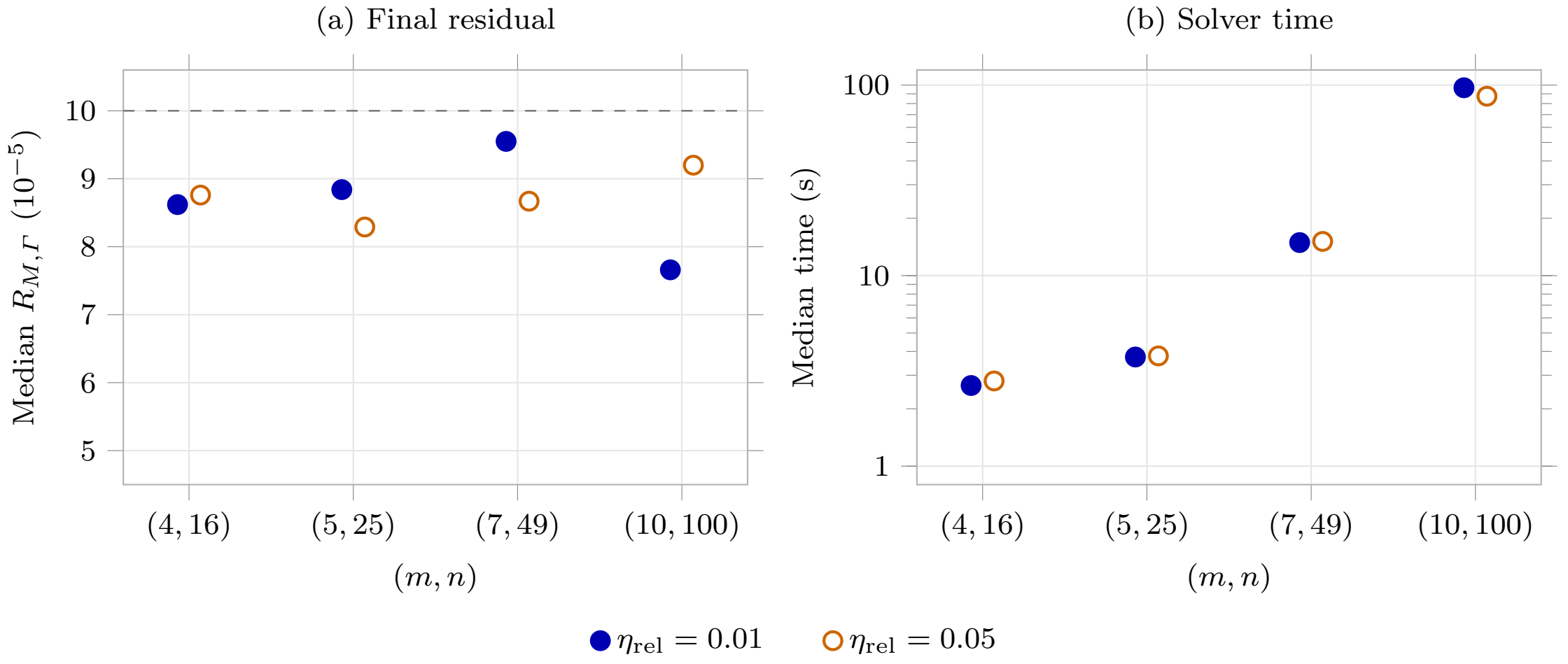

Figures~\ref{fig:family-i-endpoint-summary} and~\ref{fig:family-ii-endpoint-summary} summarize the final residuals and solver times. Median solver times range from $1.02$ to $59.44$ seconds for Family~$\mathrm{I}$ and from $2.65$ to $96.84$ seconds for Family~$\mathrm{II}$. For Family~$\mathrm{I}$, the median time is larger at $\eta_{\mathrm{rel}}=0.05$ than at $\eta_{\mathrm{rel}}=0.01$ for every tested $(n,r,d_0)$. For Family~$\mathrm{II}$, the median time increases with the Ritz order at both noise levels.

\FloatBarrier
\section{Conclusion}
\label{sec:conclusion}

SPALM applies the safeguarded augmented Lagrangian framework of Jia et al.~\cite{JiaKanzowMehlitzWachsmuth2022} to SDCMPCCs with a subproblem accuracy target defined by the limiting normal cone to $\Omega$. This target includes coupled first-order normals that a scalar orthogonality representation can omit even at strictly complementary points. Feasible slack candidates enforce complementarity, while smooth linking residuals control the distance of the matrix image to $\Omega$. For infinite outer sequences with vanishing absolute subproblem errors, feasible accumulation points are AM-stationary and, under AM-regularity, M-stationary. Active-set SDCMPCC-LICQ is sufficient for this regularity. When solving an augmented Lagrangian subproblem, the proximal alternating framework gives convergence of the full residual to zero along infinite accepted sequences under the stated level-boundedness, proximal-parameter, and error conditions. Square-summable block errors yield an $O(T^{-1/2})$ best-residual bound in accepted iterations.

Huber--Stiefel and fixed-basis procedures generate feasible slack candidates, whose acceptance requires descent tests and limiting-normal certificates. Spectral bounds quantify the objective gap of the fixed-basis candidates, but finite-budget candidate generation need not produce an accepted update. The correlation-matrix lifting has a compact feasible set and satisfies active-set SDCMPCC-LICQ at every feasible point. It preserves the projected feasible set and global optimal value.

For the two model families, explicit feasible references and the comparison rule \eqref{eq:reference_initialization} give a uniform outer AL upper bound. Under the stated model assumptions, completed inner solves, and vanishing absolute errors, the resulting outer sequence is bounded and all its accumulation points are M-stationary.

The finite-budget experiments report successful termination in all $120$ runs for the lifting and a single-equality observation-fitting inverse SDP. Success combines absolute feasibility and stationarity tests with a normalized residual evaluated using a retained normal witness. These floating-point checks provide numerical evidence for the tested implementation. Further work should quantify the work required to generate accepted updates and develop validated error bounds for the normal-witness calculations.

\appendix

\section{Complete spectral characterization of
\texorpdfstring{$\mathcal K_M$}{KM}}
\label{app:KM-auxiliary}

This appendix gives the full spectral characterization of $\mathcal K_M(u)=\Nlim_\Omega(u)$ from \cite[Proposition~3.3 and Theorem~3.1]{DingSunYe2014}. Fix $u=(Y,Z)\in\Omega$ and use the joint decomposition, spectral coordinates, and common multiplier space $\mathcal K_W(u)$ introduced in Subsection~\ref{subsec:spectral_geometry}.

Let $\beta$ be the zero spectral index set associated with $u=(Y,Z)\in\Omega$. Within the zero block, identify $\beta$ with the local indices $\{1,\ldots,|\beta|\}$. Let $\mathcal P(\beta)$ consist of ordered triples of pairwise disjoint subsets whose union is $\beta$, with empty subsets allowed. Write
\begin{equation*}
\mathbb R_{\succeq}^{|\beta|}:=
\left\{
a\in\mathbb R^{|\beta|}:
a_1\geq\cdots\geq a_{|\beta|}
\right\}.
\end{equation*}

For $a\in\mathbb R_{\succeq}^{|\beta|}$, let $D(a)$ denote the generalized first divided-difference matrix of $t\mapsto\max\{t,0\}$ at $a$, namely,
\begin{equation*}
(D(a))_{ij}
=
\begin{cases}
\displaystyle
\frac{\max\{a_i,0\}-\max\{a_j,0\}}{a_i-a_j},
& a_i\neq a_j,\\[2mm]
1,
& a_i=a_j>0,\\
0,
& a_i=a_j\leq0,
\end{cases}
\qquad
i,j=1,\ldots,|\beta|.
\end{equation*}

Define
\begin{equation*}
\mathbb U_{|\beta|}
:=
\left\{
\overline\Xi\in\mathbb S^{|\beta|}:
\overline\Xi=\lim_{k\to\infty}D(a^k),\ a^k\to0,\ a^k\in\mathbb R_{\succeq}^{|\beta|}
\right\}.
\end{equation*}

For any $\Xi_1\in\mathbb U_{|\beta|}$, choose a sequence $a^k$ realizing its divided-difference limit. After passing to a subsequence, the sign sets
\begin{equation*}
 \beta_+:=\{i:a_i^k>0\},\qquad
 \beta_0:=\{i:a_i^k=0\},\qquad
 \beta_-:=\{i:a_i^k<0\}
\end{equation*}
are independent of $k$. We call the resulting partition associated with $\Xi_1$; the same matrix $\Xi_1$ may admit more than one associated partition. Since the vectors $a^k$ are ordered, these sets occur in positive, zero, and negative order. Thus there is a partition
\begin{equation*}
\pi(\beta):=(\beta_+,\beta_0,\beta_-)\in\mathcal P(\beta)\end{equation*}
such that
\begin{equation*}\Xi_1=
\begin{bmatrix}
E_{\beta_+\beta_+} & E_{\beta_+\beta_0} & (\Xi_1)_{\beta_+\beta_-}\\
E_{\beta_+\beta_0}^{\top} & 0 & 0\\
(\Xi_1)_{\beta_+\beta_-}^{\top} & 0 & 0
\end{bmatrix},
\end{equation*}
where every entry of $(\Xi_1)_{\beta_+\beta_-}$ belongs to $[0,1]$.

For this choice of the realizing sequence and partition, define
\begin{equation*}
\Xi_2:=
\begin{bmatrix}
0 & 0 & E_{\beta_+\beta_-}-(\Xi_1)_{\beta_+\beta_-} \\
0 & 0 & E_{\beta_0\beta_-} \\
\bigl(E_{\beta_+\beta_-}-(\Xi_1)_{\beta_+\beta_-}\bigr)^{\top} &
E_{\beta_0\beta_-}^{\top} & E_{\beta_-\beta_-}
\end{bmatrix}.
\end{equation*}

For $Q\in\mathbb O^{|\beta|}$, $Q_{\beta_0}$ denotes the submatrix of $Q$ formed by the columns corresponding to the index set $\beta_0$.

With these objects, the limiting normal cone has the following spectral representation:
\begin{equation*} \begin{aligned} \mathcal K_M(u) = \Bigl\{ (\Gamma^Y,\Gamma^Z)\in\mathcal K_W(u):\;
& \text{there exist }Q\in\mathbb O^{|\beta|} \text{ and }\Xi_1\in\mathbb U_{|\beta|}\text{ with an associated partition }\\
& \pi(\beta)=(\beta_+,\beta_0,\beta_-) \text{ and matrix }\Xi_2 \text{ such that}\\
& \Xi_1\circ \bigl(Q^{\top}\widetilde\Gamma^Y_{\beta\beta}Q\bigr) + \Xi_2\circ \bigl(Q^{\top}\widetilde\Gamma^Z_{\beta\beta}Q\bigr) =0,\\
& Q_{\beta_0}^{\top} \widetilde\Gamma^Y_{\beta\beta} Q_{\beta_0} \preceq0,\qquad Q_{\beta_0}^{\top} \widetilde\Gamma^Z_{\beta\beta} Q_{\beta_0} \succeq0 \Bigr\}. \end{aligned} \end{equation*}

All conditions involving empty index sets are vacuous. If $\beta=\varnothing$, the degenerate-block conditions disappear and $\mathcal K_M(u)=\mathcal K_W(u)$. If $Y+Z=0$, then $Y=Z=0$ and $\mathcal K_W(u)=\Sn\times\Sn$; the formula reduces to the limiting normal cone at the origin. The union in this characterization ranges over all such choices of $Q$, $\Xi_1$, and associated partition. Both the divided-difference limit and its associated sign partition are required; entrywise bounds on the mixed block alone do not characterize the admissible choices.

\section{Proof of the fixed-basis objective-gap bound}
\label{app:fixed-basis-proof}

\begin{proof}[Proof of Proposition~\ref{prop:fixedA_gap}]
Set $\mathcal D:=\begin{pmatrix}A_+^2&0\\0&B_-^2\end{pmatrix}\in\mathbb S^{2n}$, with $A_+,B_-$ the spectral positive and negative parts. For each column $p_i$ of $P\in\On$, choose $q_i=(p_i,0)$ if $(p_i^\top Ap_i)_+\geq(p_i^\top Bp_i)_-$, and $q_i=(0,p_i)$ otherwise. The vectors $q_i$ are orthonormal. Jensen's inequality for $t\mapsto(t_+)^2$, applied in an eigenbasis of $A$ or $-B$, gives
\begin{equation*}
F(P)\leq\sum_{i=1}^nq_i^\top\mathcal Dq_i
\leq\sum_{j=1}^n\lambda_j^\downarrow(\mathcal D)=U_{\mathrm{spec}}.
\end{equation*}
For the second inequality, write the rank-$n$ orthogonal projector $\sum_iq_iq_i^\top$ in an eigenbasis of $\mathcal D$: its diagonal entries lie in $[0,1]$ and sum to $n$, so their weighted sum is at most the sum of the $n$ largest eigenvalues. Equation~\eqref{eq:minvalP} now gives $L_{\mathrm{spec}}\leq v_\star\leq V_A$.

To obtain \eqref{eq:optimal_spectral_pairing}, select $n$ entries attaining $U_{\mathrm{spec}}$ from the two spectral lists. If $k$ come from the $A$-list, pair these with the $k$ unselected $B$ entries, and pair the other $A$ entries with the selected $B$ entries. Each selected entry is at least its partner, so the sum of maxima is $U_{\mathrm{spec}}$. Ties are arbitrary but the orderings of $P_B$ and $\lambda(B)$ must agree.

Finally, $\diag(P_A^\top BP_A)=T\lambda(B)$. On $[-R_0,R_0]$, the functions $\phi_i(s):=\max\{(\lambda_i(A))_+^2,(s_-)^2\}$ are $2R_0$-Lipschitz. Hence
\begin{equation}\label{eq:fixedA_gap_identity}
\begin{aligned}
V_A-L_{\mathrm{spec}}
&=\sum_i\bigl[\phi_i(\lambda_i(B))-\phi_i((T\lambda(B))_i)\bigr]\\
&\leq2R_0\|(I-T)\lambda(B)\|_1
\leq2\sqrt n\,R_0\|I-T\|_2\|B\|_F.
\end{aligned}
\end{equation}
Together with $L_{\mathrm{spec}}\leq v_\star\leq V_A$, this proves the claim.
\end{proof}

\section{Verification of the assumptions for the numerical models}
\label{app:numerical-assumptions}

\subsection{Family \texorpdfstring{$\mathrm{I}$}{I}}
\label{app:familyIprime-assumptions}

\begin{proof}[Proof of Proposition~\ref{prop:familyIprime_LICQ}]
Let $((\bar X,\bar S),(\bar Y,\bar Z))\in\widetilde{\mathcal F}$ be a
feasible point of the induced slack reformulation. The linking equalities give
$\bar Y=\bar X$ and $\bar Z=-\bar S$. By
Proposition~\ref{prop:reduced_LICQ}, consider a scalar $\lambda$, a vector
$\mu\in\mathbb R^n$, and matrices $\Gamma^Y,\Gamma^Z\in\Sn$ satisfying the
homogeneous reduced system. With $D:=\Diag(\mu)$, its two stationarity
equations and spectral relation are
\begin{equation*}
 D+\Gamma^Y=0,\qquad
 2\lambda(\bar S-I)-\Gamma^Z=0,\qquad
 \bar X\Gamma^Y+\Gamma^Z\bar S=0,
\end{equation*}
where $\lambda=0$ whenever the scalar inequality is inactive and $\lambda$ is unrestricted when it is active. Eliminating
$\Gamma^Y$ and $\Gamma^Z$ gives
\begin{equation*}
 -\bar XD+2\lambda(\bar S-I)\bar S=0.
\end{equation*}
Since $\bar X\bar S=0$, the two terms are orthogonal because
\begin{equation*}
 \ip{\bar XD}{(\bar S-I)\bar S}
 =\tr\!\bigl(D\bar X(\bar S-I)\bar S\bigr)=0.
\end{equation*}
Hence
\begin{equation*}
 \norm{\bar XD}_F^2
 +4\lambda^2\norm{(\bar S-I)\bar S}_F^2=0.
\end{equation*}
Thus $\bar XD=0$; since $\diag(\bar X)=\mathbf1$,
$(\bar XD)_{ii}=\bar X_{ii}\mu_i=\mu_i$, and hence $\mu=0$. If the
inequality is inactive, then $\lambda=0$. If it is active and
$\lambda\neq0$, then $(\bar S-I)\bar S=0$, so $\bar S$ is an orthogonal
projector. Consequently $\norm{I-\bar S}_F^2$ is an integer, contradicting
$\norm{I-\bar S}_F^2=r+\vartheta \in (r, r+1)$. Therefore $\lambda=0$, and the stationarity
equations yield $\Gamma^Y=\Gamma^Z=0$.
\end{proof}

\begin{proof}[Proof of Corollary~\ref{cor:familyI_no_singularM}]
Proposition~\ref{prop:familyIprime_LICQ} excludes nonzero homogeneous tuples in the full linearized multiplier space, and hence in the smaller set of actual normal-cone multipliers.
\end{proof}

\begin{proof}[Proof of Proposition~\ref{prop:familyIprime_level_bounded}]
Fix $\rho>0$ and the multipliers. Let $c:=L_\rho(x^0)$ and take any
$x=(X,S,Y,Z)\in\mathcal S_\rho(x^0)$. All terms in \eqref{eq:AL} are
nonnegative. The objective therefore bounds $X$, while
\begin{equation*}
 \left(g(X,S)+\frac{\lambda}{\rho}\right)_+
 \leq \sqrt{\frac{2c}{\rho}}
\end{equation*}
implies
\begin{equation*}
 \|I-S\|_F^2
 \leq r+\vartheta+\sqrt{\frac{2c}{\rho}}-\frac{\lambda}{\rho}.
\end{equation*}
Thus $S$ is bounded. Keeping the multiplier shifts in the linking terms gives
\begin{equation*}
 \|X-Y\|_F\leq\sqrt{\frac{2c}{\rho}}+\frac{\|M^G\|_F}{\rho},
 \qquad
 \|S+Z\|_F\leq\sqrt{\frac{2c}{\rho}}+\frac{\|M^H\|_F}{\rho}.
\end{equation*}
Hence $Y$ and $Z$ are bounded as well. These bounds depend on
$\rho$ and the multipliers.
\end{proof}

\subsection{Family \texorpdfstring{$\mathrm{II}$}{II}}
\label{app:familyIIprime-assumptions}

\begin{proof}[Proof of Proposition~\ref{prop:familyIIprime_LICQ}]
Let $(\bar w,(\bar Y,\bar Z))\in\widetilde{\mathcal F}$ be feasible for the
induced slack reformulation, where
$\bar w=(\bar C,\bar X,\bar S,\widetilde y_0)$. Thus $\bar Y=\bar X$ and
$\bar Z=-\bar S$. In the reduced system of
Proposition~\ref{prop:reduced_LICQ}, write the equality multiplier as
$(\sigma,T)\in\mathbb R\times\Sn$. The stationarity equations with respect
to $C,X,S,\widetilde y$, respectively, are
\begin{equation*}
 T=0,\qquad
 \sigma A+\Gamma^Y=0,\qquad
 -T-\Gamma^Z=0,\qquad
 -\langle\overline A,T\rangle=0.
\end{equation*}
Thus $T=\Gamma^Z=0$ and $\Gamma^Y=-\sigma A$. The spectral relation gives
$\sigma\bar XA=0$. Taking the trace and using
$\langle A,\bar X\rangle=1$ yields $\sigma=0$. Hence all equality and
spectral multipliers vanish.
\end{proof}

\begin{proof}[Proof of Proposition~\ref{prop:familyIIprime_level_bounded}]
Fix $\rho>0$ and the multipliers, and write the equality multiplier as $\mu=(\sigma,T)$. Set $c:=L_\rho(x^0)$ and $\delta:=\sqrt{2c/\rho}$. For any point in $\mathcal S_\rho(x^0)$, the nonnegative objective terms in \eqref{eq:AL} bound $C$ and $X$. Set $Q:=S+Z$ and $E:=C-\widetilde y\,\overline A-S$. The shifted linking and equality penalties give
\begin{equation*}
 \begin{aligned}
 \|X-Y\|_F&\leq\delta+\frac{\|M^G\|_F}{\rho},\qquad
 \|Q\|_F\leq\delta+\frac{\|M^H\|_F}{\rho},\\
 \|E\|_F&\leq\delta+\frac{\|T\|_F}{\rho}.
 \end{aligned}
\end{equation*}
Thus $Y$, $Q$, and $E$ are bounded throughout the sublevel set. Since $P:=-Z\succeq0$ and $S=Q+P$,
\begin{equation*}
 R:=\widetilde y\,\overline A+P=C-Q-E
\end{equation*}
is bounded. Let $v_+$ and $v_-$ be unit eigenvectors of $\overline A$
associated with eigenvalues $\alpha_+>0$ and $\alpha_-<0$. Since $P\succeq0$,
\begin{equation*}
 \widetilde y\alpha_+\leq v_+^\top Rv_+,
 \qquad
 \widetilde y\alpha_-\leq v_-^\top Rv_-.
\end{equation*}
Choose $M$ so that $\norm{R}_2\leq M$ throughout the sublevel set. Then
\begin{equation*}
 -\frac{M}{|\alpha_-|}\leq\widetilde y\leq\frac{M}{\alpha_+}.
\end{equation*}
Consequently $P=R-\widetilde y\,\overline A$ is bounded, and hence so are
$Z=-P$ and $S=Q+P$.
\end{proof}

These sublevel-set bounds depend on both signs of the spectrum of $\overline A$ and need not be uniform across instances or outer iterations. The indefiniteness assumption cannot be
discarded from this result. For example, if $0\ne\overline A\succeq0$,
the inner-domain ray
\begin{equation*}
 C=X=Y=0,\qquad \widetilde y=-t,\qquad
 S=t\overline A,\qquad Z=-t\overline A,\qquad t\geq0,
\end{equation*}
has $(Y,Z)\in\Omega$, zero linking and matrix-equality residuals, and
scalar-equality residual $-1$. For fixed multipliers and $\rho$, the
augmented Lagrangian is constant along this unbounded ray. Hence any sublevel set whose reference point lies on the ray is unbounded.

\subsection{Uniform bounds across outer iterations}
\label{app:outer-model-bounds}
\begin{proof}[Proof of Corollary~\ref{cor:numerical_outer_boundedness}]
Lemma~\ref{lem:reference_AL_bound} gives $L_k(x^{k+1})\leq\bar L:=\bar L_{\mathrm{ref}}$ for every $k$. Both objectives are nonnegative, so $\bar L\geq0$. Using $\rho_k\geq\rho_0>0$ and $\|\bar\Upsilon^k\|\leq M_{\mathrm{saf}}$ in the squared-distance form of $L_k$ gives
\begin{equation}\label{eq:model_uniform_residual_bound}
\begin{aligned}
\operatorname{dist}(\Psi(x^{k+1}),C)
&\leq\operatorname{dist}\left(\Psi(x^{k+1})+\frac{\bar\Upsilon^k}{\rho_k},C\right)
+\frac{M_{\mathrm{saf}}}{\rho_k}\\
&\leq\sqrt{\frac{2\bar L}{\rho_k}}+\frac{M_{\mathrm{saf}}}{\rho_k}
\leq D:=\sqrt{\frac{2\bar L}{\rho_0}}+\frac{M_{\mathrm{saf}}}{\rho_0}.
\end{aligned}
\end{equation}
Thus every equality residual, linking residual and inequality positive part is uniformly bounded by $D$, without an upper bound on $\rho_k$.

For Family~$\mathrm{I}$, the objective and the inequality positive part yield
\begin{equation*}
\|X^{k+1}\|_F\leq\|C_{\mathrm{obs}}\|_F+\sqrt{2\bar L},\qquad
\|I-S^{k+1}\|_F^2\leq r+\vartheta+D.
\end{equation*}
Hence $X,S$ are bounded, and the linking residuals bound $Y,Z$. These estimates apply to the possibly infeasible outer points.

For Family~$\mathrm{II}$, the objective bounds $C,X$. Suppress the superscript $k+1$ and set $E=C-\widetilde y\overline A-S$, $Q=S+Z$ and $P=-Z\succeq0$. Then $\|E\|_F,\|Q\|_F\leq D$, and
\begin{equation*}
R:=C-Q-E=\widetilde y\overline A+P,\qquad
\|R\|_2\leq M:=\|\widehat C\|_F+\sqrt{2\bar L}+2D.
\end{equation*}
For fixed unit eigenvectors $v_+,v_-$ of $\overline A$ with eigenvalues $\alpha_+>0>\alpha_-$, positive semidefiniteness of $P$ gives
\begin{equation*}
\widetilde y\alpha_\pm\leq v_\pm^\top Rv_\pm,
\qquad -\frac{M}{|\alpha_-|}\leq\widetilde y\leq\frac{M}{\alpha_+}.
\end{equation*}
Consequently $P$, $Z=-P$, $S=Q+P$ and $Y$ are bounded. The constants depend on the fixed instance and need not be uniform as its spectral gaps or $b_0$ vary.

Boundedness gives an accumulation point. Proposition~\ref{prop:feasibility}(ii) makes every accumulation point feasible. Theorem~\ref{thm:outer_AM_stationarity}, the model LICQ results, and Proposition~\ref{thm:LICQ_AMregular} then give M-stationarity, with Proposition~\ref{prop:Mstat_equiv} transferring it through the unique slack lift. With positive stopping tolerances, any infinite run would have a feasible accumulation point and hence would satisfy Step~1 at a finite iteration by the stopping argument following Theorem~\ref{thm:outer_AM_stationarity}.
\end{proof}

\section{NP-hardness of exact slack projection}\label{app:projection-hardness}

\begin{proof} We establish the hardness claim in Remark~\ref{rem:projection-hardness} by reducing Simple Max-Cut \cite{GareyJohnsonStockmeyer1976}. Let $G$ be a simple unweighted graph on $n\geq2$ vertices with adjacency entries $g_{ij}$, and denote its maximum cut size by $\operatorname{mc}(G)$. Set \[ \begin{gathered} d_i=i^2,\qquad D=\Diag(d_1,\ldots,d_n),\qquad C_{ij}=g_{ij}\frac{|d_i-d_j|}{ij},\\ \varepsilon=(1024n^{16})^{-1},\qquad A=D+\varepsilon C,\qquad B=-D+\varepsilon C. \end{gathered} \] Throughout, $\|\cdot\|=\|\cdot\|_F$ and $[U,V]=UV-VU$. Put $b=\|C\|$, $H=DC+CD$, and $M=\|H\|$. Then $b\leq n^2$, $M\leq2n^4$, and $\varepsilon b<1$, so $A\succ0$ and $B\prec0$.

Let $\mathcal P=\{R\in\mathbb S^n:R^2=R\}$. Every complementary pair has orthogonal supports. For fixed $R\in\mathcal P$, the minimizing supported pair is $Y=RAR$, $Z=(I-R)B(I-R)$. Writing $v(A,B)$ for the value of \eqref{eq:projAB}, expansion therefore gives \[ \begin{aligned} v(A,B)&=c_0-\max_{R\in\mathcal P}h(R),\qquad c_0=\|D\|^2+\varepsilon^2b^2,\\ h(R)&=-\|[D,R]\|^2+2\varepsilon\langle H,R\rangle-\varepsilon^2\|[C,R]\|^2. \end{aligned} \]

For a cut $S$, set $R_S=\Diag(\mathbf1_S)$ and $Q_S=I-R_S$. On symmetric matrices $T=R_STQ_S+Q_STR_S$, define \[ q_S(T)=-\|[D,T]\|^2+2\varepsilon\langle H,T\rangle-\varepsilon^2\|[C,R_S]\|^2. \] Let $T_S$ have cross-block entries $(T_S)_{ij}=\varepsilon H_{ij}/(d_i-d_j)^2$ for $i\in S$, $j\notin S$, with symmetric counterparts and zero same-side blocks. Completing squares yields \[ q_S(T)=8\varepsilon^2\operatorname{cut}_G(S)-\|[D,T-T_S]\|^2,\qquad \|T_S\|\leq\varepsilon M. \]

Take a global maximizer $R$ of $h$, which exists by compactness of $\mathcal P$. Since $H$ has zero diagonal and $|d_i-d_j|\geq1$ for $i\ne j$, \[ 0\leq h(R)\leq-t^2+2\varepsilon Mt,\qquad t:=\|\operatorname{offdiag}R\|\leq2\varepsilon M\leq\tfrac12. \] Round the diagonal of $R$ to obtain $R_S$ and write $\Delta=R-R_S=T+W$, where $T$ is the cross-block part. The identities $R_{ii}(1-R_{ii})=\sum_{j\ne i}R_{ij}^2$ and $R_S\Delta+\Delta R_S+\Delta^2=\Delta$ imply $\|\Delta\|\leq2t$ and $\|W\|\leq\|\Delta\|^2\leq4t^2$. The matrices $[D,T]$ and $[D,W]$ have disjoint block supports. Using $\|[C,R_S]\|\leq2b$ and $\|[C,\Delta]\|\leq4bt$, expansion gives \[ \begin{aligned} h(R)&\leq q_S(T)+8\varepsilon Mt^2+24\varepsilon^2b^2t\\ &\leq8\varepsilon^2\operatorname{mc}(G)+32\varepsilon^3M^3+48\varepsilon^3b^2M\\ &\leq8\varepsilon^2\operatorname{mc}(G)+352\varepsilon^3n^{12}. \end{aligned} \]

For the reverse bound, fix any cut $S$ and set $K=T_SR_S-R_ST_S$ and $\widetilde R=e^KR_Se^{-K}$. Then $K^\top=-K$, $[K,R_S]=T_S$, and $q:=\|K\|=\|T_S\|\leq\varepsilon M\leq\tfrac14$. The first and second derivatives of $e^{sK}R_Se^{-sK}$ have norms at most $q$ and $2q^2$. Thus, for $\widetilde\Delta=\widetilde R-R_S$, \[ \|\widetilde\Delta\|\leq q,\qquad \|\widetilde\Delta-T_S\|\leq q^2. \] Expanding $h$ and using these bounds yields \[ \begin{aligned} h(\widetilde R)&\geq q_S(T_S)-9n^4q^3-2\varepsilon Mq^2-9\varepsilon^2b^2q\\ &\geq8\varepsilon^2\operatorname{cut}_G(S)-106\varepsilon^3n^{16}. \end{aligned} \]

Taking a maximum cut and combining both bounds gives \[ \left|v(A,B)-c_0+8\varepsilon^2\operatorname{mc}(G)\right|\leq512\varepsilon^3n^{16}=\frac{\varepsilon^2}{2}. \] Hence, for every integer $k$, \[ \operatorname{mc}(G)\geq k \quad\Longleftrightarrow\quad v(A,B)\leq c_0-8\varepsilon^2\left(k-\tfrac12\right). \] The matrices and threshold are rational, have polynomial encoding length, and can be constructed in polynomial time. This completes the reduction. \end{proof}

\bibliographystyle{spmpsci}
\bibliography{references}

\end{document}